%% file: main.tex
\documentclass[review]{siamart0516}
\input{shared_info}

\ifpdf
\hypersetup{
  pdftitle={\TheTitle},
  pdfauthor={\TheAuthors}
}
\fi

\begin{document}\sloppy

\maketitle

\begin{abstract}
In this paper, we propose an adaptive fast solver for a general class of symmetric positive definite (SPD) matrices which include the well-known graph Laplacian. We achieve this by developing an adaptive operator compression scheme and a multiresolution matrix factorization algorithm which achieve nearly optimal performance on both complexity and well-posedness. To develop our adaptive operator compression and multiresolution matrix factorization methods, we first introduce a novel notion of {\bf energy decomposition} for SPD matrix $A$ using the representation of energy elements. The interaction between these energy elements depicts the underlying topological structure of the operator. This concept of decomposition naturally reflects the hidden geometric structure of the operator which inherits the localities of the structure. By utilizing the intrinsic geometric information under this Energy framework, we propose a systematic operator compression scheme for the inverse operator $A^{-1}$. In particular, with an {\it appropriate partition} of the underlying geometric structure, we can construct localized basis by using the concept of {\bf interior} and {\bf closed energy}. Meanwhile, two important localized quantities are introduced, namely the {\bf error factor} and the {\bf condition factor}. Our error analysis results show that these two factors will be the guidelines for finding the appropriate partition of the basis functions such that prescribed compression error and acceptable condition number can be achieved. By virtue of this insight, we propose the {\bf Patch Pairing} algorithm to realize our energy partition framework for operator compression with controllable compression error and condition number.
\end{abstract}



\section{Introduction}
\label{sec:introduction}
Fast algorithms for solving symmetric positive definite (SPD) linear systems have found broad applications across both theories and practice, including machine learning \cite{belkin2006manifold, cai2011graph, scholkopf2002learning}, computer vision \cite{belkin2001laplacian, yan2007graph, bookstein1989principal}, image processing \cite{ashburner2007fast, cao2005large, myronenko2010point}, computational biology \cite{colovos1993verification, newman2006modularity}, etc. For instance, the graph laplacian, which has deep connection between the combinatorial properties of the graph $G$ and the linear algebraic properties of the laplacian $L$, is one of the foundational problems in data analysis. Performing finite element simulation of wide range of physical systems will also introduce the corresponding stiffness matrix, which is also symmetric and positive definite. 

The philosophy of this work is inspired by the Spectral Graph Theory \cite{ng2001spectral, chung1997spectral} in which the spectrum and the geometry of graphs are highly correlated. By computing the spectrum of the graph, the intrinsic geometric information can be directly obtained and various applications can be found \cite{boykov2003computing, kolluri2004spectral}. However, the price for finding the spectrum of graph is relatively expensive as it involves solving global eigen problems. On the contrary, the Algebraic Multigrid method is a purely matrix-based multiresolution type solver. It simply uses the interaction of nonzero entries within the matrix as an indicator to describe the geometry implicitly. These stimulating techniques and concepts motivate us to look into the problems from two different points of view and search for a brand new framework which can integrate the advantages from both ends. 

In this paper, we propose an adaptive fast solver for a general class of symmetric positive definite (SPD) matrices. We achieve this by developing an adaptive operator compression scheme and a multiresolution matrix factorization algorithm both with nearly optimal performance on complexity and well-posedness. These methods are developed based on a newly introduced framework, namely, the {\bf energy decomposition} for SPD matrix $A$ to extract its hidden geometric information. For the ease of discussion, we first consider $A = L$, the graph laplacian of an undirected graph $G$. Under this framework, we reformulate the connectivity of subgraphs in $G$ as the interaction between {\bf energies}. These interactions reveal the intrinsic geometric information hidden in $L$. In particular, this framework naturally leads into two important local measurements, which are the {\bf error factor} and the {\bf condition factor}. Computing these two measurements only involves solving a localized eigenvalue problem and consequently no global computation or information is involved. These two measurements serve as guidances to define an appropriate partition of the graph $G$. Using this partition, a modified coarse space and corresponding basis with exponential decaying property can be constructed. Compression of $L^{-1}$ can thus be achieved. Furthermore, the systematic clustering procedure of the graph regardless of the knowledge of the geometric information allows us to introduce a multiresolution matrix decomposition (MMD) framework for graph laplacian, and more generally, SPD linear systems. In particular, following the work in \cite{owhadi2017multigrid}, we propose a nearly-linear time fast solver for general SPD matrices. Given the prescribed well-posedness requirement (i.e., the {\bf condition factor}), every component from MMD will be a well-conditioned, lower dimensional SPD linear system. Any generic iterative solver can then be applied in parallel to obtain the approximated solution of the given arbitrary SPD matrix satisfying the prescribed accuracy.

\subsection{Overview of our results}
Given a $n \times n$ SPD matrix $A$ with $m$ nonzero entries, our ultimate goal is to develop a fast algorithm to efficiently solve $Ax = b$, or equivalently, compress the solver $A^{-1}$ with desired compression error. We make the following assumptions on the matrix $A$. First of all, $\lambda_{\min}(A) = O(1)$ for well-posedness, where $\lambda_{\min}(A)$ is the minimum eigenvalue of the matrix $A$. Secondly, the spectrum of the matrix $A$ is broad-banded. Thirdly, it is stemmed from summation of symmetric and positive semidefinite (SPSD) matrices. We remark that the second assumption can be interpreted as the sparsity requirement of $A$, which is, the existence of some intrinsic, localized geometric information. For instance, if $A = L$ is a graph laplacian and $\mathcal{A}$ is the corresponding matrix, such sparsity can be described by the requirement
\begin{equation*}
\#\mathcal{N}_k(i) = O(k^d), \ \forall i,
\end{equation*}
where $\#\mathcal{N}_k(i)$ is the number of vertices near the vertex $i$ with logic distance smaller than $k$ (i.e., number of nonzero off-diagonal entries on row $i$ of $\mathcal{A}^k$) and $d$ is the geometric dimension of the graph (i.e., the optimal embedding dimension). This is equivalent to assuming that the portion of long interaction edges is small. The third assumption, in many concerning cases, is a natural consequence during the assembling of the matrix $A$. In particular, a graph laplacian $L$ can be viewed as a summation of 2-by-2 matrices representing edges in the graph $G$. These 2 by 2 matrices are SPSD matrices and can be obtained automatically if $G$ is given. Another illustrative example is the patch-wise stiffness matrix of a finite element from the discretization of PDEs using FE type methods.

\begin{figure}[!h]
\centering
\includegraphics[width=1.0\textwidth]{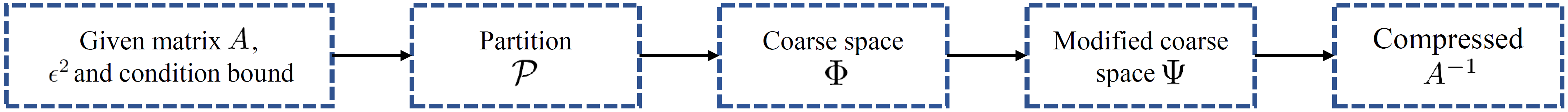}
\caption{Process flowchart for compressing $A^{-1}$.}
\label{fig:flowchart}
\end{figure}

To compress the solver $A^{-1}$ (where $A$ satisfies the above assumptions) with a desired error bound, we adopt the idea of constructing modified coarse space as proposed in \cite{maalqvist2014localization,owhadi2017multigrid,hou2016sparse}. The procedure is summarized in \Cref{fig:flowchart}. Noted that one common strategy of these PDE approaches is to make use of the natural partition under some a priori geometric assumption on the computational domain. In contrast, we adaptively construct an appropriate partition using the energy decomposition framework, which requires no a priori knowledge related to the underlying geometry of the domain. This partitioning idea is prerequisite and advantageous in the scenario when no explicit geometry information is provided, especially in the case of graph laplacian. Therefore, one of our main contributions is to develop various criteria and systematic procedures to obtain an appropriate partition $\mathcal{P} = \{P_j\}_{j=1}^M$ (i.e., graph partitioning in the case of graph laplacian) which reveals the structural property of $A$.

Leaving aside the difficulties of finding an appropriate partition, our next task is to define a coarse space $\Phi$ such that $\|x-P_\Phi x\|_2\leq \epsilon\|x\|_A$. As shown in \cite{owhadi2017multigrid,hou2016sparse}, having such requirement, together with the modification of coarse space $\Phi$ into $\Psi = A^{-1}(\Phi)$, we have $\|A^{-1} - P^A_{\Psi}A^{-1} \|_2 \leq \epsilon^2$. This affirms us that $\Phi$ must be constructed carefully in order to achieve the prescribed compression accuracy. Further, under the energy decomposition setting, we can ensure such requirement by simply considering local accuracy $\| x - P_{\Phi_j}x\|_2$, which in turns gives a local computation procedure for checking the qualification of the local $\Phi_j$. Specifically, we introduce the {\bf error factor} $\varepsilon(\mathcal{P},q) = \max_{P_j \in \mathcal{P}}\frac{1}{\sqrt{(\lambda_{q+1}(P_j))}}$ (where $\lambda_{q+1}(P_j)$ corresponds to the $(q+1)^{\text{th}}$ smallest eigenvalue of some eigen problem defined on the patch $P_j$) as our defender to strictly control the overall compression error. The {\bf error factor} guides us to construct a subspace $\Phi_j^q \in \text{span}(P_j)$ for every patch $P_j$ satisfying the local accuracy, and eventually, the global accuracy requirement. Afterwards, we apply the formulation of the modified coarse space $\Psi$ to build up the exponential decaying basis $\psi_j$ for the operator compression. To be general, we reformulate the coarse space modification procedure by using purely matrix-based arguments. In addition, this reformulation immediately introduces another criterion, called the {\bf condition factor} $\delta(P_j)$ over every patch $P_j \in \mathcal{P}$. This second measurement serves as another defender to control the well-posedness of our compressed operator. Similar to the {\bf error factor}, the {\bf condition factor} is a local measurement which can be obtained by solving a partial local eigen problem. This local restriction can naturally convey in a global sense to bound the maximum eigenvalue of the compressed operator. In particular, we prove that the compressed operator $A_{\text{st}}$ satisfies $\kappa(A_{\text{st}}) \leq \max_{P_j \in \mathcal{P}} \delta(P_j) \| A^{-1} \|_2$. 

Up to this point, we can see that the choice of the partition has a direct influence on both the accuracy and the well-posedness of the compressed operator. In this work, we propose a nearly-linear time partitioning algorithm which is purely matrix-based, and with complexity $O(d\cdot s^2 \cdot \log s \cdot n) + O(\log s \cdot n \cdot \log n)$, where $s$ is the average patch size and $d$ is the intrinsic geometric dimension. With the relationship between the {\bf error factor} and {\bf condition factor}, we can reversely treat the local properties as the blueprint to govern the process of partitioning. This in turn regularizes the condition number of the compressed operator such that computational complexity of solving the original linear system can be magnificently reduced.

Having a generic operator compression scheme, we follow the idea in \cite{owhadi2017multigrid} to extend the compression scheme hierarchically to form a multiresolution matrix decomposition algorithm. However, instead of using a precedent nested partitioning of the domain, we perform the decomposition level-by-level in a recursive manner. In other words, every new level is generated inductively from the previous level (subject to some uniform well-posedness constraints) by applying our adaptive partitioning technique. This provides more flexibility and convenience to deal with various, and even unknown multiresolution behavior appearing in the matrix. This decomposition further leads us to develop a fast solver for SPD matrices with time complexity $O(m \log n \cdot \left(\log \epsilon^{-1} + \log n \right)^c \log \epsilon^{-1})$, where $m$ is the number of nonzero entries of $A$ and $c$ is some absolute constant depending only on the geometric property of $A$. We would like to emphasize that the construction of the appropriate partition is essential to the formation of the hierarchical decomposition procedure. Principally, the hidden geometric information of $A$ can be subtly recovered from the {\bf inherited energy decomposition} of compressed operator using our framework. The interaction between these inherited energies serve similar purposes as the energy elements of $A$. Therefore we can recognize the compressed operator as an initial operator in the second level of the decomposition procedure and proceed to the next level repeatedly. We also remark that with our right choice of partitioning, the sparsity and the well-posedness properties of the compressed operator can be inherited among layers. This nice property enables us to decompose the original problem of solving $A^{-1}$ into sets of independent problems with similar complexity and condition, which favors the parallel implementation of the solver.

\subsection{Previous Works} 
Recently, several works relevant to the compression of elliptic operators with heterogeneous and highly varying coefficients have been proposed. M{\aa}lqvist and Peterseim et. al. \cite{maalqvist2014localization, larson2007adaptive} construct localized multiscale basis functions from the modified coarse space $V_{H}^{ms}=V_{H}-\mathfrak{F}V_{H}$, where $V_{H}$ is the original coarse space spanned by conforming nodal basis, and $\mathfrak{F}$ is the energy projection onto the space $(V_{H})^\perp$. The exponential decaying property of these modified basis has also been shown both theoretically and numerically. Meanwhile, a beautifully insightful work from Owhadi \cite{owhadi2017multigrid} reformulates the problem from the perspective of Decision Theory using the idea of {\it Gamblets} as the modified basis. In particular, a coarse space $\Phi$ of measurement functions is constructed from Bayesian perspective, and the gamblet space is explicitly given as $\Psi=A^{-1}(\Phi)$, which turns out to be a counterpart of the modified coarse space in \cite{maalqvist2014localization}. In addition, the basis of $\Phi$ is generalized to non-conforming measurement functions and the gamblets are still proven to decay exponentially such that localized computation is made possible. Hou and Zhang in \cite{hou2016sparse} extend these works such that localized basis functions can also be constructed for higher order strongly elliptic operators. Owhadi further generalizes these frameworks to a more unified methodology for arbitrary elliptic operators on Sobolev spaces in \cite{owhadi2017universal} using the Gaussian process interpretation. Noted that for the above-mentioned works, since the problems they considered are originated from PDE-type modeling, the computational domains are naturally assumed to be given, that is, the partition $\mathcal{P}$ can be obtained directly (which is not available for graph laplacians or general SPD matrices). This assumption greatly helps the immersion of the nested coarse spaces with different scales into the computational domain. In other words, the exponential decaying property of the basis can be precisely achieved.


Recall that for solving linear systems exhibiting multiple scales of behavior, the class of multiresolution methods decomposes the problem additively in terms of different scales of resolution. This captures the features of different scales and allows us to treat these components differently. For instance, the enlightening Geometric Multigrid (GMG) methods \cite{chan1994additive, wesseling2001geometric, sampath2010parallel} provide fast solvers for linear systems which are stemmed from discretization of linear elliptic differential equations. The main idea is to accelerate the convergence of basic iterative methods by introducing a nested structure on the computational domain so that successive subspace correction can be performed. However, the performance is hindered when the regularity of the coefficients is lost. To overcome this deficiency, an enormous amount of significant progress has been achieved. Numerous methods ranging from geometry specific to purely algebraic/ matrix-based approach have been developed (See \cite{stuben2001review, barrett1994templates, hackbusch2013multi} for review). Using the tools of compressing the operator $A^{-1}$ possessing multiple scale features, Owhadi in \cite{owhadi2017multigrid} also proposes a straightforward but intelligible way to solve the roughness issue. By introducing a natural nested structure on the given computational domain, a systematic multiresolution algorithm for hierarchically decomposing elliptic operators is proposed. This in turn derives a near-linear complexity solver with guaranteed prescribed error bounds. The efficiency of this multilevel solver is guaranteed by carefully choosing a nested structure of measurement functions $\mathrm{\Phi}$, which satisfies (i) the Poincar\'e inequality; (ii) the inverse Poincar\'e inequality; and (iii) the frame inequality. In \cite{owhadi2017universal}, Owhadi and Scovel extend the result to problems with general SPD matrices, where the existence of $\mathrm{\Phi}$ satisfying (i), (ii) and (iii) is assumed. In particular, for discretization of continuous linear bijections from $H^s(\Omega)$ to $H^{-s}(\Omega)$ or $L^2(\Omega)$ space these assumptions are shown to hold true using prior information on the geometry of the computational domain $\Omega$. However, the practical construction of this nested global structure $\mathrm{\Phi}$ is an essentially hard problem when no intrinsic geometric information is provided a priori. To solve this generic problem, we introduce the energy decomposition and the inherited system of energy elements. Instead of a priori assuming the existence of such nested structure $\Phi$, we use the idea of inherited energy decomposition to level-wisely construct $\Phi$ and the corresponding energy decomposition by using local spectral information and an adaptive clustering technique.

On the other hand, to mimic the functionality and convergence behavior of GMG without providing the nested meshes, the algebraic multigrid (AMG) methods \cite{stuben2001review, vanvek1996algebraic, trottenberg2000multigrid} speculate the coefficients through the connectivity information in the given matrix to define intergrid transfer operators, which avoids the direct construction of the restriction and relaxation operators in GMG methods. Intuitively, the connectivity information discloses the hidden geometry of the problem subtly. This purely algebraic framework bypasses the ``geometric" requirement in GMG, and is widely used in practice on graphs with sufficiently nice topologies. In particular, a recent AMG method called {\it LAMG} has been proposed by Livne and Brandt \cite{livne2012lean}, where the run time and storage of the algorithm are empirically demonstrated to scale linearly with the number of edges. We would like to emphasize that the difference between our proposed solver and a general recursive-typed iterative solver is the absence of nested iterations. Our solver decomposes the matrix adaptively according to the inherited multiple scales of the matrix itself. The matrix decomposition divides the original problem into components of controllable well-conditioned, lower dimensional SPD linear systems, which can then be solved in parallel using any generic iterative solver. In other words, this decomposition also provides a parallelizable framework for solving SPD linear systems. 

Another inspiring stream of nearly-linear time algorithm for solving graph laplacian system was given by Spielman and Teng \cite{spielman2004nearly, spielman2008local, spielman2011spectral, spielman2013local, spielman2014nearly}. With the innovative discoveries in spectral graph theory and graph algorithms, such as the fast construction low-stretch spanning trees and clustering scheme, they successfully employ all these important techniques in developing an effective preconditioned iterative solver. Later, Koutis, Miller and Peng \cite{koutis2011nearly} follow these ideas and simplify the solver with computation complexity $O(m \log n \log \log n \log \epsilon^{-1})$, where $m$ and $n$ are the number of edges and vertices respectively. In contrast, we employ the idea of modified coarse space to compress a general SPD matrix (i.e., the graph laplacian in this case) hierarchically with the control of sparsity and well-posedness. 

\subsection{Outline}
In \Cref{sec:energy_decomposition}, we will introduce the foundation of our work, which is the notion of {\bf Energy Decomposition} of general SPD matrices. \Cref{sec:operator_compression} discusses the construction of the coarse space and its corresponding modified coarse space, which serves to construct the basis with exponential decaying property. Concurrently, the local measurements {\bf error factor} and the {\bf condition factor} are introduced. The analysis in this section will guide us to design the systematic algorithm for constructing the partition $\mathcal{P}$, which is described in \Cref{Sec:partition}. Discussion of the computational complexity is also included. To demonstrate the efficacy of our partitioning algorithm, two numerical results are reported in \Cref{sec:numerical1}. Furthermore, we extend the idea of operator compression into Multiresolution Matrix Decomposition (MMD) in \Cref{Sec:Multiresolution}. In the meanwhile, we propose the concept of localization of MMD and the inherited locality of the compressed operator. These essential ingredients guides us to develop the parallelizable solver with nearly-linear time complexity. Error estimate and numerical results are reported to show the efficacy of this proposed algorithm. Conclusion and discussion of future works are included in \Cref{sec:conclusion}. For better readability, most of the proofs are moved to the Appendix section.

\section{Preliminaries}
In this paper we aim to solve the linear system $Lx=b$, where $L$ is the Laplacian of a undirected, positive-weighted graph $\bm{G}=\{\bm{V};\bm{E},\bm{W}\}$, i.e.
\begin{equation}
L_{ij} = \begin{cases} \sum_{(i,j')\in \bm{E}}w_{ij'} & \text{ if } i = j; \\ -w_{ij} & \text{ if } i \neq j\ \text{and}\ (i,j)\in \bm{E}; \\ 0 & \text{ otherwise}. \end{cases}
\end{equation}
We allow for the existence of selfloops $(i,i)\in \bm{E}$. When $L$ is singular we mean to solve $x=L^{\dagger}b$, where $L^\dagger$ is the pseudo-inverse of $L$. Our algorithm will base on a fast clustering technique using local spectral information to give a good partition of the graph, upon which special local basis will be constructed and used to compress the operator $L^{-1}$ into a low dimensional approximation $L^{-1}_{\text{com}}$ subject to a prescribed accuracy.

As we will see, our clustering technique exploits local path-wise information of the graph $\bm{G}$ by operating on each single edge in $\bm{E}$, which can be easily adapted to a larger class of linear systems with symmetric, positive semidefinite matrix. Notice that the contribution of an edge $(i,j)\in \bm{E}$ with weight $w_{ij}$ to the laplacian matrix $L$ is simply
\begin{equation}
E_{ii} \triangleq 
\begin{small}
\begin{blockarray}{cccc}
& i & \\
\begin{block}{(ccc)c}
  \textbf{\large $ 0 $} &  &  & \\
   & w_{ii} & & i \\
  & & \textbf{\large $0$} & \\
\end{block}
\end{blockarray}
\end{small},\ i=j;\quad \text{or}\quad
E_{ij} \triangleq 
\begin{small} 
\begin{blockarray}{cccccc}
& i & & j & \\
\begin{block}{(ccccc)c}
  \textbf{\large $ \ 0 $} &  &  &  &  &   \\
   & w_{ij} &  & -w_{ij} &  & i \\
   &  & \ddots &  &  &   \\
   & -w_{ij} &  & w_{ij} &  & j \\
   &  &  &  & \textbf{\large $ \ 0 $} & \\
\end{block}
\end{blockarray}
\end{small},\ i\neq j,
\label{Graph_energy_element}
\end{equation}
and we have $L=\sum_{(i,j)\in \bm{E}}E_{ij}$. In view of such matrix decomposition, our algorithm works for any symmetric, positive semi-definite matrix $A$ that has a similar decomposition $A=\sum_{k=1}^mE_k$ with each $E_k\succeq 0$. Therefore, we will theoretically develop our method for general decomposable SPD matrices. Also we assume that $A$ is invertible, as we can easily generalize our method to the case when $A^\dagger b$ is pursued.

\subsection{Energy Decomposition}
\label{sec:energy_decomposition}

In this section, we will introduce the idea of energy decomposition and the corresponding mathematical formulation which motivates the methodology for solving linear systems with energy decomposable linear operator. Let $A$ be a $n\times n$ symmetric, positive definite matrix. We define the {\bf Energy Decomposition} as follows:
\label{subsec:notation}

\begin{definition}[Energy Decomposition]
We call $\{ E_k \}_{k=1}^m$ an \textbf{energy decomposition} of $A$ and $E_k$ to be an \textbf{energy element} of $A$ if 
\begin{equation}
A = \sum_{k=1}^m E_k,\quad E_k \succeq 0\ \forall k = 1,\ldots, m,
\end{equation}
\label{def:energydecomp}
\end{definition}
where $E_k \succeq 0$ means $E_k$ is positive semidefinite. Intuitively, the underlying structural(geometric) information of the original matrix $A$ can be realized through an appropriate energy decomposition. And to preserve as much detailed information of $A$ as possible, it's better to use the finest energy decomposition that we can have, which actually comes naturally from the generating of $A$ as we will see in  some coming examples. More precisely, for an energy decomposition $\mathcal{E}=\{E_k\}_{k=1}^m$ of $A$, if there is some $E_k$ that has its own energy decomposition $E_k=E_{k,1}+E_{k,2}$ that comes naturally, then the finer energy decomposition $\mathcal{E}_{fine}=\{E_k\}_{k=2}^m\cup\{E_{k,1},E_{k,2}\}$ is more preferred as it gives us more detailed information of $A$. However one would see that any $E_k$ can have some trivial decomposition $E_k=\frac{1}{2}E_k+\frac{1}{2}E_k$, which makes no essential difference. To make it clear what should be the finest underlying energy decomposition of $A$ that we will use in our algorithm, we first introduce the neighboring relation between energy elements and basis.

Let $\mathcal{E}=\{E_k\}_{k=1}^m$ be an energy decomposition of $A$, and $\mathcal{V}=\{v_i\}_{i=1}^n$ be an orthonormal basis of $\mathbb{R}^n$. we introduce the following notation:
\begin{itemize}[leftmargin=*]
\item For any $E \in \mathcal{E}$ and any $v \in \mathcal{V}$, we denote $E \sim v$ if $v^T Ev > 0$ ( or equivalently $Ev \neq \bm{0}$, since $E \succeq \bm{0}$ );
\item For any $u,v\in \mathcal{V}$, we denote $u \sim v$ if $u^T Av \neq 0$ ( or equivalently $\exists E\in \mathcal{E}$ such that $u^TEv\neq0$ ).
\end{itemize}

As an immediate example, if we take $\mathcal{V}$ to be the set of all distinct eigen vectors of $A$, then $v\not\sim u$ for any two $v,u\in \mathcal{V}$, namely all basis functions are isolated and everything is clear. But such choice of $\mathcal{V}$ is not trivial in that we know everything about $A$ if we know its eigen vectors. Therefore, instead of doing things in the frequency space, we assume the least knowledge of $A$ and work in the physical space, that is we will choose $\mathcal{V}$ to be the natural basis $\{\bm{e}_i\}_{i=1}^n$ of $\mathbb{R}^n$ in all practical use. But for theoretical analysis, we still use the general basis notation $\mathcal{V}=\{v_i\}_{i=1}^n$.

Also, for those who are familiar with graph theory, it's more convenient to understand the sets $\mathcal{V},\mathcal{E}$ from graph perspective. Indeed, one can keep in mind that $G=\{\mathcal{V},\mathcal{E}\}$ is the generalized concept of undirected graphs, where $\mathcal{V}$ stands for the set of vertices, and $\mathcal{E}$ stands for the set of edges. For any vertices (basis) $v,u\in \mathcal{V}$, and any edge (energy) $E\in \mathcal{E}$, $v\sim E$ means that $E$ is an edge of $v$, and $v\sim u$ means that $v$ and $u$ share some common edge. However, different from the traditional graph setting, here one edge(energy) $E$ may involve multiple vertices instead of just two, and two vertices(basis) $v,u$ may share multiple edges that involve different sets of vertices. Further, the spectrum magnitude of the ``multi-vertex edge" $E$ can be viewed as an counterpart of edge weight in graph setting. Conversely, if the problem comes directly from a weighted graph, then one can naturally construct the sets $\mathcal{V}$ and $\mathcal{E}$ from the vertices and edges of the graph as we will see in \Cref{example:Example2}.

\begin{definition}[Neighboring] Let $\mathcal{E}=\{E_k\}_{k=1}^m$ be an energy decomposition of $A$, and $\mathcal{V}=\{v_i\}_{i=1}^n$ be an orthonormal basis of $\mathbb{R}^n$. For any $E\in \mathcal{E}$, We define $\mathcal{N}(E;\mathcal{V}) :=\{v\in \mathcal{V}: E\sim v\}$ to be the set of $v \in \mathcal{V}$ \textbf{neighboring} $E$. Similarly, for any $v\in \mathcal{V}$, we define $\mathcal{N}(v;\mathcal{E}) := \{E\in \mathcal{\mathcal{E}}: E\sim v\}$ and $\mathcal{N}(v) := \{u\in \mathcal{V}: u\sim v\}$
to be the set of $E \in \mathcal{E}$ and the set of $u \in \mathcal{V}$ \textbf{neighboring} $v \in \mathcal{V}$ respectively. Furthermore, for any $\mathcal{S} \subset \mathcal{V}$ and any $E\in \mathcal{E}$, we denote $E\sim \mathcal{S}$ if $\mathcal{N}(E;\mathcal{V})\cap \mathcal{S}\neq \emptyset$ and $E \in \mathcal{S}$ if $\mathcal{N}(E;\mathcal{V})\subset \mathcal{S}$.
\label{def:neighboring}
\end{definition}

In what follows, we will see that if two energy elements $E_k,E_{k'}$ have the same neighbor basis, namely $\mathcal{N}(E_k;\mathcal{V})=\mathcal{N}(E_{k'};\mathcal{V})$, then there is no need to distinguish between them, since it is the neighboring relation between energy elements and basis that matters in how we make use of the energy decomposition. Therefore we say an energy decomposition $\mathcal{E}=\{E_k\}_{k=1}^m$ is the finest underlying energy decomposition of $A$ if no $E_k\in \mathcal{E}$ can be further decomposed as 
\[E_k=E_{k,1}+E_{k,2},\]
where either $\mathcal{N}(E_{k,1};\mathcal{V})\varsubsetneqq \mathcal{N}(E_{k};\mathcal{V})$ or $\mathcal{N}(E_{k,2};\mathcal{V})\varsubsetneqq \mathcal{N}(E_{k};\mathcal{V})$. From now on, we will always assume that $\mathcal{E}=\{E_k\}_{k=1}^m$ is the finest underlying energy decomposition of $A$ that comes along with $A$.

Using the {\it neighboring} concept between energy elements and orthonormal basis, we can then define various energies of a subset $\mathcal{S} \subset \mathcal{V}$ as follows:

\begin{definition}[Restricted, Interior and Closed energy]

Let $\mathcal{E}=\{E_k\}_{k=1}^m$ be a energy decomposition of $A$. Let $\mathcal{S}$ be a subset of $\mathcal{V}$, and $P_\mathcal{S}$ be the orthogonal projection onto $\mathcal{S}$. The \textbf{restricted energy} of $\mathcal{S}$ with respect to $A$ is defined as
\begin{equation}
A_\mathcal{S} := P_\mathcal{S} A P_\mathcal{S};
\end{equation}
The \textbf{interior energy} of $\mathcal{S}$ with respect to $A$ and $\mathcal{E}$ is defined as
\begin{equation}
\underline{A}^\mathcal{E}_\mathcal{S}=\sum_{E\in \mathcal{S}}E;
\end{equation}
The \textbf{closed energy} of $\mathcal{S}$ with respect to $A$ and $\mathcal{E}$ is defined as
\begin{equation}
\overline{A}^\mathcal{E}_\mathcal{S}=\sum_{E\in \mathcal{S}}E+\sum_{E\notin \mathcal{S},E\sim S} P_\mathcal{S}E^{d}P_\mathcal{S},
\end{equation}
where
\begin{equation}
E^{d}=\sum_{v\in \mathcal{V}}\Big(\sum_{u\in \mathcal{V}} \big|v^TEu\big|\Big) vv^T=\sum_{v\sim E}\Big(\sum_{u\sim v} \big|v^TEu\big|\Big) vv^T
\end{equation}
is called the \textbf{diagonal concentration} of $E$, and we have
\begin{equation}
P_\mathcal{S}E^{d}P_\mathcal{S}=\sum_{v\in \mathcal{S},v\sim E}\Big(\sum_{u\sim v} \big|v^TEu\big|\Big) vv^T
\end{equation}
\label{def:energies}
\end{definition}

\begin{remark}
The {\it restricted energy} of $\mathcal{S}$ can be simply viewed as the restriction of $A$ on the subset $\mathcal{S}$. The {\it interior energy} ({\it closed energy}) of $\mathcal{S}$ is $A_{\mathcal{S}}$ excluding (including) contributions from other energy elements $E \not\in \mathcal{S}$ neighboring $\mathcal{S}$. The following example illustrates the idea of various energies introduced in \Cref{def:energies} by considering the 1-dimensional discrete Laplace operator with Dirichlet boundary conditions.
\end{remark}

\bigskip

\begin{example}
Consider $A$ to be the $(n+1) \times (n+1)$ tridiagonal matrix with entries -1 and 2 on off-diagonals and diagonal respectively. Let 
\begin{equation}
E_1 = \left( \begin{smallmatrix} 2 & -1 & \\ -1 & 1 & \\ & & \textbf{\Large $0 \ $} \end{smallmatrix} \right),\ 
E_{n} = \left( \begin{smallmatrix} \textbf{\Large $\ 0 $} & & \\ & 1 & -1 \\ & -1 & 2\end{smallmatrix} \right),\ 
E_{k} = \left( \begin{smallmatrix} \textbf{\large $ \ 0 $} & & \\ & 1 & -1 & \\ & -1 & 1 & \\ & & & \textbf{\large $ 0 \ $}\end{smallmatrix} \right)
\end{equation}
for $k = 2,\ldots n-1$. Let $\mathcal{V} = \left\{\bf{e}_i\right\}_{i=0}^{n}$ to be the standard orthonormal basis for Euclidean space $\mathbb{R}^{n+1}$. Formally $E_k$ is the edge between $e_{k-1}$ and $e_k$. If $\mathcal{S} = \left\{\bf{e}_3,{\bf e}_4,{\bf e}_5,{\bf e}_6\right\}$, then we have
\begin{equation*}
A_{\mathcal{S}} = \left( \begin{smallmatrix} \textbf{\large $ \ 0 $} & & & & & \\ & 2 & -1 & & & \\ & -1 & 2 & -1 & & \\ & & -1 & 2 & -1 & \\ & & & -1 & 2 & \\ & & & & & \textbf{\large $ \ 0 $} \end{smallmatrix} \right),\ \underline{A}^\mathcal{E}_\mathcal{S} = \left( \begin{smallmatrix} \textbf{\large $ \ 0 $} & & & & & \\ & 1 & -1 & & & \\ & -1 & 2 & -1 & & \\ & & -1 & 2 & -1 & \\ & & & -1 & 1 & \\ & & & & & \textbf{\large $ \ 0 $} \end{smallmatrix} \right),\text{ and }\
\overline{A}^\mathcal{E}_\mathcal{S} = \left( \begin{smallmatrix} \textbf{\large $ \ 0 $} & & & & & \\ & 3 & -1 & & & \\ & -1 & 2 & -1 & & \\ & & -1 & 2 & -1 & \\ & & & -1 & 3 & \\ & & & & & \textbf{\large $ \ 0 $} \end{smallmatrix} \right). 
\end{equation*}
Recall that the {\it interior energy} $\underline{A}^\mathcal{E}_\mathcal{S} = \sum_{E_k \in \mathcal{S}} E_k = \sum_{k=4}^6 E_k$, while the {\it closed energy}
\begin{align*} \overline{A}^\mathcal{E}_\mathcal{S} 
=& \underline{A}^\mathcal{E}_\mathcal{S} + \sum_{E\not\in \mathcal{S}, E \sim \mathcal{S}} P_{\mathcal{S}} E^d P_{\mathcal{S}} \\
=& \underline{A}^\mathcal{E}_\mathcal{S} +  \vert {\bf e}_3^T E_3 {\bf e}_2 \vert {\bf e}_3 {\bf e}_3^T +  \vert {\bf e}_3^T E_3 {\bf e}_3 \vert {\bf e}_3 {\bf e}_3^T +  \vert {\bf e}_7^T E_7 {\bf e}_6 \vert {\bf e}_6 {\bf e}_6^T + \vert {\bf e}_6^T E_7 {\bf e}_6 \vert {\bf e}_6 {\bf e}_6^T
\end{align*}
includes the partial contributions from other energy elements $E \not\in \mathcal{S}$ neighboring $\mathcal{S}$, which are $E_3$ and $E_7$ respectively.
\label{ex:example1}
\end{example}

\begin{remark}
$ $\linebreak
\vspace{-2mm}
\begin{itemize}
\item[-] Notice that any eigenvector $x$ of $A_\mathcal{S}$ (or $\underline{A}^\mathcal{E}_\mathcal{S}$, $\overline{A}^\mathcal{E}_\mathcal{S}$) corresponding to non-zero eigenvalue must satisfy $x\in \mathrm{span}(\mathcal{S})$. In this sense, we also say $A_\mathcal{S}$ (or $\underline{A}^\mathcal{E}_\mathcal{S}$, $\overline{A}^\mathcal{E}_\mathcal{S}$) is local to $\mathcal{S}$.
\item[-] For any energy $E$, we have $E \preceq E^d$ since for any $x=\sum_{i=1}^nc_iv_i$, we have
\begin{align*}
x^TEx=&\ \sum_{i} c_i^2v_i^TEv_i+\sum_{i\neq j}2 c_ic_jv_i^TEv_j\\
& \leq\ \sum_{i} c_i^2v_i^TEv_i+\sum_{i\neq j}(c_i^2+c_j^2)\big|v_i^TEv_j\big|\ =\ \sum_{i} \sum_{j}c_i^2\big|v_i^TEv_j\big| = x^TE^dx.
\end{align*}
\end{itemize}
\end{remark}

\begin{proposition}
For any $\mathcal{S}\subset\mathcal{V}$, we have that $\underline{A}^\mathcal{E}_\mathcal{S} \preceq A_\mathcal{S} \preceq \overline{A}^\mathcal{E}_\mathcal{S}$.
\end{proposition}
\begin{proof}
We have
\begin{align*}
\underline{A}^\mathcal{E}_\mathcal{S} &= \sum_{E\in\mathcal{S}} E \preceq \sum_{E\in \mathcal{S}}E + \sum_{\mathclap{\substack{E\notin \mathcal{S}, \\ E\sim S}}} P_\mathcal{S} E P_\mathcal{S} \preceq \sum_{E\in \mathcal{S}}E + \sum_{\mathclap{\substack{E \notin \mathcal{S},\\E\sim S}}} P_\mathcal{S}E^d P_\mathcal{S} = \overline{A}^\mathcal{E}_\mathcal{S}.
\end{align*}
Notice that $P_\mathcal{S}EP_\mathcal{S}=E$ for $E\in \mathcal{S}$, and $P_\mathcal{S}EP_\mathcal{S}=0$ for $E\not\sim \mathcal{S}$, thus 
\[A_\mathcal{S}=P_\mathcal{S}AP_\mathcal{S}=\sum_{E\in \mathcal{S}}E+\sum_{\mathclap{\substack{E\notin \mathcal{S}, \\ E\sim S}}} P_\mathcal{S}E P_\mathcal{S},\]
and the desired result follows.
\end{proof}

\begin{definition}[Partition of basis]
Let $\mathcal{V}=\{v_i\}_{i=1}^n$ be an orthonormal basis of $\mathbb{R}^n$. We say $\mathcal{P}=\{P_j\}_{j=1}^M$ is a \textbf{partition} of $\mathcal{V}=\{v_i\}_{i=1}^n$ if (i) $P_j\subset \mathcal{V} \ \forall j$; (ii) $P_j\cap P_{j'}=\emptyset$ if $j\neq j'$; and (iii) $\bigcup_{j=1}^M P_j = \mathcal{V}$.
\label{def:partition}
\end{definition}

Again one can see the partition of basis as partition of vertices. This partition $\mathcal{P}$ is the key to construction of local basis for operator compression purpose. The following proposition serves to bound the matrix $A$ from both sides with blocked(patched) matrices, which will further serve to characterize properties of local basis. 

\begin{proposition}
Let $\mathcal{E}=\{E_k\}_{k=1}^m$ be an energy decomposition of $A$, and $\mathcal{P}=\{P_j\}_{j=1}^M$ be a partition of $\mathcal{V}$. Then 
\begin{equation}
\sum_{j=1}^M\underline{A}^\mathcal{E}_{P_j} \preceq A \preceq \sum_{j=1}^M\overline{A}^\mathcal{E}_{P_j}.
\end{equation}
\label{prop:OED}
\end{proposition}
\begin{proof} Let $\mathcal{E}_\mathcal{P}=\{ E\in \mathcal{E}: \exists P_j\in\mathcal{P}\ \text{such that} \ E\in P_j\}$, and $\mathcal{E}_\mathcal{P}^c=\mathcal{E}\backslash\mathcal{E}_\mathcal{P}$. Recall that $E \in P_j$ if $\mathcal{N}(E,\mathcal{V}) \subset P_j$ (See \cref{def:neighboring}). We will use $P_j$ to denote the orthogonal projection onto $P_j$. Since $P_j\cap P_{j'}=\emptyset$ for $j\neq j'$, we have $\sum_{j}P_j= \textbf{Id}$. Then 
\begin{align*}
\sum_{j}\underline{A}^\mathcal{E}_{P_j}=&\ \sum_{j}\sum_{E\in P_j}E\ \preceq \ \sum_{j}\sum_{E\in P_j}E+\sum_{E\in \mathcal{E}^c_\mathcal{P}}E\\
& \preceq \ \sum_{j}\sum_{E\in P_j}E+\sum_{E\in \mathcal{E}^c_\mathcal{P}}\Big[\big(\sum_{j}P_j\big)E^d\big(\sum_{j'}P_{j'}\big)\Big] \ =\ \sum_{j}\sum_{E\in P_j}E+\sum_{E\in \mathcal{E}^c_\mathcal{P}}\big(\sum_{j}P_jE^dP_j\big)\\
&\quad =\ \sum_{j}\Big(\sum_{E\in P_j}E+\sum_{E\notin P_j,E\sim P_j}P_jE^dP_j\Big) \ =\ \sum_{j}\overline{A}^\mathcal{E}_{P_j}.
\end{align*}
We have used the fact that $P_jE^dP_{j'}=0$ for $j\neq j'$. Notice that 
\[A=\sum_{E\in \mathcal{E}_\mathcal{P}}E+\sum_{E\in \mathcal{E}^c_\mathcal{P}}E=\sum_{j}\sum_{E\in P_j}E+\sum_{E\in \mathcal{E}^c_\mathcal{P}}E,\]
and the desired result follows.
\end{proof}

Throughout the paper, we will always assume that $A$ has a finest energy decomposition $\mathcal{E}=\{E_k\}_{k=1}^m$, and all the other discussed energies of $A$ are constructed from $\mathcal{E}$ with respect to some orthonormal basis $\mathcal{V}$ (by taking interior or closed energy). Therefore we will simply use $\underline{A}_\mathcal{S},\overline{A}_\mathcal{S}$ to denote $\underline{A}^\mathcal{E}_\mathcal{S},\overline{A}^\mathcal{E}_\mathcal{S}$ for any $S\subset\mathcal{V}$.

\begin{example}
\begin{figure}[h!]
\centering
\begin{subfigure}[b]{0.4\textwidth}
        \centering
        \includegraphics[height=1.8in]{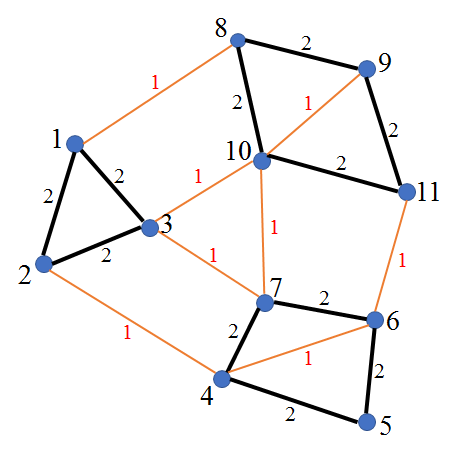}
        \caption{}
        \label{fig:example2a}
    \end{subfigure}%
    ~ 
    \begin{subfigure}[b]{0.4\textwidth}
        \centering
        \includegraphics[height=1.8in]{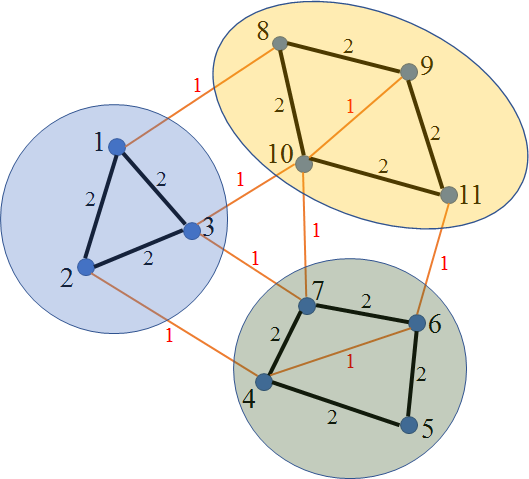}
        \caption{}
        \label{fig:example2b}
    \end{subfigure}
    \caption{(a) An illustration of a graph example . (b) An illustration of a partition $\mathcal{P} = \{ \{1,2,3\},\{4,5,6,7\},\{8,9,10,11\}\}$.}
    \label{fig:example2}
\end{figure}
Consider $L$ to be the graph laplacian matrix of the graph given in \Cref{fig:example2}. For graph laplacian, an intrinsic energy decomposition arises during the assembling of the matrix in which the energy element is defined over each edge(see \Cref{Graph_energy_element}). Now suppose we have given the partition $\mathcal{P} = \left\{ P_j \right\}_{j=1}^3$ with $P_1 = [{\bf e}_1, {\bf e}_2, {\bf e}_3]$, $P_2 = [{\bf e}_4, {\bf e}_5, {\bf e}_6, {\bf e}_7]$ and $P_3 = [{\bf e}_8, {\bf e}_9, {\bf e}_{10}, {\bf e}_{11}]$, where ${\bf e}_i$ are the standard basis of $\mathbb{R}^{11}$. Then we can obtain $\underline{L}_{P_j}$ and $\overline{L}_{P_j}$ as follows:
\begin{equation*}
\underline{L}_{P_1} = \left( \begin{smallmatrix} 4 & -2 & -2 \\ -2 & 4 & -2 \\ -2 & -2 & 4 \end{smallmatrix} \right)_{P_1},\ \underline{L}_{P_2} = \left( \begin{smallmatrix} 5 & -2 & -1 & -2 \\ -2 & 4 & -2 & 0 \\ -1 & -2 & 5 & 2 \\ -2 & 0 & -2 & 4 \end{smallmatrix} \right)_{P_2}, \ \underline{L}_{P_3} = \left( \begin{smallmatrix} 4 & -2 & -2 & 0 \\ -2 & 5 & -1 & -2 \\ -2 & -1 & 5 & -2 \\ 0 & -2 & -2 & 4 \end{smallmatrix} \right)_{P_3}
\end{equation*}
\begin{equation*}
\overline{L}_{P_1} = \left( \begin{smallmatrix} 6 & -2 & -2 \\ -2 & 6 & -2 \\ -2 & -2 & 8 \end{smallmatrix} \right)_{P_1},\ \overline{L}_{P_2} = \left( \begin{smallmatrix} 7 & -2 & -1 & -2 \\ -2 & 4 & -2 & 0 \\ -1 & -2 & 7 & 2 \\ -2 & 0 & -2 & 8 \end{smallmatrix} \right)_{P_2}, \ \overline{L}_{P_3} = \left( \begin{smallmatrix} 6 & -2 & -2 & 0 \\ -2 & 5 & -1 & -2 \\ -2 & -1 & 9 & -2 \\ 0 & -2 & -2 & 6 \end{smallmatrix} \right)_{P_3}
\end{equation*}
Here we denote the matrix $\left( \cdot \right)_{P_j}$ to be the matrix in $\mathbb{R}^{11\times11}$ but with non-zero entries on $P_j$ only. 
\label{example:Example2}
\end{example}

\section{Operator Compression}
\label{sec:operator_compression}

As mentioned in the \Cref{sec:introduction}, inspired by the FEM approach for solving partial differential equation (PDE) in which the variational formulation naturally gives the energy decomposition of the operator, we adopt a similar strategy of FEM to find a subspace $\Phi$ approximating the solution space of a linear system involving $A$ that are energy decomposable. In particular, approximation of $A^{-1}$ can also be obtained.

For traditional FEM, the accuracy of these approximations relies on the regularity of the given coefficients. Without assuming any smoothness on coefficients, one promising way to approximate the operator is to consider projecting the operator into the modified subspace $\Psi = \Phi - P_{U}^A\Phi$ in \cite{maalqvist2014localization}, or $\Psi = A^{-1} (\Phi)$ as proposed in \cite{owhadi2017multigrid,hou2016sparse}. Here $U = (\Phi)^{\perp}$ is the $l_2$-orthogonal complement space and $P_{U}^A$ is the A-orthogonal projection operator. In the case where $A$ is invertible, these two modified spaces are equivalent. Therefore, we propose to employ a similar methodology for compressing a general symmetric, positive definite matrix $A$.

We first obtain a general error estimate for projecting the matrix $A$ into a subspace $\Psi = A^{-1}(\Phi)$ of $\mathbb{R}^n$, given the projection type approximation property of the subspace $\Phi$. With this observation, the operator compression problem is narrowed down into choosing an appropriate $\Phi$ which satisfies condition \cref{eqt:Phi_errortolerance}. The following lemma also gives us a general idea on how we can control the errors introduced during the compression of the operator $A$. 

\begin{lemma}
Let $\Phi$ be a subspace of $\mathbb{R}^n$, and $P_{\Phi}$ be the orthogonal projection onto $\Phi$ with respect to $\langle\cdot,\cdot\rangle_2$. Let $\Psi$ be the subspace of $\mathbb{R}^n$ given by $\Psi=A^{-1}(\Phi)$, and $P_{\Psi}^{A}$ be the orthogonal projection onto $\Psi$ with respect to $\langle\cdot,\cdot\rangle_A$. If 
\begin{equation}
\|x-P_\Phi x\|_2\leq \epsilon\|x\|_A,\quad \forall x\in \mathbb{R}^n,
\label{eqt:Phi_errortolerance}
\end{equation}
for some $\epsilon>0$, then 
\begin{enumerate}
\item For any $x \in \mathbb{R}^n$, and $b=Ax$, we have
\begin{equation}
\|x-P_\Psi^A x\|_A\leq \epsilon\|b\|_2.
\label{eqt:Aorthogonalprojection_error}
\end{equation}
\item For any $x \in \mathbb{R}^n$, and $b=Ax$, we have
\begin{equation}
\|x-P_\Psi^A x\|_2\leq \epsilon^2\|b\|_2.
\end{equation}
\item We have
\begin{equation}\|A^{-1}-P_\Psi^A A^{-1}\|_2\leq \epsilon^2.
\end{equation}
\end{enumerate}
\label{lemma:phi_error}
\end{lemma}

Interchangeably, we will use $\Phi$ and $\Psi$ to denote the basis matrix of the space $\Phi$ and $\Psi$ respectively, such that $\Phi^T\Phi=I_N$ and $\Psi=A^{-1}\Phi T$. Here $N$ is the dimension of $\Phi$, and $T$ is some $N\times N$ nonsingular matrix to be determined. Then we have
\begin{equation}
P_\Psi^A=\Psi(\Psi^TA\Psi)^{-1}\Psi^TA=A^{-1}\Phi(\Phi^TA^{-1}\Phi)^{-1}\Phi^T,
\label{eqt:Aorthogonal_projection}
\end{equation}
which is the $A$-orthogonal projection matrix into the subspace $\Psi$. In \cite{maalqvist2014localization}, M{\aa}lqvist and Petersein proposed the use of modified coarse space in order to handle roughness of coefficients when solving elliptic equations with FEM. Assuming that the finite elements are conforming and if we see $\Phi$ as the original coarse space $V_H$ in \cite{maalqvist2014localization}, then $\Psi$ is exactly the modified coarse space $V_H^{ms}$ as they proposed, and the first error estimate in \Cref{lemma:phi_error} is consistent to their error analysis. More generally, Owhadi in \cite{owhadi2017universal} makes use of the Gamblet framework to construct the basis of modified coarse space such that the conforming properties of those basis is no long required. In particular, \Cref{eqt:Aorthogonal_projection} is an analogy of $\Psi\Phi = K \Psi^T (\Phi K \Phi^T)^{-1} \Phi$ in page 9 of \cite{owhadi2017multigrid} and the error estimate in \cref{eqt:Aorthogonalprojection_error} is correspondingly the Proposition 3.6 in that paper.

As a FEM type method, the choice of $\Phi$ determines the operator compression error or the solution approximation error. We know that the optimal rank-$N$ approximator of $A^{-1}$ is given by taking $\Phi$ to be the eigenspace of $A$ corresponding to the first $N$ smallest eigenvalues, which is essentially the Principal Component Analysis (PCA) \cite{jolliffe2002principal}. And the optimal compression error is given by $\epsilon^2=(\lambda_{N+1}(A))^{-1}$. Though with optimal approximation property, the drawback of the PCA is non-negligible in that the eigenvectors of $A$ are almost always dense even when $A$ has strong local properties. While the sparse PCA \cite{moghaddam2006spectral, zou2006sparse, d2007direct} provides a strategy to obtain a sparse approximation of $A^{-1}$, it implicitly assumes that the operators inherit a low rank characteristics such that $l_1$ minimization approach is effective. To fully use the local properties of $A$, we would prefer to choose $\Phi$ that can be locally computed but still has good approximation property, namely satisfying condition \cref{eqt:Phi_errortolerance} with a pretty good error $\epsilon$ and a nearly optimal dimension $N$. Also we hope the a priori error bound $\epsilon$ can be estimated locally.

Indeed when solving elliptic PDEs using FEM, the nodal basis can be chosen as (discretized) piece-wise polynomials with compact local supports, and the error is given by the resolution of the partition of the computational domain \cite{bathe1976numerical, brenner2004finite}. However, such choices of partition of the computational domain do not depend on the operator $A$ in traditional FEM. Yet it only depends on the geometry of the computational domain, and thus the performance relies on the regularity of $A$. So a natural question arises: can we do better if we choose the partition and the nodal basis using the local information of $A$?

Furthermore, what can we do if we don't a priori have the computational domain? Such scenario arises, for example, when the underlying geometry of some operators $A$, like graph laplacian, is unknown and no embedding maps to the physical domain can be found easily. In this case, one of the promising ways to accomplish such task lies in the deep connection between our {\it energy representation} of the operator $A$ and its hidden geometric structure. More specifically, the energy decomposition of the operator introduced in \Cref{sec:energy_decomposition} reveals the intrinsic locality of the underlying geometry in an algebraic way, so that we can construct an optimal partition of the computational space and choose a proper subspace/basis $\Phi$ using only local information.

After constructing the partition and the basis $\Phi$, the next mission is to find a good basis $\Psi$ of the space $A^{-1}(\Phi)$. The choice of $\Psi$ serves to preserve the locality of the stiffness matrix $A_{\text{st}} = \Psi^T A \Psi$ inherited from $A$, and to give a reasonable bound on the condition number of $A_{\text{st}}$.

\begin{algorithm}[!h]
\caption{\it Operator Compression}
\label{alg:operator_compression}
\begin{algorithmic}[1]
\REQUIRE{Energy decomposition $\mathcal{E}$, underlying basis $\mathcal{V}$, desire accuracy $\epsilon$}
\STATE{Construct partition $\mathcal{P}$ subject to $\epsilon$ using \Cref{alg:pair_clustering};}
\STATE{Construct $\Phi$ using \Cref{alg:construct_phi};}
\STATE{Construct $\widetilde{\Psi}$ using \Cref{alg:construct_tilde_psi} subject to $\epsilon$;}
\STATE{Compute $P_{\widetilde{\Psi}}^A A^{-1}=\widetilde{\Psi}(\widetilde{\Psi}^TA\widetilde{\Psi})^{-1}\widetilde{\Psi}^T$ as the compressed operator;}
\end{algorithmic}
\end{algorithm}

In summary, as mentioned in \Cref{sec:introduction}, our approach is to (i) construct a partition of the computational space/basis using local information of $A$; (ii) construct $\Phi$ that is locally computable in each patch of the partition and satisfies error condition \cref{eqt:Phi_errortolerance}; (iii) construct $\Psi$ that provides stiffness matrix $A_{\text{st}}$ with locality and reasonable condition number. The whole process can be summarized as the \Cref{alg:operator_compression}, and each step will be discussed in following sections.

But to theoretically develop our approach, we first assume that we are given an imaginary partition $\mathcal{P}$, and then derive proper constructions of $\Phi$ and $\Psi$ serving the desired purposes based on this partition. In the derivation process, we come up with some desired conditions that will, in return, guide us how to construct the adaptive partition $\mathcal{P}$ with the desirable properties.

\subsection{Choice of $\Phi$}
\label{subsec:choose_phi}
As discussed in the last subsection, the underlying geometry of the operator may not be given. Therefore determining $\Phi$ which archives the condition \cref{eqt:Phi_errortolerance} is not a trivial task. Instead of tackling this problem directly, the following proposition provides us a more apparent and local criterion on choosing $\Phi$.  

\begin{proposition}
Let $\mathcal{P}=\{P_j\}_{j=1}^M$ be a partition of $\mathcal{V}$, and $\{\underline{A}_{P_j}\}_{j=1}^M$ be the corresponding interior energies as defined in \Cref{def:energies}. For each $1\leq j\leq M$, let $\Phi_j$ be some subspace of $\mathrm{span}\{P_j\}$ such that 
\begin{equation}
\|x-P_{\Phi_j}x\|_2\leq \epsilon \|x\|_{\underline{A}_{P_j}},\quad \forall x\in \mathrm{span}\{P_j\},
\label{eqt:localphi_constraint}
\end{equation}
for some constant $\epsilon$. Then we have
\begin{equation}
\|x-P_{\Phi}x\|_2\leq \epsilon \|x\|_A,\quad \forall x\in \mathbb{R}^n,
\label{eqt:phi_constraint}
\end{equation}
where $\Phi=\bigoplus_{j}\Phi_j$.
\label{prop:localphi}
\end{proposition}

\begin{proof} Since $\mathcal{P}=\{P_j\}_{j=1}^M$ is a partition of $\mathcal{V}$, we have $P_\Phi=\sum_jP_{\Phi_j}$, and thus
\[\|x-P_\Phi x\|_2^2=\|\sum_j(P_jx-P_{\Phi_j}x)\|_2^2=\sum_j\|P_jx-P_{\Phi_j}x\|_2^2 \leq\epsilon^2\sum_j\|P_jx\|^2_{\underline{A}_{P_j}}.\]
Notice that 
\[\sum_j\|P_jx\|^2_{\underline{A}_{P_j}}=\sum_j\|x\|^2_{\underline{A}_{P_j}}\leq\|x\|_A^2,\]
and the conclusion follows.
\end{proof}

Intuitively, given a partition $\mathcal{P}$ of $\mathcal{V}$, we can construct $\Phi$ locally by choosing $\Phi_j$ that satisfies \Cref{eqt:localphi_constraint} for each $P_j$. Apparently, the choice of $\Phi_j$ depends on the partition $\mathcal{P}$ and the feasibility of this problem is guaranteed since we can always set $\mathcal{P} = \mathcal{V}$ to fulfill \Cref{eqt:localphi_constraint}. But this choice is not optimal. We should adaptively choose $\mathcal{P}$ and $\Phi$ in such a way that it minimizes $N$, the dimension of $\Phi$.

Suppose we are given the partition $\mathcal{P} = \{ P_j \}_{j=1}^M$ (in other words, the number of patches, $M$, is fixed), minimizing $N$ is equivalent to minimizing the dimension of each $\Phi_j$. In the following, we will first define the notion of {\it interior spectrum} of interior energy $A_{\mathcal{S}}$. \Cref{lemma:int_spectrum} will then show the relationship between the interior spectrum and the minimum dimension that can be achieved for each $\Phi_j$.
\begin{definition}[Interior Spectrum]
Let $\mathcal{S}$ be a subset of $\mathcal{V}$. We define the interior spectrum $\Lambda_{int}(\mathcal{S};A)$ as the set of eigenvalues of $\underline{A}_\mathcal{S}$, where $\underline{A}_\mathcal{S}$ is the {\it interior energy} of $\mathcal{S}$ with respect to $A$ (or we can view it as an operator restricted to the space $\Span\{\mathcal{S}\}$).
\end{definition}

In what follows, since $A$ is generally given and fixed, we will write $\Lambda_{int}(\mathcal{S};A)$ as $\Lambda_{int}(\mathcal{S})$. Also we will simply use $0\leq\lambda_1(\mathcal{S})\leq\lambda_2(\mathcal{S})\leq\cdots\leq\lambda_s(\mathcal{S})$ to denote the ordered elements of $\Lambda_{int}(\mathcal{S})$, where $s=\#\mathcal{S}=\mathrm{dim}(\Span\{\mathcal{S}\})$.

\begin{lemma} Given an $\mathcal{S}\subset \mathcal{V}$, and a constant $\epsilon$, let $q(\epsilon)$ be the smallest integer such that $\frac{1}{\epsilon^2}\leq \lambda_{q(\epsilon)+1}(\mathcal{S})$. Also define $\mathcal{G}(\epsilon)=\{\Theta\subset \Span\{\mathcal{S}\}: \|x-P_\Theta x\|_2\leq\epsilon\|x\|_{\underline{A}_\mathcal{S}},\ \forall x\in \Span\{\mathcal{S}\} \}$, and let $p(\epsilon)=\min_{\Theta\in \mathcal{G}(\epsilon)}\dim \Theta$. Then we have $q(\epsilon)=p(\epsilon)$.
\label{lemma:int_spectrum}
\end{lemma}

By \Cref{lemma:int_spectrum}, one optimal way to minimize $\dim \Phi_j$ for each $P_j$ subject to condition \cref{eqt:localphi_constraint}, is to take $\Phi_j=\Phi_j^{q_j(\epsilon)}$, the eigenspace corresponding to interior eigenvalues $\lambda_1(P_j)\leq \lambda_2(P_j)\leq \cdots\leq\lambda_{q_j(\epsilon)}(P_j)$, where $q_j(\epsilon)$ is the smallest integer such that $\frac{1}{\epsilon^2}\leq \lambda_{q_j(\epsilon)+1}(P_j)$. Recall that this criterion for choosing $\Phi_j$ is based on the fact that the partition $\mathcal{P}$ is given. Then one shall ask a more practical question: how do we construct an ``optimal'' partition $\mathcal{P}$, in the sense that it has a smallest total dimension of $\Phi$?

Instead of answering this question directly, we consider the problem in a more tractable way. We fix an integer $q$, and choose a $q$-dimensional local space $\Phi_j$ for each $P_j$. Then the problem of minimizing $\dim\Phi$ subject to the condition \cref{eqt:localphi_constraint} is reduced to finding a partition $\mathcal{P}=\{P_j\}_{j=1}^M$ with a minimal patch number. Still guided by \Cref{lemma:int_spectrum}, we know that we should choose $\Phi_j=\Phi_j^{q}$, and the condition \cref{eqt:localphi_constraint} is satisfied if and only if $\frac{1}{\epsilon^2}\leq \lambda_{q+1}(P_j)$ for each $P_j$. Define the \textbf{error factor} of a partition $\mathcal{P}$ as
\begin{equation}
\varepsilon(P_j,q) = \frac{1}{\sqrt{\lambda_{q+1}(P_j)}},\quad 1 \leq j \leq M, \quad \text{ and }\quad \varepsilon(\mathcal{P},q)=\max_{j}\frac{1}{\sqrt{\lambda_{q+1}(P_j)}},
\label{eqt:error_factor}
\end{equation}
and so given a constant $\epsilon$, we need to minimize the patch number of $\mathcal{P}$ subject to $\varepsilon(\mathcal{P},q)\leq\epsilon$.

\begin{construction}[Construction of $\Phi$] We choose $\Phi=\bigoplus_{j=1}^M\Phi_j^q$, where $\Phi^q_j\subset \Span\{P_j\}$ is the eigenspace corresponding to the first $q$ interior eigenvalues of patch $P_j$. We also require $(\Phi_j^q)^T\Phi_j^q=I_q$, i.e. $\Phi^T\Phi=I_N$. Then the condition \cref{eqt:localphi_constraint} is satisfied if $\varepsilon(\mathcal{P},q)\leq \epsilon$.
\label{construction:phi}
\end{construction}

We propose \Cref{alg:construct_phi} to construct $\Phi$ guided by the \Cref{construction:phi}. Notice that it also computes the compliment space $U_j$ of $\Phi_j$ in each $\Span(P_j)$, which will serve for the purpose of performing multi-resolution matrix decomposition in \Cref{Sec:Multiresolution}.

\begin{algorithm}[!h]
\caption{\it Construction of $\Phi$}
\label{alg:construct_phi}
\begin{algorithmic}[1]
\REQUIRE{Energy decomposition $\mathcal{E}$, partition $\mathcal{P}$ subject to $\varepsilon(\mathcal{P},q)\leq \epsilon$.}
\FOR{each $P_j\in\mathcal{P}$}
\STATE{Extract $\underline{A}_{P_j}$ from $\mathcal{E}$;}
\STATE{Find the first $q$ normalized eigenvectors of $\underline{A}_{P_j}$ as $\Phi_j$;}
\STATE{Find $U_j$ such that $[\Phi_j,U_j]$ is an orthonormal basis of $\Span(P_j)$;}
\ENDFOR
\STATE{Collect all $\Phi_j$ as $\Phi$, and all $U_j$ as $U$.}
\end{algorithmic}
\end{algorithm}

\begin{remark}
$ $ \linebreak
\vspace{-3mm}
\begin{itemize}
\item[-] The construction of $U_j$ can be implicitly done, for example, by extending $\Phi_j$ to an orthonormal basis of $\Span(P_j)$ with local QR factorization, where only $q$ Householder vectors $[h_1,h_2,\cdots,h_q]$ need to be stored. In fact, we can apply economic QR factorization to $\Phi_j$ to obtain $(I-h_1h_1^T)(I-h_2h_2^T)\cdots(I-h_qh_q^T)=[Q_j,U_j]$ where $[Q_j,U_j]$ is orthogonal and $\mathrm{Span}(\Phi_j)=\mathrm{Span}(Q_j)$. In following algorithms there are only two kinds of operation that involve $U_j$, namely $U_j^Tx$ for some $x\in \mathbb{R}^{s}$ and $U_jx$ for some $x\in \mathbb{R}^{s-q}$. The former one can be done by computing $y=(I-h_qh_q^T)\cdots(I-h_2h_2^T)(I-h_1h_1^T)x$ and then taking the last $s-q$ entries of $y$; the latter one can be done by extending $x$ to $\tilde{x}=[\mathbf{0},x]$ with additional $q$ $0$s in front and then computing $(I-h_1h_1^T)(I-h_2h_2^T)\cdots(I-h_qh_q^T)\tilde{x}$.
\item[-] The integer $q$ is given before the partition is constructed, and the choice of $q$ will be discussed in \Cref{Sec:partition}.
\end{itemize}
\end{remark}

\subsection*{Complexity of \Cref{alg:construct_phi}} For simplicity, we assume that all patches in partition $\mathcal{P}$ have the same patch size $s$. Then number of patches is $\#\mathcal{P}=\frac{n}{s}$. Let $F(s)$ denote the local patch-wise complexity of solving partial eigen problem and extending $\Phi_j$ to $[\Phi_j,U_j]$. Then the complexity of \Cref{alg:construct_phi} is 
\begin{equation}
O(\frac{F(s)}{s}\cdot n).
\label{eqt:complexity_phi}
\end{equation}

\subsection{Choice of $\Psi$}
\label{sec:choose_psi}
Suppose that we have determined the space $\Phi=[\varphi_1,\varphi_2,\cdots,\varphi_N]$, the next step is to find $\Psi=[\psi_1,\psi_2,\cdots,\psi_N]=A^{-1}\Phi T$, namely to determine $T$, so that 
\begin{enumerate}
\item each $\psi_i$ is locally computable, or can be approximated by some $\tilde{\psi}_i$ that is locally computable;
\item the stiffness matrix $A_{\text{st}}=\Psi^TA\Psi$ has relatively small condition number, or the condition number can be bounded by some local information.
\end{enumerate}
Generally each $A^{-1}\phi_i$ is not local (sparse), so it may be impossible to find even one $\psi\in \Span\{A^{-1}\Phi\}$ that is locally computable. A more promising idea is to find $\psi$ that can be well approximated by some $\tilde{\psi}_i$ which is locally computable.
\begin{lemma} 
Assume that $\Psi = [\psi_1,\psi_2,\cdots,\psi_N]$ satisfies $\|x - P^{A}_{\Psi}x \|_A \leq \epsilon \| A x \|_2$ and $\|A_{\text{st}}^{-1}\|_2\leq \|A^{-1}\|_2$, and that $\widetilde{\Psi}=[\tilde{\psi}_1,\tilde{\psi}_2,\cdots,\tilde{\psi}_N]$ satisfies $\|\psi_i-\tilde{\psi}_i\|_A\leq \frac{C\epsilon}{\sqrt{N}},\ 1\leq i\leq N$ for some constant $C$. Then we have
\begin{enumerate}
\item For any $x\in\mathbb{R}^n$, and $b=Ax$, we have
\[\|x-P_{\widetilde{\Psi}}^A x\|_A\leq (1+C\|A^{-1}\|_2)\epsilon\|b\|_2.\]
\item For any $x\in\mathbb{R}^n$, and $b=Ax$, we have
\[\|x-P_{\widetilde{\Psi}}^A x\|_2\leq (1+C\|A^{-1}\|_2)^2\epsilon^2\|b\|_2.\]
\item We have
\[\|A^{-1}-P_{\widetilde{\Psi}}^A A^{-1}\|_2\leq (1+C\|A^{-1}\|_2)^2\epsilon^2.\]
\end{enumerate}
\label{lemma:localization_psi}
\end{lemma}
Guided by \Cref{lemma:localization_psi}, in order to preserve the compression accuracy, we require that each $\psi_i$ be approximated accurately in energy norm $\|\cdot\|_A$ by some $\tilde{\psi}_i$ that is locally computable. To implement this idea, we consider the problem reversely. Suppose we already have some $\tilde{\psi}_i$ that is locally computable, so the construction of $\Psi$ is to find $\psi_i\in A^{-1}(\Phi)$ so that $\|\psi_i-\tilde{\psi}_i\|_A$ is small for each $i$. Since $\tilde{\psi}_i$ is given, minimizing $\|\psi_i-\tilde{\psi}_i\|_A$ can be simply solved by taking $\psi_i=P_\Psi^A\tilde{\psi}_i$. Thanks to the expression $P_\Psi^A=A^{-1}\Phi(\Phi^TA^{-1}\Phi)^{-1}\Phi^T$, we can perform the energy projection $P_\Psi^A$ as long as we know $\Phi$. Therefore we have
\begin{align}
\Psi &= P_\Psi^A\widetilde{\Psi}=A^{-1}\Phi(\Phi^TA^{-1}\Phi)^{-1}\Phi^T\widetilde{\Psi}, \\
\Longrightarrow \Phi^T\Psi &=\Phi^TA^{-1}\Phi(\Phi^TA^{-1}\Phi)^{-1}\Phi^T\widetilde{\Psi}=\Phi^T\widetilde{\Psi}.
\end{align}
Then we shall discuss how to describe the locality of each $\tilde{\psi}_i$. Similar to the locality of $\Phi$, though seems greedy, we can also require that $\tilde{\psi}_i\in \Span\{P_{j_i}\}$ for some $j_i$, and this requirement implies that $\varphi_{i'}^T\tilde{\psi}_i=0$, for all $\varphi_{i'}\notin \Span\{P_{j_i}\}$. Then to determine $\Phi^T\Psi=\Phi^T\widetilde{\Psi}$, we still need to determine $\varphi_{i'}^T\tilde{\psi}_i$ for each $\varphi_{i'}\in \Span\{P_{j_i}\}$. But actually, in the following proof of exponential decay of $\psi_i$, we can see that the value of $\varphi_{i'}^T\tilde{\psi}_i$ for each $\varphi_{i'}\in \Span\{P_{j_i}\}$ does not essentially change the decay property of $\psi_i$. We only need to make sure that $\widetilde{\Psi}$ has the same dimension as $\Phi$. So for simplicity, we require that 
\begin{equation}
\varphi^T_{i'}\tilde{\psi}_i=\delta_{i',i},\quad 1\leq i' \leq N,\quad i.e. \quad \Phi^T\Psi=\Phi^T\widetilde{\Psi}=I_N.
\end{equation}
Adding this extra localization constraint to the form of $\Psi = A^{-1}\Phi T$, we can choose $\Psi$ as follows:
\begin{construction}[Construction of $\Psi$] 
We choose $\Psi=A^{-1}\Phi T$ so that $\Phi^T\Psi=I_N$, that is 
\begin{equation}
\Psi=A^{-1}\Phi(\Phi^TA^{-1}\Phi)^{-1},\quad T=(\Phi^TA^{-1}\Phi)^{-1},
\label{eqt:construction_psi}
\end{equation}
and we have
\begin{equation}
A_{\text{st}}=\Psi^TA\Psi=(\Phi^TA^{-1}\Phi)^{-1}.
\label{eqt:stiffness}
\end{equation}
\label{construction:psi}
\end{construction}
\begin{remark} Our choice of $\Psi$ is inspired by the result proposed by Owhadi in \cite{owhadi2017multigrid}, where the author obtained the same format of $\Psi$ from a marvelous probabilistic perspective. In this work, the idea of {\it Gamblet Transformation} is introduced. Such transformation gives a particular choice of basis in the modified coarse space, which ensures the exponential decay feature of $\Psi$. Our derivation of the choice of $\Psi$ can be seen as a algebraic interpretation of Owhadi's probabilistic construction.
\end{remark}
Though we construct each $\psi_i$ from some local vector $\tilde{\psi}_i$, the error $\|\psi_i-\tilde{\psi}_i\|_A$ is not necessarily small. To have both good locality and small error, we need to use some thing in between. The following lemma (see also Section 3.2 in \cite{owhadi2017multigrid}) shows that the construction of $\Psi$ in \Cref{eqt:construction_psi} is equivalent to the optimizer of a minimization problem. 
\begin{lemma} 
Let $\Psi$ be constructed as in \Cref{eqt:construction_psi}. Then for each $i$, $\psi_i$ satisfies
\begin{equation*}
\begin{array}{rcl}
\psi_i &=& \arg\min_{x\in \mathbb{R}^n} \|x\|_A,\\\\
\text{subject to } \quad \varphi_{i'}^T x &=& \delta_{i',i},\quad \forall\ i' = 1,\ldots, N.
\end{array}
\end{equation*}
\label{lemma:psi_optimization}
\end{lemma}
\begin{proof}
Notice that $A\psi_i\in \Phi$, thus for any $x$ that satisfies $\varphi_{i'}^Tx=\delta_{i',i}$, we have
$\Phi^T(x-\psi_i)=0$, and $\psi_i^TA(x-\psi_i)=0$. Then we have
\begin{equation}
\|x\|_A^2=\|x-\psi_i+\psi_i\|_A^2=\|\psi_i\|_A^2+\|x-\psi_i\|_A^2+2\psi_i^TA(x-\psi_i)\geq\|\psi_i\|_A^2.
\end{equation}
\end{proof}
By the construction given in \Cref{eqt:construction_psi} and guided by \Cref{lemma:psi_optimization}, we can obtain every $\psi_i$ by solving the optimization problem. Our next step is to make use of this minimal property to construct local $\tilde{\psi}_i$ that will be proved exponentially convergent to $\psi_i$.
\begin{definition}[Layers of neighbors]
Let $\mathcal{P}=\{P_j\}_{j=1}^M$ be a partition of $\mathcal{V}$. For any $P_j\in \mathcal{P}$, we recursively define $S_0(P_j)=P_j$, and
\begin{equation}
S_{k+1}(P_j)=\bigcup_{P_{j'}\sim S_k(P_j)}P_{j'},\quad k=0,1,2,\cdots.
\end{equation}
$S_k(P_j)$ is called the $k_\mathrm{th}$ neighbor patch ball of $P_j$, and $S_k(P_j)/S_{k-1}(P_j)$ the $k_\mathrm{th}$ neighbor patch layer of $P_j$.
\label{eqt:layer_neighbors}
\end{definition}

\begin{remark}
By making use of the notion of {\it neighboring} introduced in \Cref{def:neighboring}, we can construct the ``algebraic neighbor layers'' starting from any initial patch $P_j$. Still we do not implicitly assume any underlying physical domain to the operator $A$.
\end{remark}

\begin{definition}[Local approximator] For each $\psi_i$, let $P_{j_i}$ be the patch such that $\varphi_i\in \Phi_{j_i}\subset\Span\{P_{j_i}\}$. Then for each $k\geq0$, we define the {\it $k$-local approximator} of $\psi_i$ as 
\begin{equation}
\begin{array}{rcl}
 \psi_i^k &=& \underset{{x\in \Span\{S_k(P_{j_i})\}}}{\argmin} \|x\|_A,\\\\
\text{subject to } \quad \varphi_{i'}^Tx &=& \delta_{i',i},\quad \forall\ i' = 1,\ldots, N.
\end{array}
\label{eqt:optimization_psi}
\end{equation}
\label{local_psi}
\end{definition}

\begin{remark}
Here $k$ is called the radius of $\psi^k_i$. The condition $\varphi_{i'}^T\psi_i^k=\delta_{i',i}$ is equivalent to $\Phi^T\psi_i^k=\Phi^T\psi_i$. By \Cref{lemma:psi_optimization} and the definition of $\psi_i^k$, we have
\begin{equation}
(\psi_i^k-\psi_i)^TA\psi_i=0,\quad (\psi_i^{k-1}-\psi_i^k)^TA\psi_i^k=0,\quad \forall k,
\end{equation}
and hence
\begin{align}
\|\psi_i^k\|_A^2&=\|\psi_i\|_A^2+\|\psi_i^k-\psi_i\|_A^2,\\
\|\psi_i^{k-1}\|_A^2&=\|\psi_i^k\|_A^2+\|\psi_i^{k-1}-\psi_i^k\|_A^2.
\label{eqt:psi_k_orthogonality}
\end{align}
\end{remark}

\begin{definition}[Condition factor of partition]
Let $\mathcal{P}=\{P_j\}_{j=1}^M$ be a partition of $\mathcal{V}$. Writing $\overline{A}_{P_j}^{-1}$ the inverse of $\overline{A}_{P_j}$ as an operator restricted on $\Span\{P_j\}$, we define
\begin{equation}
\delta(P_j, \Phi_j)=\max_{x\in \Phi_j}\ \frac{x^Tx}{x^T\overline{A}_{P_j}^{-1}x},\quad 1\leq j\leq M,\qquad \text{and} \qquad \delta(\mathcal{P}, \Phi)=\max_{P_j\in\mathcal{P}}\ \delta(P_j,\Phi_j).
\end{equation}
\label{def:delta}
\end{definition}

\begin{remark}
$ $ \linebreak
\vspace{-3mm}
\begin{itemize}
\item[-] In what follows, since we always fix a choice of $\Phi$ for a partition $\mathcal{P}$, we will simply use $\delta(P_j)$ and $\delta(\mathcal{P})$ to denote $\delta(P_j,\Phi_j)$ and $\delta(\mathcal{P},\Phi)$ respectively. In particular, when we use the \Cref{construction:phi} for $\Phi$ with some integer $q$, we correspondingly use the notations $\delta(\mathcal{P},q)$.
\item[-] If we follow the construction $\Phi_j^T\Phi_j=I_{q_j}$, where $q_j$ is the dimension of $\Phi_j$,  then we have
\begin{equation}
\delta(P_j, \Phi_j)=\max_{c\in \mathbb{R}^{q_j}}\ \frac{c^T\Phi_j^T\Phi_jc}{c^T\Phi_j^T\overline{A}_{P_j}^{-1}\Phi_jc}=\max_{c\in \mathbb{R}^{q_j}}\ \frac{c^Tc}{c^T\Phi_j^T\overline{A}_{P_j}^{-1}\Phi_jc}=\|(\Phi_j^T\overline{A}_{P_j}^{-1}\Phi_j)^{-1}\|_2,
\label{eqt:delta&Ast}
\end{equation}
that is 
\begin{equation}
(\Phi_j^T\overline{A}_{P_j}^{-1}\Phi_j)^{-1}\preceq \delta(P_j)I_{q_j}\preceq \delta(\mathcal{P})I_{q_j}.
\end{equation}
Moreover, by block-wise inequalities we have
\begin{equation}
(\Phi^T(\sum_{j=1}^M\overline{A}_{P_j})^{-1}\Phi)^{-1} \preceq \delta(\mathcal{P})I_N.
\label{eqt:delta&Ast2}
\end{equation}
This analysis will help us to bound the maximum eigenvalue of the stiffness matrix $A_{\text{st}}=\Psi^TA\Psi$ by $\delta(\mathcal{P})$.
\end{itemize}
\end{remark}

\begin{example} 
In this example, we consider the operator $A$ to be the discretization of 2-D second-order elliptic operator by standard 5-point Finite Difference scheme. Similar to the case of graph laplacian in \Cref{example:Example2}, we have a natural energy decomposition inherited from the assembling of such discretization. Specifically, for every pair of vertices $e_{hori} := [(i,j),(i,j+1)]$ and $e_{vert} := [(i,j),(i+1,j)]$ in the finite difference grid, the energy elements are
\begin{equation*}
E^{ij}_{hori} = -\frac{1}{|e_{hori}|^2} 
\begin{scriptsize}
\begin{blockarray}{cccccc}
& (i,j) & & (i,j+1) & \\
\begin{block}{(ccccc)c}
  \textbf{\large $ \ 0 $} &  &  &  &  &   \\
   & a_{i,j+\frac{1}{2}} &  & -a_{i,j+\frac{1}{2}} &  & (i,j) \\
   &  & \ddots &  &  &   \\
   & -a_{i,j+\frac{1}{2}} &  & a_{i,j+\frac{1}{2}} &  & (i,j+1) \\
   &  &  &  & \textbf{\large $ 0 \ $} & \\
\end{block}
\end{blockarray}
\end{scriptsize}, \text{ and }
\end{equation*}
\begin{equation*}
E^{ij}_{vert} = -\frac{1}{|e_{vert}|^2}
\begin{scriptsize}
\begin{blockarray}{cccccc}
& (i,j) & & (i+1,j) & \\
\begin{block}{(ccccc)c}
  \textbf{\large $ \ 0 $} &  &  &  &  &   \\
   & a_{i+\frac{1}{2},j} &  & -a_{i+\frac{1}{2},j} &  & (i,j) \\
   &  & \ddots &  &  &   \\
   & -a_{i+\frac{1}{2},j} &  & a_{i+\frac{1}{2},j} &  & (i+1,j) \\
   &  &  &  & \textbf{\large $ 0 \ $} & \\
\end{block}
\end{blockarray}
\end{scriptsize}.
\end{equation*}
Now suppose we are given a partition $\mathcal{P}$, we focus on a particular local patch $P_j$ to study how mesh size and contrast affect the {\bf error factor} and the {\bf condition factor}. \Cref{fig:example3a} shows the high-contrast field (colored in black) in $P_j$. For simplicity, we set $P_j$ as a square domain with fixed length $H = |e_{hori}| = |e_{vert}|$. We also set the coefficient in high-contrast field to be $10^3$ (and 1 otherwise). \Cref{fig:example3b} shows the decreasing trend of the {\bf error factor} $\epsilon(P_j,1)$ as the vertex number $\#V$ inside the patches (i.e. the vertex density in $P_j$) increases. Here we choose $q = 1$ for illustration purpose. \Cref{fig:example3c} shows a similar decreasing trend of the {\bf condition factor} $\delta(P_j,1)$ and \Cref{fig:example3d} plots $\epsilon(P_j,1)^2 \cdot \delta(P_j,1)$ versus $\#V$.  Fixing the vertex number $\#V$ in the patch $P_j$, we also study the relationship of $\epsilon(P_j,1)$, $\delta(P_j,1)$ and the contrast. In particular, we double the contrast by 2 in each single computation and investigate the trend of $\epsilon(P_j,1)$ and $\delta(P_j,1)$. \Cref{fig:example3e} shows the decrease of $\epsilon(P_j,1)$ as contrast increases. For $\delta(P_j,1)$, although it also increases as the contrast increases, we can clearly see that there is an upper bound (around 220 in this example), even when the contrast jumps up to $2^{20}$. \Cref{fig:example3g} plots $\epsilon(P_j,1)^2 \cdot \delta(P_j,1)$ versus contrast. 

\begin{figure}[h!]
\centering
\begin{subfigure}[b]{0.2\textwidth}
        \centering
        \includegraphics[width=1.0\textwidth]{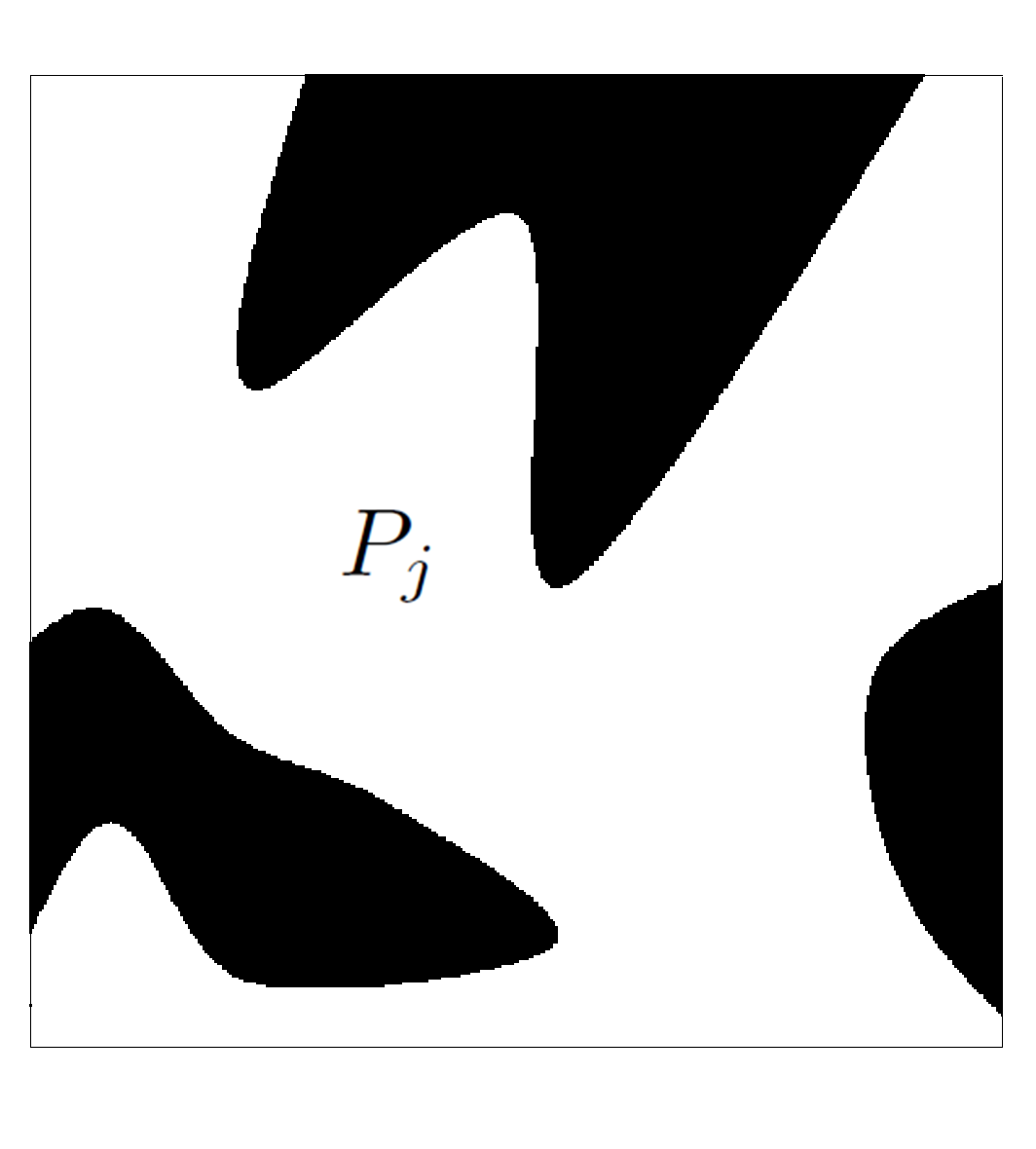}
        \caption{}
        \label{fig:example3a}
    \end{subfigure}%
    \begin{subfigure}[b]{0.24\textwidth}
        \centering
        \includegraphics[width=1.0\textwidth]{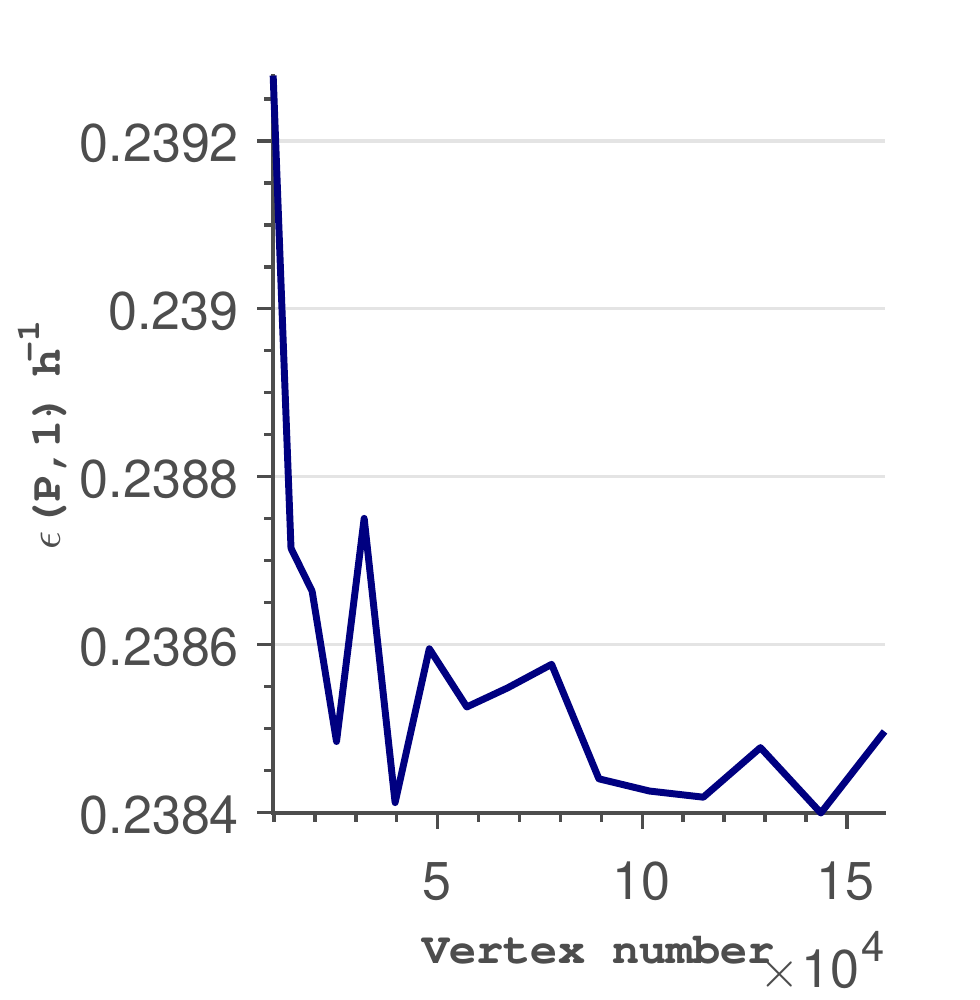}
        \caption{}
        \label{fig:example3b}
    \end{subfigure}
    \begin{subfigure}[b]{0.24\textwidth}
        \centering
        \includegraphics[width=1.0\textwidth]{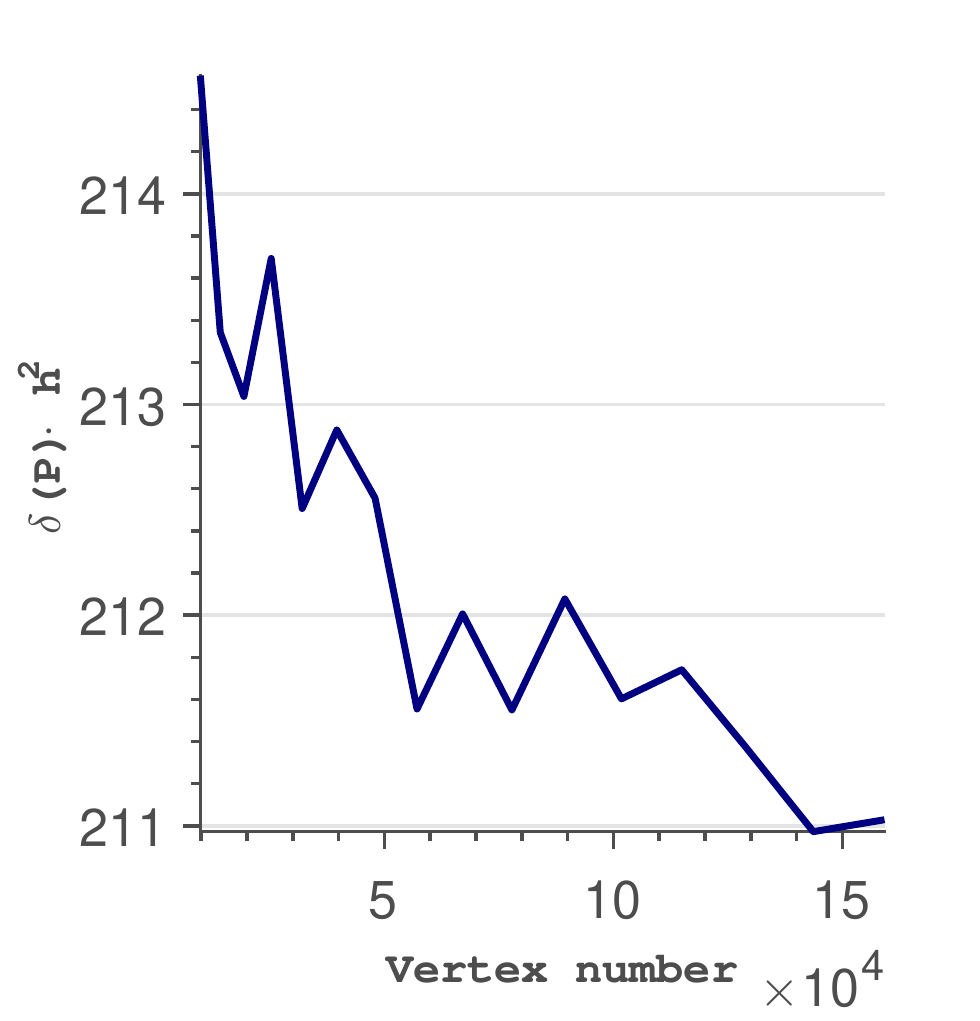}
        \caption{}
        \label{fig:example3c}
    \end{subfigure}
    \begin{subfigure}[b]{0.24\textwidth}
        \centering
        \includegraphics[width=1.0\textwidth]{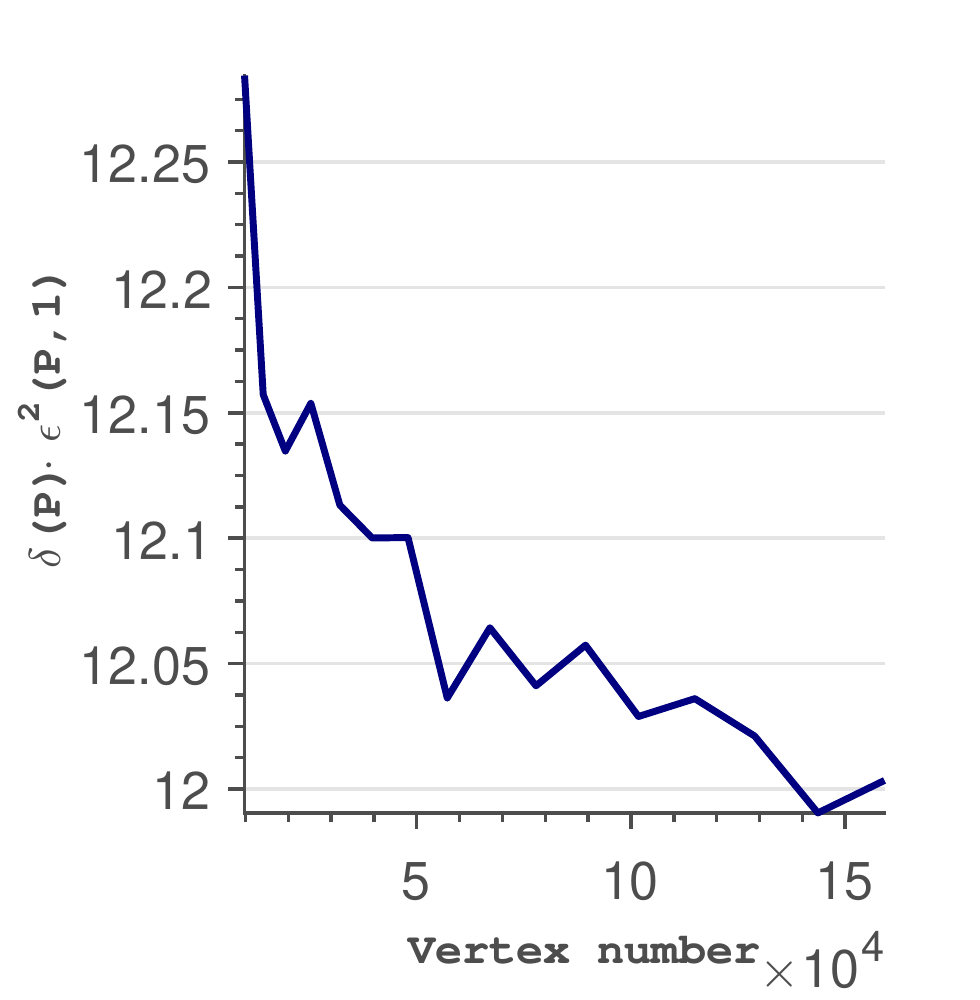}
        \caption{}
        \label{fig:example3d}
    \end{subfigure}
    
    \begin{subfigure}[b]{0.24\textwidth}
        \centering
        \includegraphics[width=1.0\textwidth]{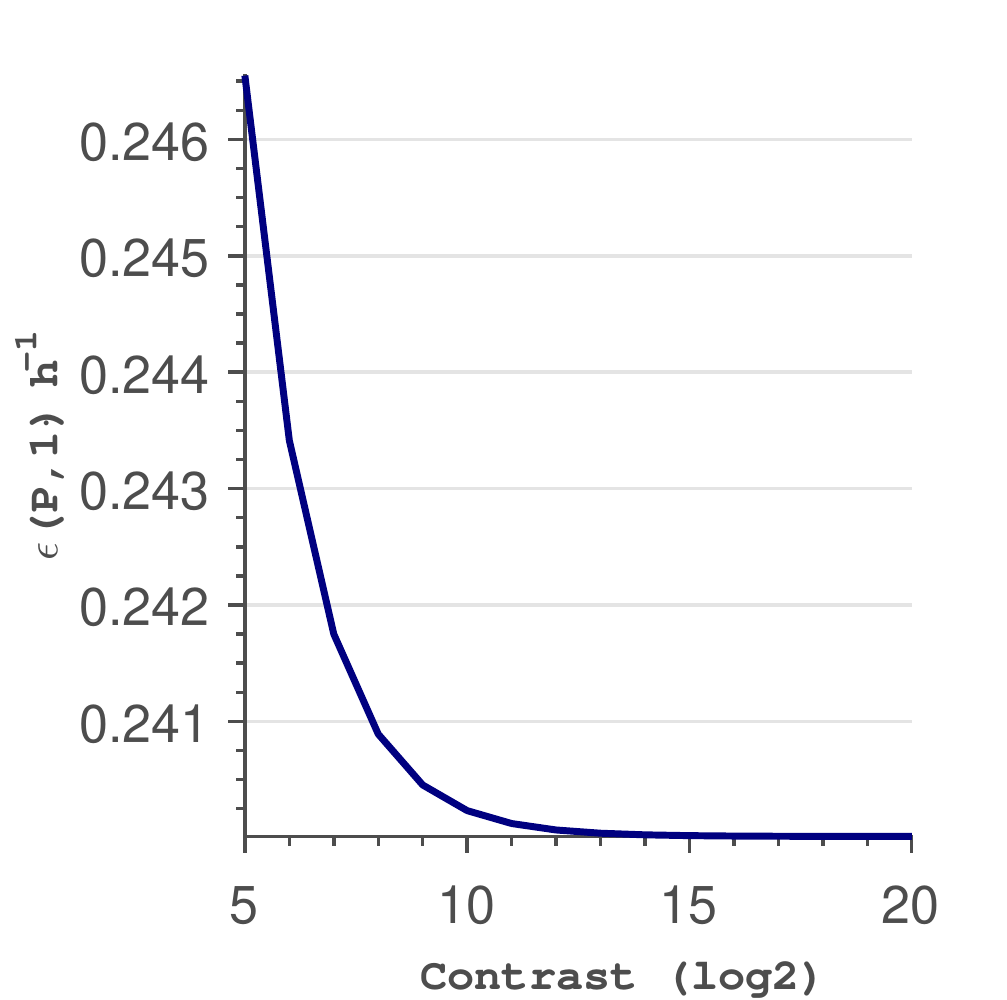}
        \caption{}
        \label{fig:example3e}
    \end{subfigure}
    \begin{subfigure}[b]{0.24\textwidth}
        \centering
        \includegraphics[width=1.0\textwidth]{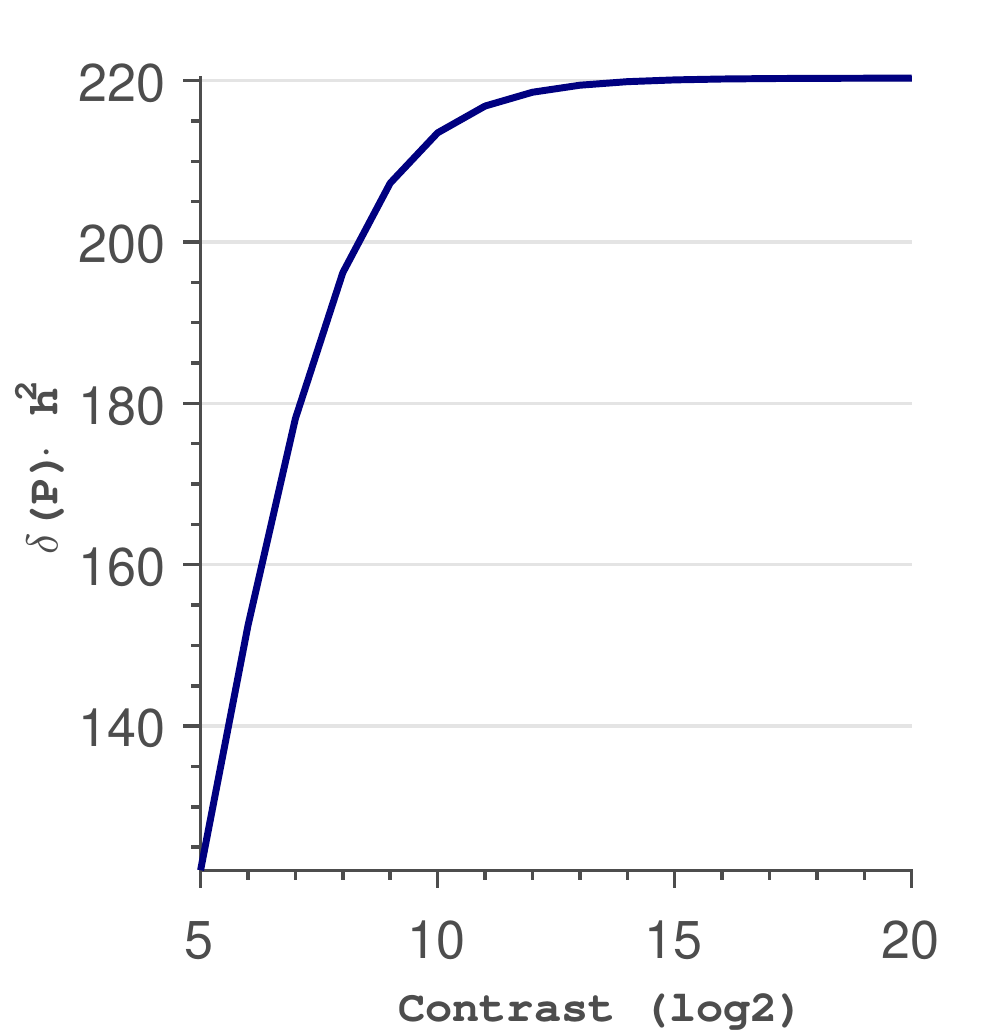}
        \caption{}
        \label{fig:example3f}
    \end{subfigure}
    \begin{subfigure}[b]{0.24\textwidth}
        \centering
        \includegraphics[width=1.0\textwidth]{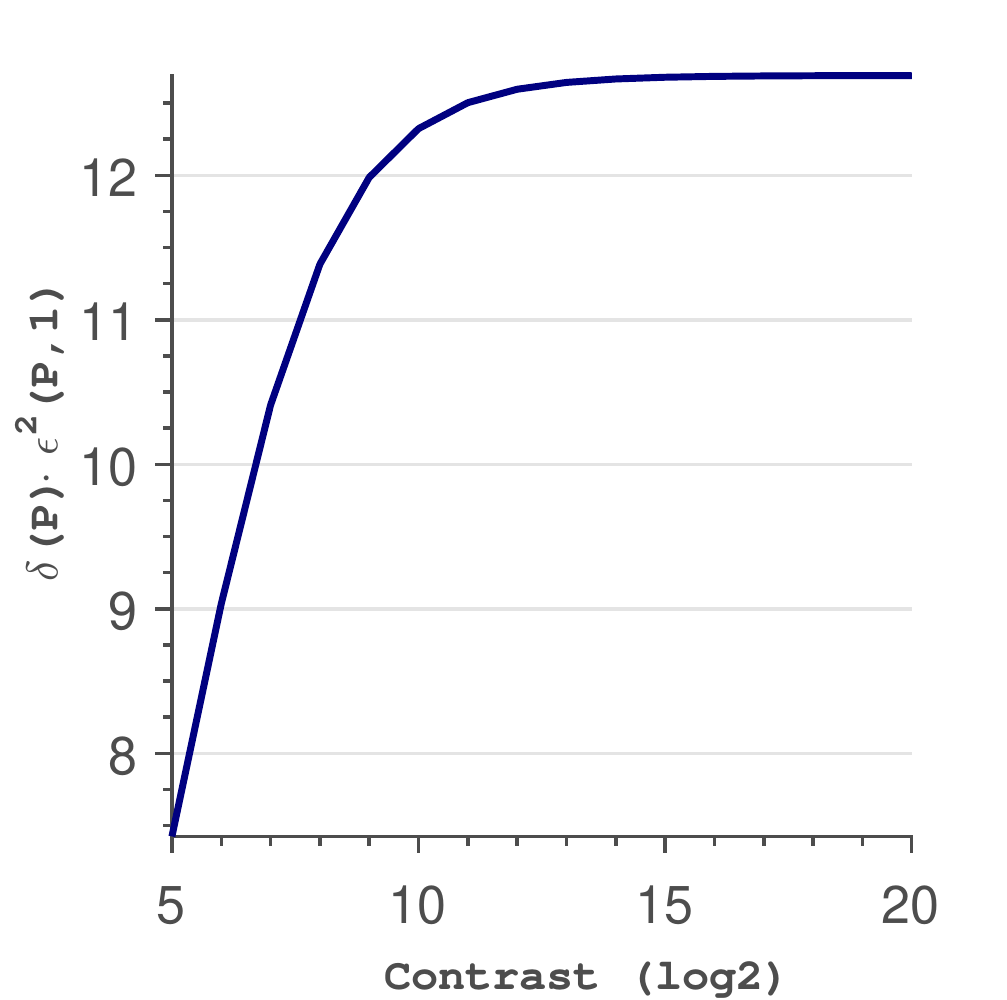}
        \caption{}
        \label{fig:example3g}
    \end{subfigure}
    \caption{Example showing the relationship between mesh size, the {\bf error factor} $\epsilon(P_j,1)$, the {\bf condition factor} $\delta(P_j,1)$ and contrast.}
    \label{fig:example3}
\end{figure}
\label{example:example3}
\end{example}

The following theorem shows the scaling properties of $\psi_i$, $\psi_i^k$ under construction \Cref{construction:psi} and definition \Cref{local_psi}, which will help to prove the exponential decay of the basis function $\psi_i$.

\begin{thm} 
For each $\psi_i$, we have 
\begin{equation}
\|\psi_i\|_A\leq \|\psi_i^k\|_A\leq\|\psi_i^0\|_A\leq\sqrt{\delta(P_{j_i})}.
\end{equation}
\label{thm:psi_delta}
\end{thm}
\subsubsection*{Compliment space} 
For each $P_j$, without causing any ambiguity, we use $U_j$ interchangeably to denote both the orthogonal compliment of $\Phi_j$ with respect to $\Span\{P_j\}$, or an orthonormal basis matrix of $U_j$. Namely we have $U_j\subset\Span\{P_j\}$, and $\Phi_j^TU_j= {\bf 0}$. Then we define 
\begin{equation}
\alpha(P_j)=\max_{x\in U_j}\ \frac{x^T\overline{A}_{P_j}x}{x^T\underline{A}_{P_j}x},\quad 1\leq j\leq M,\qquad \text{and}\qquad \alpha(\mathcal{P})=\max_{P_j\in \mathcal{P}}\ \alpha(P_j).
\label{def:alpha}
\end{equation}
\begin{remark}
$ $ \linebreak
\vspace{-3mm}
\begin{itemize}
\item[-] If we choose $\Phi_j$ so that it satisfies condition \cref{eqt:localphi_constraint}, then we have
\[x^T\underline{A}_{P_j}x=\|x\|^2_{\underline{A}_{P_j}}\geq\frac{1}{\epsilon^2}\|x-P_{\Phi_j}x\|_2^2=\frac{1}{\epsilon^2}\|x\|_2^2,\quad \forall x\in U_j,\]
and $x^T\overline{A}_{P_j}x\leq \|\overline{A}_{P_j}\|_2\|x\|_2^2$. Thus $\alpha(P_j)\leq \epsilon^2\|\overline{A}_{P_j}\|_2$. This argument is meant to show that we can have $\alpha(\mathcal{P})<+\infty$ if we choose $\mathcal{P}$ and $\Phi$ properly. But this bound is not tight, as $\alpha(\mathcal{P})$ can be much smaller in general.
\item[-] An immediate result of the definition of $\alpha(P_j)$ is that 
\[U_j^T\overline{A}_{P_j}U_j\preceq \alpha(P_j)U_j^T\underline{A}_{P_j}U_j,\quad \forall\ j.\]
\end{itemize}
\end{remark}

The following theorem shows that the local basis function $\tilde{\psi}_i^k$ is exponentially convergent to $\psi_i$ as its support $S_k(P_{j_i})$ extends(or as $k$ increases). Indeed, the exponential decay of $\psi_i$ has been proved in \cite{maalqvist2014localization,owhadi2017multigrid,owhadi2017universal,hou2016sparse} in different manners based on a common observation that the energy of $\psi_i$ in the region beyond a certain single layer of patches is comparable to its energy only on this layer, which reflects the local interacting feature of the operator $A$ itself. Also based on this observation, we modify the proof in Section 6 of \cite{owhadi2017universal} using matrix framework coherent to our energy settings.

\begin{thm}[Exponential decay]
For each $\psi_i$, we have
\begin{equation}
\|\psi_i^k-\psi_i\|_A^2\leq \Big(\frac{\alpha(\mathcal{P})-1}{\alpha(\mathcal{P})}\Big)^k\|\psi_i^0-\psi_i\|_A^2\leq \Big(\frac{\alpha(\mathcal{P})-1}{\alpha(\mathcal{P})}\Big)^k\delta(P_{j_i}).
\end{equation}
\label{thm:exponential_decay}
\end{thm}
\begin{proof}
For simplicity, we will write $\psi_i$ as $\psi$, $\psi_i^k$ as $\psi^k$, $P_{j_i}$ as $P$, and $S_k(P_{j_i})$ as $S_k$. Let $Y_k$ denote the joint space of all $U_j$ such that $P_j\subset S_k$, and $Z_k$ the joint space of all $U_j$ such that $P_j\subset \mathcal{V}\backslash S_k$. We still use $Y_k, Z_k$ as the basis matrix for the spaces $Y_k, Z_k$, so that each $U_j$ is a bunch of columns of either $Y_k$ or $Z_k$. We use $U$ to denote $Y_\infty$. Notice that we can always arrange $U_j$ in a particular order so that the matrix form $U=[Y_k,Z_k]$ holds. We define 
\begin{equation}
r^k=\psi^k-\psi,\quad k\geq0;\qquad w^k=\psi^{k-1}-\psi^k,\quad k\geq1,
\end{equation}
then according to \Cref{eqt:psi_k_orthogonality} we have
\begin{equation}
\|r^{k-1}\|_A^2=\|r^k\|_A^2+\|w^k\|_A^2.
\end{equation}
Since $\Phi^T(\psi^k-\psi)=\Phi^T(\psi^{k-1}-\psi^k)= {\bf 0}$, we have $r^k\in U$ and $w^k\in Y_k$. Then by the minimal properties of $\psi^k$ and $\psi$, we actually have
\begin{equation}
r^{k-1}=P_{U}^A\psi^{k-1}=U(U^TAU)^{-1}U^TA\psi^{k-1},
\end{equation}
\begin{equation}
w^k=P_{Y_k}^A\psi^{k-1}=Y_k(Y_k^TAY_k)^{-1}Y_k^TA\psi^{k-1}=Y_k(Y_k^TA_{S_k}Y_k)^{-1}Y_k^TA_{S_k}\psi^{k-1}.
\end{equation}
By the definition of $S_k$, we know that $S_{k-1}\not\sim\mathcal{V}\backslash S_k$, and therefore $Z_k^TA\psi^{k-1}=0$. Then we get
\begin{align*}
\|r^{k-1}\|_A^2=&\ \psi^{k-1,T}AU(U^TAU)^{-1}U^TA\psi^{k-1}\\
=&\ \psi^{k-1,T}A\left[\begin{array}{cc}
Y_k& 0
\end{array}\right](U^TAU)^{-1}\left[\begin{array}{c}
Y_k^T\\
0
\end{array}\right]
A\psi^{k-1}.
\end{align*}
Due to the locality of $\overline{A}_{P_j},\underline{A}_{P_j}$, we obtain
\begin{align*}
U^TAU\succeq U^T\Big(\sum_{j=1}^M\underline{A}_{P_j}\Big)U&=\sum_{j=1}^M\left(\begin{smallmatrix} \textbf{\large $ \ 0 $} & & \\ & U_j^T\underline{A}_{P_j}U_j & \\ & & \textbf{\large $ \ 0 $}\end{smallmatrix}\right) \succeq \frac{1}{\alpha(\mathcal{P})}\sum_{j=1}^M\left(\begin{smallmatrix} \textbf{\large $ \ 0 $} & & \\ & U_j^T\overline{A}_{P_j}U_j & \\ & & \textbf{\large $ \ 0 $} \end{smallmatrix}\right).
\end{align*}
As a simple inference of \Cref{prop:OED}, we have
\[\sum_{P_j\subset S_k}\overline{A}_{P_j}\succeq A_{S_k},\quad \sum_{P_j\subset \mathcal{V}\backslash S_k}\overline{A}_{P_j}\succeq A_{\mathcal{V}\backslash S_k},\]
and therefore
\begin{equation*}
\sum_{j=1}^M\left(\begin{smallmatrix} \textbf{\large $ \ 0 $} & & \\ & U_j^T\overline{A}_{P_j}U_j & \\ & & \textbf{\large $ \ 0 $} \end{smallmatrix}\right) = \ \left(\begin{smallmatrix}Y_k^T\Big(\sum_{P_j\subset S_k}\overline{A}_{P_j}\Big)Y_k & \\ &Z_k^T\Big(\sum_{P_j\subset \mathcal{V}\backslash S_k}\overline{A}_{P_j}\Big)Z_k\end{smallmatrix}\right)
\succeq \ \left(\begin{smallmatrix}Y_k^TA_{S_k}Y_k & \\ &Z_k^TA_{\mathcal{V}\backslash S_k}Z_k\end{smallmatrix}\right).
\end{equation*}
Combining all results, we have
\begin{equation*}
(U^TAU)^{-1} \preceq \alpha(\mathcal{P})\left(\begin{array}{cc}Y_k^TA_{S_k}Y_k & \\ &Z_k^TA_{\mathcal{V}\backslash S_k}Z_k\end{array}\right)^{-1} = \alpha(\mathcal{P})\left(\begin{array}{cc}(Y_k^TA_{S_k}Y_k)^{-1} & \\ & (Z_k^TA_{\mathcal{V}\backslash S_k}Z_k)^{-1}\end{array}\right),
\end{equation*}
and thus 
\begin{align*}\|r^{k-1}\|_A^2\leq&\ \alpha(\mathcal{P})\psi^{k-1,T}A\left[\begin{array}{cc}
Y_k& 0
\end{array}\right]\left(\begin{smallmatrix}(Y_k^TA_{S_k}Y_k)^{-1} & \\ & (Z_k^TA_{\mathcal{V}\backslash S_k}Z_k)^{-1}\end{smallmatrix}\right)\left[\begin{array}{c}
Y_k^T\\
0
\end{array}\right]
A\psi^{k-1}\\
=&\ \alpha(\mathcal{P})\psi^{k-1,T}AY_k(Y_k^TA_{S_k}Y_k)^{-1}Y_k^TA\psi^{k-1}\\
=&\ \alpha(\mathcal{P})\psi^{k-1,T}A_{S_k}Y_k(Y_k^TA_{S_k}Y_k)^{-1}Y_k^TA_{S_k}\psi^{k-1}\\
=&\ \alpha(\mathcal{P})\|w^k\|_A^2.
\end{align*}
This gives us 
\begin{equation*}
\|r^{k-1}\|_A^2=\|r^k\|_A^2+\|w^k\|_A^2\geq \|r^k\|_A^2 +\frac{1}{\alpha(\mathcal{P})}\|r^{k-1}\|_A^2 \ \Longrightarrow \ \|r^k\|_A^2\leq \frac{\alpha(\mathcal{P})-1}{\alpha(\mathcal{P})}\|r^{k-1}\|_A^2.
\end{equation*}
Applying this recursively, we have
\begin{equation}
\|\psi^k-\psi\|_A^2\leq \Big(\frac{\alpha(\mathcal{P})-1}{\alpha(\mathcal{P})}\Big)^k\|\psi^0-\psi\|_A^2.
\end{equation}
Notice that $\|\psi^0-\psi\|_A^2=\|\psi^0\|_A^2-\|\psi\|_A^2\leq \|\psi^0\|_A^2\leq \delta(P)$, and this completes our proof.
\end{proof}

\begin{remark}
\label{rmk:alpha_remark}
$ $ \\
\vspace{-3mm}
\begin{itemize}
\item[-] Recall that in \Cref{lemma:localization_psi}, to make a k-layer approximator $\Psi^k$ become a good approximator, we need $\|\psi_i-\psi_i^k\|_A\leq \frac{C\epsilon}{\sqrt{N}}$ for some constant $C$, and thus \Cref{thm:exponential_decay} guides us to choose 
\begin{equation}
k= O(\log\frac{1}{\epsilon}+\log N+\log \delta(\mathcal{P})).
\label{eqt:radius}
\end{equation}
And the locality of $\psi_i^k$ lies in the local connection property of the matrix $A$. 
\item[-] By the definition in \Cref{def:alpha}, $\alpha(\mathcal{P})$ is locally scaling invariant. Therefore the layer-wise decay rate is unchanged when $A$ is locally multiplied by some scaling constant.
\end{itemize}
\end{remark}

\begin{corollary} For any $\psi_i$, the interior energy of $\psi_i$ on $\mathcal{V}/S_k(P_{j_i})=S_k^c(P_{j_i})$ decays exponentially with $k$, namely
\[\|\psi_i\|_{\underline{A}_{S_k^c(P_{j_i})}}^2\leq \Big(\frac{\alpha(\mathcal{P})-1}{\alpha(\mathcal{P})}\Big)^k\delta(P_{j_i}).\]
Moreover, for any two $\psi_i,\psi_{i'}$, we have
\[|\psi_i^TA\psi_{i'}|\leq \Big(\frac{\alpha(\mathcal{P})-1}{\alpha(\mathcal{P})}\Big)^{\frac{k_{ii'}}{4}-\frac{1}{2}}\delta(\mathcal{P}),\]
where $k_{ii'}$ is the largest integer such that $P_{j_{i'}}\subset S_{k_{ii'}}^c(P_{j_i})$( or equivalently $P_{j_i}\subset S_{k_{ii'}}^c(P_{j_{i'}})$).
\label{cor:expdecay}
\end{corollary}

\begin{remark}
Recall that in \Cref{lemma:localization_psi} for the compression error with localization $\tilde{\epsilon}_{\text{com}}^2$ to be bounded by some prescribed accuracy $\epsilon^2$, we need the localization error $\epsilon^2_{loc} \leq \frac{\epsilon^2}{N}$. But now with the exponential decaying feature of $\Psi$, empirically, we observe that we can relax the requirement of the localization error to be $\epsilon^2_{loc} \leq O(\epsilon^2)$ in practice.
\end{remark}

Now we have constructed a $\Psi$ that can be approximated by local computable basis in energy norm. So the remaining task is to tackle with the second criteria: to give a control on the condition number of $A_{\text{st}}=\Psi^TA\Psi$.

\begin{thm} Let $\lambda_{\min}(A_{\text{st}})$ and $\lambda_{\max}(A_{\text{st}})$ denote the smallest and largest eigenvalues of $A_{\text{st}}$ respectively, then we have
\begin{equation}
\lambda_{\min}(A_{\text{st}})\geq \lambda_{\min}(A),\qquad \lambda_{\max}(A_{\text{st}})\leq \delta(\mathcal{P}),
\label{eqt:lambdamin_control} 
\end{equation}\label{eqt:lambdamax_control}
and so we have
\begin{equation}\kappa(A_{\text{st}})=\frac{\lambda_{\max}(A_{\text{st}})}{\lambda_{\min}(A_{\text{st}})}\leq \delta(\mathcal{P})\|A^{-1}\|_2.
\end{equation}
\label{thm:condition_number}
\end{thm}

\begin{corollary} Let $\widetilde{\Psi}$ be the local approximator of $\Psi$ defined in \Cref{lemma:localization_psi} such that $\|\psi_{i}-\tilde{\psi}_i\|_A\leq \frac{\epsilon}{\sqrt{N}}$, then we have
\begin{equation}\lambda_{\min}(\widetilde{A}_{\text{st}})\geq \lambda_{\min}(A),\qquad \lambda_{\max}(\widetilde{A}_{\text{st}})\leq \left(1+\frac{\epsilon}{\sqrt{\delta(\mathcal{P})}}\right)^2\delta(\mathcal{P}),
\end{equation}
and thus 
\begin{equation}\kappa(\widetilde{A}_{\text{st}})=\frac{\lambda_{\max}(\widetilde{A}_{\text{st}})}{\lambda_{\min}(\widetilde{A}_{\text{st}})}\leq \left(1+\frac{\epsilon}{\sqrt{\delta(\mathcal{P})}}\right)^2\delta(\mathcal{P})\|A^{-1}\|_2.
\end{equation}
\label{cor:localpsi_conditionnum}
\end{corollary}

Guided by \Cref{thm:condition_number}, we obtain a simple methodology on the control of condition number of the stiffness matrix $A_{\text{st}}$. As shown in \Cref{eqt:lambdamin_control} and \Cref{eqt:lambdamax_control}, the only variable is the choice of partition $\mathcal{P}$ and thus the burden again falls to the construction of the partition $\mathcal{P}$. Nevertheless, this new criterion allows us to regulate the quality of partitions directly by avoiding large $\delta(\mathcal{P})$.

Now we can design an algorithm to construct the local approximator $\widetilde{\Psi}$ of $\Psi$ subject to a desired localization error $\epsilon_{loc}$. Intuitively, a straightforward way is to choose a large enough uniform decay radius $r$ and directly compute $\widetilde{\Psi}=\Psi^r$. The localization error can then be guaranteed by \Cref{lemma:localization_psi}. But redundant computation will probably occur since some $\psi$ may decay much faster than the others. Instead, we propose to compute each $\tilde{\psi}_i$ hierarchically from the center patch $P_{j_i}$ by making use of optimization property \cref{eqt:optimization_psi}. Suppose that we already obtain $\psi_i^{k-1}$, then by optimization property \cref{eqt:optimization_psi}, one can check that $w_i^{k}=\psi_i^{k-1}-\psi_i^{k}$ satisfies the following optimization problem
\begin{equation}
\begin{array}{rcl}
 w_i^k &=& \underset{{w\in \Span\{S_k(P_{j_i})\}}}{\argmin} \|\psi_i^{k-1}-w\|_A,\\\\
\text{subject to } \quad \Phi^Tw &=& 0.
\end{array}
\label{eqt:optimization_w}
\end{equation}
Similar to the proof of \Cref{thm:exponential_decay}, let $Y_{i,k}$ denote the joint space of all $U_j$ such that $P_j\subset \{S_k(P_{j_i})\}$. Then the constraints in optimization problem \cref{eqt:optimization_w} imply that $w_i^k\in \Span\{Y_{i,k}\}$. Therefore we can explicitly compute $w_i^k$ as 
\begin{equation}
w_i^k=P_{Y_{i,k}}^A\psi_i^{k-1}=Y_{i,k}\big(Y_{i,k}^TA_{S_k(P_{j_i})}Y_{i,k}\big)^{-1}Y_{i,k}^TA_{S_k(P_{j_i})}\psi_i^{k-1},
\end{equation}
and thus 
\begin{equation}
\psi_i^k=\psi_i^{k-1}-Y_{i,k}\big(Y_{i,k}^TA_{S_k(P_{j_i})}Y_{i,k}\big)^{-1}Y_{i,k}^TA_{S_k(P_{j_i})}\psi_i^{k-1}.
\label{eqt:compute_psi}
\end{equation}
Specially, we can compute $\psi_i^0$ by \Cref{eqt:compute_psi} with an initial guess $\psi_i^{0-1}\in \Span \{P_{j_i}\}$ satisfying $\varphi_{i'}^T\psi_i^{0-1} = \delta_{i',i},\ \forall\ i' = 1,\ldots, N$. Notice that the main cost of computation of $\psi_i^k$ involves inverting the matrix $Y_{i,k}^TAY_{i,k}$, whose condition number can be bounded by $\varepsilon(\mathcal{P},q)^2\lambda_{\max}(A)\kappa(Y_{i,k}^TY_{i,k})$ as we will see in \cref{lemma:Bst_eigenvalue}. By choosing all $U_j$ orthonormal, we have $\kappa(Y_{i,k}^TY_{i,k})=1$, then the computation efficiency is measured by $\varepsilon(\mathcal{P},q)^2\lambda_{\max}(A)$ if we use CG type method. When we prescribe some certain accuracy $\varepsilon(\mathcal{P},q)^2$ but $\lambda_{\max}(A)$ is really large, a multiresolution strategy will be adopted to ensure the efficiency of computing $\psi_i^k$. 

To summarize, the process of computing a sufficient approximator $\tilde{\psi}_i$ starts with the formation of $\psi_i^0$, then inductively computes $\psi_i^{k}$ by solving inverse problem \cref{eqt:compute_psi} with initializer $\psi_i^{k-1}$, and finally ends with $\tilde{\psi}_i=\psi_i^{r}$ when some stopping criterion is attained for $k=r$. Such inductive computation suggests us to use the CG method to take advantage of the exponential convergence of $\psi_i^k$. Having faith in the exponential decay of $\|\psi_i^k-\psi_i\|_A$, we choose the stopping criterion as
\[\frac{\eta_k^2}{1-\eta_k^2}\|\psi_i^{k-1}-\psi_i^k\|_A^2\leq \epsilon_{loc}^2 ,\quad \mathrm{for}\quad \eta_k=\frac{\|\psi_i^{k-1}-\psi_i^k\|_A}{\|\psi_i^{k-2}-\psi_i^{k-1}\|_A}.\]
The reason is that if $\|\psi_i^k-\psi_i\|_A$ does decay as $\|\psi_i^k-\psi_i\|_A=O(\eta^k)$ for some constant $\eta\in(0,1)$, then 
\[\|\psi_i^k-\psi_i\|_A^2=\frac{\eta^2}{1-\eta^2}O(\|\psi_i^{k-1}-\psi_i\|_A^2-\|\psi_i^k-\psi_i\|_A^2)=\frac{\eta^2}{1-\eta^2}O(\|\psi_i^{k-1}-\psi_i^k\|_A^2),\]
where we have used \Cref{eqt:psi_k_orthogonality}. With the analysis above, we propose \Cref{alg:construct_tilde_psi} for constructing $\widetilde{\Psi}$.

\subsection*{Complexity of \Cref{alg:construct_tilde_psi}} For simplicity, we assume that all patches in partition $\mathcal{P}$ have the same patch size $s$. Let $r$ be the necessary number of layers for $\|\psi_i^r-\psi_i\|_A\leq \frac{\epsilon}{\sqrt{N}}$, where $N$ is the dimension of $\Phi$(or $\Psi$), then we have 
\begin{equation}
r = O(\log\frac{1}{\epsilon}+\log N+\log\delta(\mathcal{P})).
\end{equation} 
Since we are actually compressing $A^{-1}$, we can bound $N$ by the original dimension $n$. Further we assume that locality condition \cref{condition:localityA}, \cref{condition:scalesproperty} and \cref{condition:localregularity} are satisfied, then the support size of each $\psi_i^r$ is $O(s\cdot r^d)$. Since we also only need to solve \Cref{eqt:compute_psi} up to the same relative accuracy $O(\frac{\epsilon}{\sqrt{N}})$ using the CG method, the cost of computing $\psi_i^r$ can be estimated by 
\begin{align}
O(\kappa(Y_{i,k}^TAY_{i,k})\cdot s\cdot r^d \cdot &(\log\frac{1}{\epsilon}+\log N+\log\delta(\mathcal{P})))\nonumber\\
\leq&\ O(\varepsilon(\mathcal{P},q)^2\cdot \lambda_{\max}(A)\cdot s \cdot (\log\frac{1}{\epsilon}+\log n+\log\delta(\mathcal{P}))^{d+1} )
\end{align}
Finally the total complexity of \Cref{alg:construct_tilde_psi} is $N$ times the cost for every $\psi_i^r$, i.e.
\begin{equation}
O(n\cdot q\cdot \varepsilon(\mathcal{P},q)^2\cdot\lambda_{\max}(A)\cdot(\log\frac{1}{\epsilon}+\log n+\log\delta(\mathcal{P}))^{d+1} ),
\label{eqt:complexity_psitilde}
\end{equation}
where we have used the relation $N=nq/s$.
 
\begin{algorithm}[!h]
\caption{\it Construction of $\widetilde{\Psi}$}
\label{alg:construct_tilde_psi}
\begin{algorithmic}[1]
\REQUIRE{Energy decomposition $\mathcal{E}$, partition $\mathcal{P}$, $\Phi$, desired accuracy $\epsilon_{loc}$}
\ENSURE{$\widetilde{\Psi}$}
\FOR{$i=1,2,\cdots,\dim(\Phi)$}
\STATE{Compute $\psi_i^0$ by solving \Cref{eqt:compute_psi} on $\mathcal{S}_0(P_{j_i})$;}
\STATE{Compute $\psi_i^1$ by solving \Cref{eqt:compute_psi} on $\mathcal{S}_1(P_{j_i})$ with initializer $\psi_i^{0}$;}
\REPEAT
\STATE{Compute $\psi_i^k$ by solving \Cref{eqt:compute_psi} on $\mathcal{S}_k(P_{j_i})$ with initializer $\psi_i^{k-1}$;}
\STATE{$\eta \leftarrow \frac{\|\psi_i^k-\psi_i^{k-1}\|_A}{\|\psi_i^{k-1}-\psi_i^{k-2}\|_A}$;}
\UNTIL{$\frac{\eta^2}{1-\eta^2}\|\psi_i^k-\psi_i^{k-1}\|_A^2 < \epsilon_{loc}^2$;}
\STATE{$\tilde{\psi}_i \leftarrow \psi_i^k$;}
\ENDFOR
\STATE{$\widetilde{\Psi} = [ \tilde{\psi}_0, \tilde{\psi}_1, \ldots, \tilde{\psi}_{\dim(\Phi)} ]$.}
\end{algorithmic}
\end{algorithm}

\section{Construction of partition $\mathcal{P}$}
\label{Sec:partition}
With the analysis in the previous sections, we now have a blueprint for the construction of partition $\mathcal{P}=\{P_j\}_{j=1}^M$. Given an underlying energy decomposition $\mathcal{E}=\{E_k\}_{k=1}^m$ of a SPD matrix $A$ and an orthonormal basis $\mathcal{V}$ of $\mathbb{R}^n$, the basic idea is to find a partition $\mathcal{P}$ of $\mathcal{V}$ with small patch number $\#\mathcal{P}$ and small {\bf condition factor} $\delta(\mathcal{P},q)$, while subject to a prescribed error bound on the {\bf error factor} $\varepsilon(\mathcal{P},q)$. In particular, our goal is to find the optimizer of the following problem: 
\begin{equation*}
\begin{array}{rl}
 \mathcal{P}=\underset{\widetilde{\mathcal{P}}}{\argmin} & f_1(\#\widetilde{\mathcal{P}})+f_2(\delta(\widetilde{\mathcal{P}},q)),\\\\
\text{subject to } & \varepsilon(\widetilde{\mathcal{P}},q)\leq \epsilon,
\end{array}
\end{equation*}
where $f_1,f_2$ are some penalty functions, $q$ is a chosen integer, and $\epsilon$ is the desired accuracy. This ideal optimization problem is intractable, since in general such discrete optimization means to search over all possible combinations. Instead, we propose to use local clustering approach to ensure efficiency. 

Generally, if we have a priori knowledge of the underlying computational domain of the problem, like $\Omega\subset\mathbb{R}^d$, one of the optimal choices of partition will be the uniform regular partition. For instance, in \cite{brenner2004finite}, regular partitions are used in the sense that each patch (finite element) has a circumcircle of radius $H$ and an inscribed circle of radius $\rho H$ for some $\rho\in(0,1)$. The performance under regular partitioning relies on the regularity of the coefficients of $A$ (low contrast, strong ellipticity), and the equivalence between energy norm defined by $A$ and some universal norm independent of $A$. In particular, since regular partitioning of the computational domain is simply constructed regardless of the properties of $A$, its performance cannot be ensured when $A$ loses some regularity in some local or micro-scaled regions. 

In view of this, a more reasonable approach is to construct a partition $\mathcal{P}$ based on the information extracted from $A$, which is represented by the local energy decomposition $\mathcal{E}$ of $A$ in our proposed framework. For computational efficiency, the construction procedure should rely only on local information (rather than global spectral information as in the procedure of Eigendecomposition). This explains why we introduce the local measurements in \Cref{sec:operator_compression}: the {\bf error factor} $\epsilon(\mathcal{P},q)$ and the {\bf condition factor} $\delta(\mathcal{P},q)$, which keep track of the performance of partition in our searching approach. These measurements are locally (patch-wisely) computable and thus provide the operability of constructing partition with local operations interacting with only neighbor data.

To make use of the local spectral information, we propose to construct the desired partition $\mathcal{P}$ of $\mathcal{V}$ by iteratively clustering basis functions in $\mathcal{V}$ into patches. In particular, small patches(sets of basis) are combined into larger ones, and the scale of the partition becomes relatively coarser and coarser. For every such newly generated patch $P_j$, we check if $\varepsilon(P_j,q)$ still satisfies the required accuracy (See \Cref{eqt:localphi_constraint}). The whole clustering process stops when no patch combination occurs, that is, when the partition achieves the resolution limit. Also, for patch $P_j$ to be well-conditioned, we set a bound $c$ on $\delta(P_j,q)\varepsilon(P_j,q)^2$. The motivation of such bound will be explained in \Cref{Sec:Multiresolution}. And for large $\delta(P_j,q)$ to diminish, patches with large {\bf condition factor} are combined first. To realize the partitioning procedure and maintain the computation efficiency, we combine patches pair-wisely. Our proposed clustering algorithm is summarized in \Cref{alg:pair_clustering} and \Cref{alg:find_match}.

\begin{algorithm}[!h]
\caption{\it Pair-Clustering}
\label{alg:pair_clustering}
\begin{algorithmic}[1]
\REQUIRE{energy decomposition $\mathcal{E}$, underlying basis $\mathcal{V}$, desired accuracy $\epsilon$, condition bound $c$.}
\ENSURE{Partition $\mathcal{P}$.}
\STATE{\textbf{Initialize:} $P_j=\{v_j\}$, $\delta(P_j,q)=\overline{A}_{P_j}$ (scalar), $1\leq j\leq n$;}
\WHILE{Number of {\it active} patches $>0$,} \label{line:pair-clustering_line3}
\STATE{Sort {\it active} $\{P_j\}$ with respect to $\delta(P_j,q)$ in descending order;}\label{alg:line3_pairclustering}
\STATE{Mark all ({\it active} and {\it inactive}) $P_j$ as {\it unoperated};}
\FOR{each active $P_j$ in descending order of $\delta(P_j,q)$,}
\STATE{\textbf{Find\_Match}$(P_j,\epsilon,c)$;} \label{alg:line6_pairclustering}
\IF{\textbf{Find\_Match} succeeds,}
\STATE{Mark $P_j$ as {\it operated};}
\ELSIF{all neighbor patches of $P_j$ are {\it unoperated},}
\STATE{Mark $P_j$ as {\it inactive};}
\ENDIF
\ENDFOR
\ENDWHILE
\end{algorithmic}
\end{algorithm}

\begin{algorithm}[!h]
\caption{\it Find\_Match}
\label{alg:find_match}
\begin{algorithmic}[1]
\REQUIRE{$P_j$, $\epsilon$, $c$.}
\ENSURE{Succeeds or Fails.}
\FOR{$P_{j''}\sim P_j$}
\STATE{Find largest $\mathrm{Con}(P_j,P_{j''})$ among all {\it unoperated} $P_{j''}$ (stored as $P_{j'}$);}
\ENDFOR
\STATE{Compute $\varepsilon(P_j\cup P_{j'},q)$ and $\delta(P_j\cup P_{j'},q)$;}
\IF{$\varepsilon(P_j\cup P_{j'},q)\leq \epsilon$ \& $\delta(P_j\cup P_{j'},q)\varepsilon(P_j\cup P_{j'},q)^2\leq c$,}
\STATE{combine $P_j$ and $P_{j'}$ to form $P_j$ ($P_{j'}$ no longer exists);}
\STATE{update $\delta(P_j,q)$;}
\RETURN{\textbf{Find\_Match} succeeds.}
\ELSE
\STATE{\textbf{return} \textbf{Find\_Match} Fails.}
\ENDIF
\end{algorithmic}
\end{algorithm}

\begin{remark}
$ $ \linebreak
\vspace{-3mm}
\begin{itemize}
\item[-] If we see $1/\varepsilon(P_j\cup P_{j'},q)^2$ as the gain, and $\delta(P_j\cup P_{j'},q)$ as the cost, the well-conditioning bound $\varepsilon(P_j\cup P_{j'},q)^2\delta(P_j\cup P_{j'},q)\leq c$ implies that the cost is proportional to the gain. 
\item[-] If we want the patch sizes to grow homogeneously, we can take patch size into consideration when sorting the patches (Line~\ref{alg:line3_pairclustering} in \Cref{alg:pair_clustering}).
\item[-] The local basis functions, $\Phi_j$, are also computed in the sub-function \textbf{Find\_Match}, and can be stored for future use.
\end{itemize}
\end{remark}

The sub-function \textbf{Find\_Match} in Line~\ref{alg:line6_pairclustering} of \Cref{alg:pair_clustering} takes a patch $P_j$ as input and finds another patch $P_{j'}$ that will be absorbed by $P_j$. As a local operation, the possible candidates for $P_{j'}$ are just the neighboring patches of $P_j$. To further accelerate the algorithm, we avoid checking the {\bf error factor} for all possible pair $(P_j,P_{j'})$ with $P_{j'}\sim P_j$. Alternatively, we check the patch $P_{j'}$ that has the largest ``connection" (correlation) with $P_j$. Undoubtedly, this quantity can be defined in different ways. Here we propose the {\it connection} between $P_j$ and $P_{j'}$ as:
\begin{equation}
\mathrm{Con}(P_j,P_{j'})=\sum_{E\sim P_j,E\sim P_{j'}}\Big(\sum_{\begin{subarray}{c}u\in P_j,v\in P_{j'}\\ u\sim v\end{subarray}}|u^TEv|\Big).
\end{equation}
On the one hand, noted that $\mathrm{Con}(P_j,P_{j'})$ can be easily computed and inherited directly after patch combination since one can check that $\mathrm{Con}(P_j\cup P_{j'},P_{j''})=\mathrm{Con}(P_j,P_{j''})+\mathrm{Con}(P_{j'},P_{j''})$. On the other hand, we observe that 
\begin{equation}
\underline{A}_{P_j\cup P_{j'}}=\underline{A}_{P_j}+\underline{A}_{P_{j'}}+\sum_{\begin{subarray}{c} E\sim P_j,E\sim P_{j'}\\ E\in P_j\cup P_{j'}\end{subarray}}E=\underline{A}_{P_j}+\underline{A}_{P_{j'}}+\mathrm{Cross\ Energy}.
\end{equation}
In other words, a larger cross energy implies larger interior eigenvalues of $\underline{A}_{P_j\cup P_{j'}}$, which means $\underline{A}_{P_j\cup P_{j'}}$ is less likely to violate the accuracy requirement. One can also recall the similarity of this observation to the findings in Spectral Graph Theory, where stronger connectivity of the graph corresponds to larger eigenvalues of the graph Laplacian $L$. These motivate us to simplify the procedure by examining the patch candidate $P_j'$ with largest connection to $P_j$.

Though our algorithm does not assume any a priori structural information of $A$, its efficiency and effectiveness may rely on the hidden locality properties of $A$. To perform a complexity analysis of \Cref{alg:pair_clustering}, we first introduce some notations. Similar to the {\it layers of neighbors} defined in \Cref{eqt:layer_neighbors}, we define $\mathcal{N}_1(v)=\mathcal{N}(v)$, and 
\[\mathcal{N}_{k+1}(v)=\mathcal{N}(\mathcal{N}_k(v))=\{u\in \mathcal{V}: u\sim \mathcal{N}_k(v)\},\] 
that is, for any $u\in \mathcal{N}_k(v)$, there is a path of length $k$ that connects $u$ and $v$ with respect to the connection relation ``$\sim$'' defined in \Cref{def:neighboring}. The following definition describes the local interaction property of $A$:

\begin{definition}[Locality/Sparsity of $A$]
$A$ is said to be local of dimension $d$ with respect to $\mathcal{V}$, if 
\begin{equation}
\#\mathcal{N}_k(v)= O(k^d),\quad \forall k\geq 1,\quad \forall v\in \mathcal{V}.
\label{condition:localityA}
\end{equation}
\end{definition}

The following definition describes the local spectral properties of an energy decomposition of $A$. It states that a smaller local patch corresponds to a smaller scale, and that $\varepsilon(P,q)$ tends to increase and $\delta(P,q)$ tends to decrease as the patch size of $P$ increases. This explains why we combine patches from finer scales to coarser scales to construct the desired partition $\mathcal{P}$. 

\begin{definition}[Local energy decomposition]
$\mathcal{E}=\{E_k\}_{k=1}^m$ is said to be a {\bf local energy decomposition} of $A$ of order $(q,p)$ with respect to $\mathcal{V}$, if there exists some constant $h$, such that 
\begin{equation}
\varepsilon(\mathcal{N}_k(v),q)= O((hk)^p) ,\quad \forall k\geq 1,\quad \forall v\in \mathcal{V}.
\label{condition:scalesproperty}
\end{equation}
Moreover, $\mathcal{E}$ is said to be well-conditioned if there is some constant $c$ such that 
\begin{equation}
\varepsilon(\mathcal{N}_k(v),q)^2\delta(\mathcal{N}_k(v),q)\leq c,\quad \forall k\geq 1,\quad \forall v\in \mathcal{V}.
\label{condition:localregularity}
\end{equation}
\end{definition}

\begin{remark}
\label{rmk:remark_locality}
$ $ \linebreak
\vspace{-3mm}
\begin{itemize}
\item[-] The locality of $A$ implies that $\#\mathcal{N}_{k+1}(v)-\#\mathcal{N}_k(v)=O(d\cdot k^{d-1})$. In particular $\#\mathcal{N}_1(v)=O(d)$, and thus the number of nonzero entries of $A$ is $m=O(d\cdot n)$.
\item[-] Let $\mathcal{P}$ be a partition of $\mathcal{V}$ such that each patch $P\in \mathcal{P}$ satisfies $\mathrm{diam}(P)=O(r)$ and $\# P=O(r^d)$, where $r$ is an integer, and ``$\mathrm{diam}$'' is the path diameter with respect to the adjacency relation ``$\sim$'' defined in \Cref{def:neighboring}. Let $S_k(P)$ be the {\it layers of neighbors} (patch layers) defined in \Cref{eqt:layer_neighbors}, and $\#^{P}S_k(P)$ denote the number of patches in $S_k(P)$. Then the locality of $A$ implies that 
\[\#^{P}S_k(P)=O(\frac{\#S_k(P)}{r^d})=O(\frac{(rk)^d}{r^d})=O(k^d).\]
This means that a $\mathcal{V}$ with adjacency relation defined by $A$ has a self-similar property between fine scale and coarse scale.
\end{itemize}
\end{remark}

These abstract formulations/notations actually summarize a large class of problems of interest. For instance, suppose $A$ is assembled from the FEM discretization of a well-posed elliptic equation with homogeneous Dirichlet boundary conditions:
\[\mathcal{L}u=\sum_{0\leq |\sigma|, |\gamma|\leq p}(-1)^{|\sigma|}D^\sigma(a_{\sigma\gamma}D^\gamma u)=f,\quad u\in H_0^p(\Omega),\]
where $\Omega\subset\mathbb{R}^d$ is a bounded domain. Let $\mathcal{V}$ be the nodal basis of the discretization, and each energy element $E$ in $\mathcal{E}$ be the energy inner product matrix (i.e. the stiffness matrix) of the neighbor nodal functions on a fine mesh patch. The locality of $\mathcal{L}$ and the underlying dimension of $\Omega$ ensure that $A$ is local of dimension $d$ with respect to $\mathcal{V}$. With a consistent discretization of $\mathcal{L}$ on local domains, the interior energy corresponds to a Neumann boundary condition, while the closed energy corresponds to a Dirichlet boundary condition. In this sense, using a continuous limit argument and the strong ellipticity assumption, Hou and Zhang in \cite{hou2016sparse} prove that if $q\geq\left(\begin{smallmatrix}p+d-1\\d \end{smallmatrix}\right)$, then $\varepsilon(\mathcal{N}_k(v),q)\lesssim (hk)^p$ (generalized Poincar\'e inequality) where $h$ is the fine mesh scale and $p$ is half the order of the elliptic equation; and $\varepsilon(\mathcal{N}_k(v),q)^2\delta(\mathcal{N}_k(v),q)\leq c$ (inverse estimate) for some scaling-invariant constant $c$. Unfortunately these arguments would be compromised if strong ellipticity is not assumed, especially when high contrast coefficients are present. However, as we see in \Cref{example:example3}, $\varepsilon(\mathcal{N}_k(v),q)$ and $\delta(\mathcal{N}_k(v),q)$ actually converge when the contrast of the coefficient becomes large, which is not explained by general analysis. So we could still hope that the matrix $A$ and the energy decomposition $\mathcal{E}$ have the desired locality that can be numerically learned, even when conventional analysis fails. 

\subsection*{Estimate of patch number} Intuitively, if $A$ is local of dimension $d$, and $\mathcal{E}$ is well-conditioned and local of order $q$, then the patch number of an ideal partition $\mathcal{P}$ subject to accuracy $\epsilon$ should be 
\begin{equation}
\#\mathcal{P}=O\Big(\frac{n}{(\epsilon^{1/p}/h)^d}\Big)=O\Big(\frac{nh^d}{\epsilon^{d/p}}\Big),\quad \delta(\mathcal{P},q)=O\big(\frac{c}{\epsilon^2}\big),
\label{eqt:patchnumber_dimension}
\end{equation}
where we estimate the path diameter of each patch in $\mathcal{P}$ by $O(\epsilon^{1/p}/h)$, and thus the patch size by $O((\epsilon^{1/p}/h)^d)$.

\subsection*{Inherited locality} As we have made the locality assumption on $A$, we would hope that the compressed operator $P_{\widetilde{\Psi}}^{A}A^{-1}=\widetilde{\Psi}\widetilde{A}_{st}^{-1}\widetilde{\Psi}^T$ can also take advantage of such locality. In fact, the localization of $\Psi$ not only ensures the efficiency of the construction of $\widetilde{\Psi}$, but also conveys the locality of $A$ to the stiffness matrix $\widetilde{A}_{st}$. Suppose $A$ is local of dimension $d$. Let $\widetilde{\Psi}$ be the local approximator obtained in \Cref{alg:construct_tilde_psi} such that $\widetilde{\psi}_i=\psi_i^r,\ 1\leq i\leq N$ for some uniform radius $r$. Let $\mathcal{\widetilde{V}} = \{\tilde{v}_i\}_{i=1}^N$ be the orthonormal basis of $\mathbb{R}^{N}$ such that 
\[ \tilde{v}_i^T\widetilde{A}_{\text{st}}\tilde{v}_j=\tilde{v}_i^T\widetilde{\Psi}^TA\widetilde{\Psi}\tilde{v}_j=\widetilde{\psi}_i^TA\widetilde{\psi}_j,\quad \forall\ 1\leq i,j\leq N. \]
We then similarly define the adjacency relation between basis vectors in $\mathcal{\widetilde{V}}$ with respect to $\widetilde{A}_{\text{st}}$. Since each localized basis $\tilde{\psi}_i$ interacts with patches in $r$ patch layers (thus interacts with other $\tilde{\psi}_{i'}$ corresponding to patches in $2r$ layers), and each patch corresponds to $q$ localized basis functions, using the result in \Cref{rmk:remark_locality} we have 
\[\#\mathcal{N}_k(\tilde{v})=O((rk)^d\cdot q)=q\cdot r^d\cdot O(k^d),\quad \forall k\geq 1,\quad \forall \tilde{v}\in \mathcal{\widetilde{V}}.\]
That means $\widetilde{A}_{\text{st}}$ inherits the locality of dimension $d$ from $A$. In addition, by using the same argument as in \Cref{rmk:remark_locality}, we have $\#\mathcal{N}_1(\tilde{v})=O(d\cdot q\cdot r^d)$, which is the number nonzero entries (NNZ) of one single column of $\widetilde{A}_{\text{st}}$. Therefore the number of nonzero entries(NNZ) of $\widetilde{A}_{\text{st}}$ can be bounded by $O(N\cdot d \cdot q \cdot r^d)=O(m\cdot q^2\cdot r^d/s)$, where $m=O(d\cdot n)$ is the NNZ of $A$, $s$ is the average patch size, and we have used the relation $N=nq/s$. In particular, if the localization error is subject to $\epsilon_{loc}=\frac{\epsilon}{\sqrt{N}}$, then the radius has estimate $r=O(\log\frac{1}{\epsilon}+\log n+\log\delta(\mathcal{P},q))$, and thus the bound on the NNZ of $\widetilde{A}_{\text{st}}$ becomes
\begin{equation}
O(m\cdot q^2\cdot \frac{1}{s} \cdot (\log\frac{1}{\epsilon}+\log n+\log\delta(\mathcal{P},q))^d).
\end{equation}

\subsection*{Choice of $q$} Recall that, instead of choosing a larger enough $q_j$ for each patch $P_j$ to satisfied $\varepsilon(P_j,q_j)\leq \epsilon$, we use a uniform integer $q$ for all patches, and leave the mission of accuracy to the construction of partition. So before we proceed to the algorithm, we still need to know what $q$ we should choose. In some problems, $q$ can be determined by theoretical analysis. For example, when solving elliptic equation of order $p$ with FEM, we should at least choose $q=\left(\begin{smallmatrix}p+d-1\\d \end{smallmatrix}\right)$ to obtain an optimal rate of convergence. And as for graph laplacians that are generally considered as a discrete second order elliptic problem ($p=1$), we can thus choose $q=1$. But when the problem is more complicated and has no intrinsic order, the choice of $q$ can be tricky. So one practical strategy is to start from $q=1$, and increase $q$ when the partition obtained is not acceptable. 

\subsection*{Complexity of \Cref{alg:pair_clustering}} 
\label{subsec:complexity_paircluster}
For simplicity, we assume that all patches in the final output partition $\mathcal{P}$ have the same patch size $s$. Under locality assumption of $A$ given in \cref{condition:localityA}, the local operation cost of \textbf{Find\_Match}($P_j$) is approximately 
\[O(d\cdot F(\mathrm{size}(P_j)))\leq O(d\cdot F(s)).\]
Here $F(\#P_j)$ is a function of $\#P_j$ that depends only on complexity of solving local eigen problems(with respect to $\underline{A}_{P_j}$) and local inverse problems(with respect to $\overline{A}_{P_j}$) on patch $P_j$, thus we can bound $F(s)$ by $O(s^3)$. Since patches are combined pair-wise, the number of while-loops starting at Line~\ref{line:pair-clustering_line3} is of order $O(\log s)$, and in each while-loop, the operation cost can be bounded by 
\[O(d\cdot \frac{F(s)}{s}\cdot n)+O(n\cdot\log n)=O(d\cdot s^2\cdot n)+O(n\cdot\log n),\] 
where $O(d\cdot \frac{F(s)}{s}\cdot n)$ comes from operating \textbf{Find\_Match}($P_j$) for all surviving active patch $P_j$, and $O(n\cdot\log n)$ comes from sorting operations. Therefore the total complexity of \Cref{alg:pair_clustering} is 
\begin{equation}
O\big(d\cdot s^2\cdot\log s\cdot n\big)+O(\log s\cdot n\cdot\log n).
\label{eqt:pair_clustering_complexity}
\end{equation}

\subsection*{Complexity of \Cref{alg:operator_compression}} Combining all procedures together and noticing that \Cref{alg:construct_phi} can be absorbed to \Cref{alg:pair_clustering}, the complexity of \Cref{alg:operator_compression} is 
\begin{equation}
O\big(d\cdot s^2\cdot\log s\cdot n\big)+O(\log s\cdot n\cdot\log n)+O\big(q\cdot n\cdot \epsilon^2\cdot \|A\|_2\cdot (\log\frac{1}{\epsilon}+\log n+\log\delta(\mathcal{P}))^{d+1}\big),
\label{eqt:operatorcompression_complexity}
\end{equation}
where $n$ is the original dimension of basis, $s$ is the maximal patch size, $\epsilon$ is the prescribed accuracy, and we have used the fact $\varepsilon(\mathcal{P},q)\leq \epsilon$.

\begin{remark}
$ $ \linebreak
\vspace{-3mm}
\begin{itemize}
\item[-] The maximum patch size $s$ can be viewed as the compression rate since $s\sim \frac{n}{M}$ where $M$ is the patch number. The complexity analysis above implies that for a fixed compression rate $s$, complexity of \Cref{alg:pair_clustering} is linear in $n$. However, when the locality conditions \cref{condition:localityA} and \cref{condition:scalesproperty} are assumed, one can see that $s\sim \big(\epsilon^\frac{1}{p}/h\big)^d$, where $\epsilon$ is the desired (input) accuracy. As a consequence, while the desired accuracy is fixed, if $n$ increases, then the complexity is no longer linear in $n$ since $h$ (finest scale of the problem) may change as $n$ changes. In other words, \Cref{alg:pair_clustering} loses the near linear complexity when the desired accuracy $\varepsilon$ is set too large compared to the finest scale of the problem. To overcome such limitation, we should consider a hierarchical partitioning introduced in \Cref{Sec:Multiresolution}.
\item[-] As we mentioned before, the factor $\epsilon^2\|A\|_2$ also suggests a hierarchical compression strategy, when $\|A\|_2$ is too large compared to the prescribed accuracy $\epsilon^2$.
\end{itemize}
\end{remark}

\section{Numerical Examples}
\label{sec:numerical1}
In this section, two numerical examples are reported to demonstrate the efficacy and effectiveness of our proposed operator compression algorithm. For consistency, all the experiments are performed on a single machine equipped with Intel(R) Core(TM) i5-4460 CPU with 3.2GHz and 8GB DDR3 1600MHz RAM.

\subsection{Numerical Example 1}
\label{subsec:numerical_example1}
\begin{figure}[!h]
    \begin{subfigure}[b]{0.24\textwidth}
        \centering
        \includegraphics[width=1.0\textwidth]{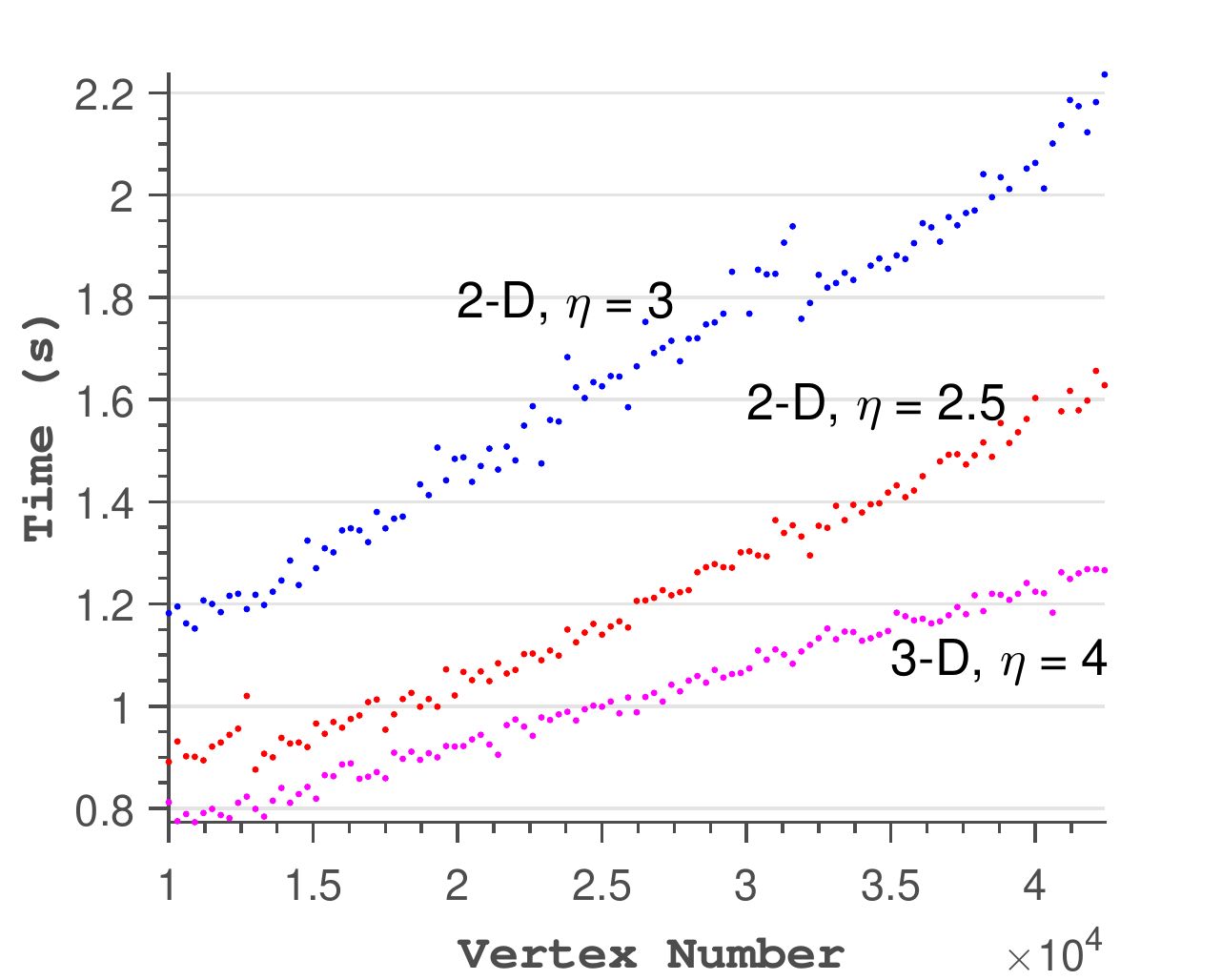}
        \caption{Running time}
        \label{fig:complexity_all}
    \end{subfigure}
    \begin{subfigure}[b]{0.24\textwidth}
        \centering
        \includegraphics[width=1.0\textwidth]{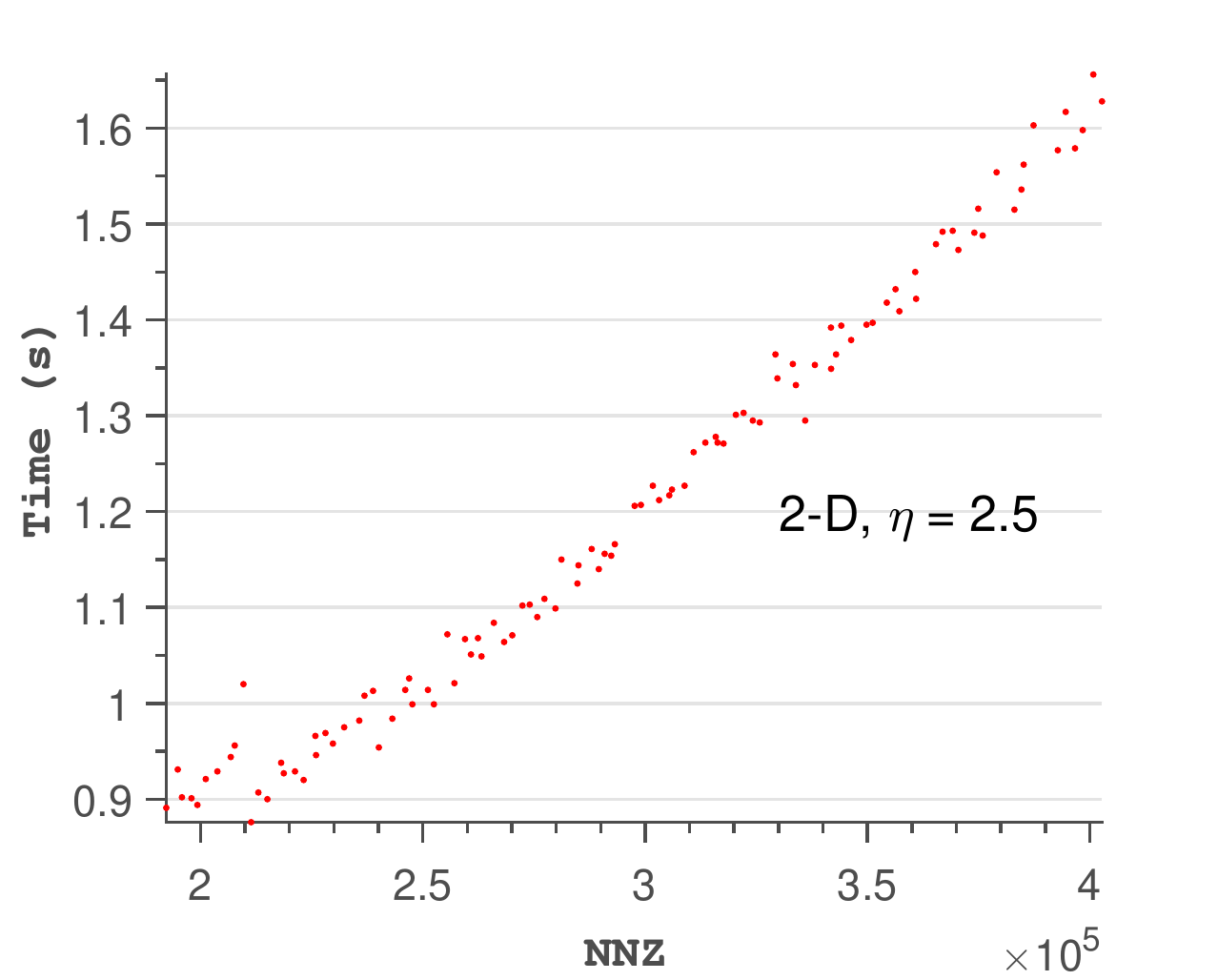}
        \caption{2-D graphs with $\eta = 2.5$}
        \label{fig:complexity_2d25}
    \end{subfigure}
        \begin{subfigure}[b]{0.24\textwidth}
        \centering
        \includegraphics[width=1.0\textwidth]{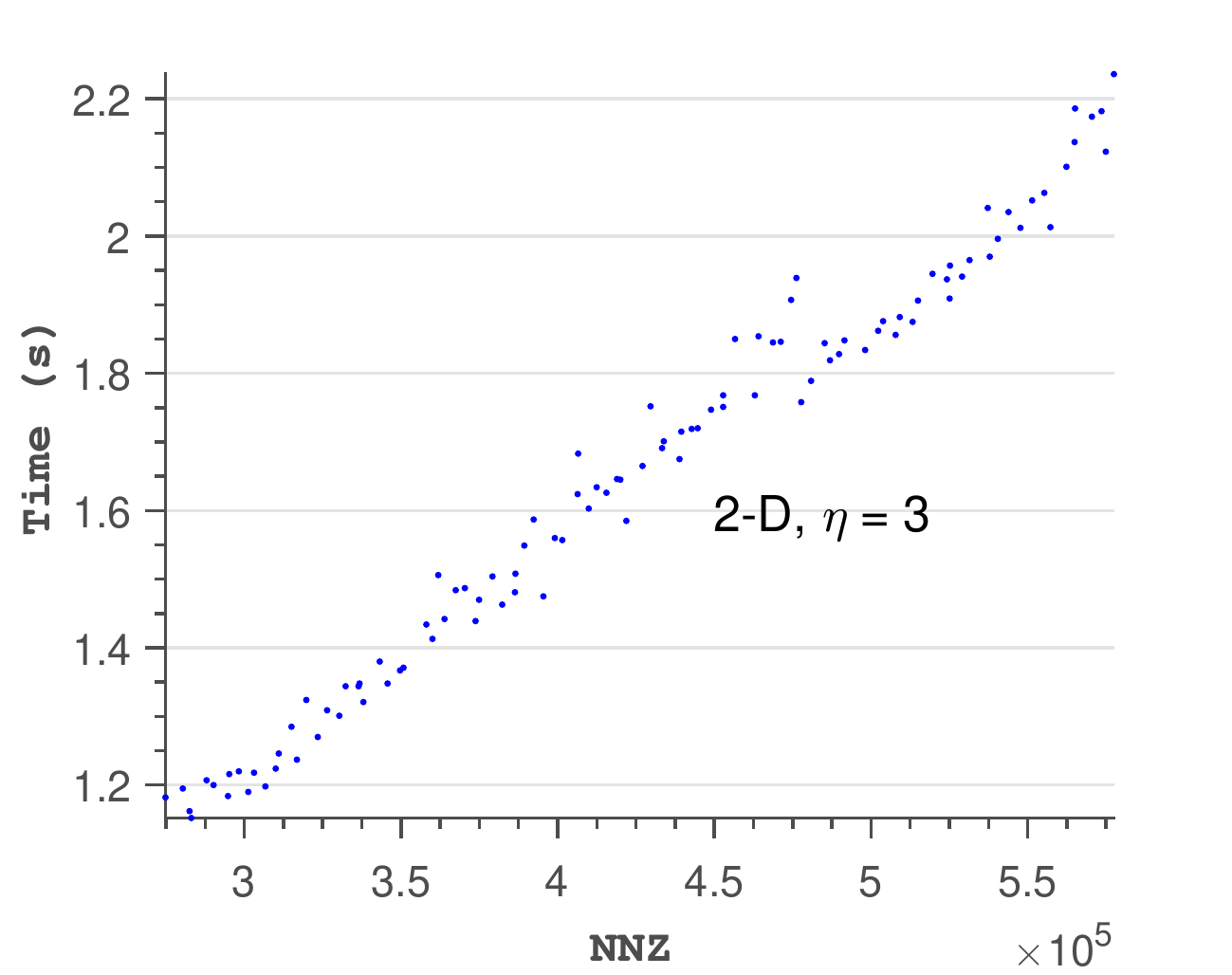}
        \caption{2-D graphs with $\eta = 3$}
        \label{fig:complexity_2d3}
    \end{subfigure}
    \begin{subfigure}[b]{0.24\textwidth}
        \centering
        \includegraphics[width=1.0\textwidth]{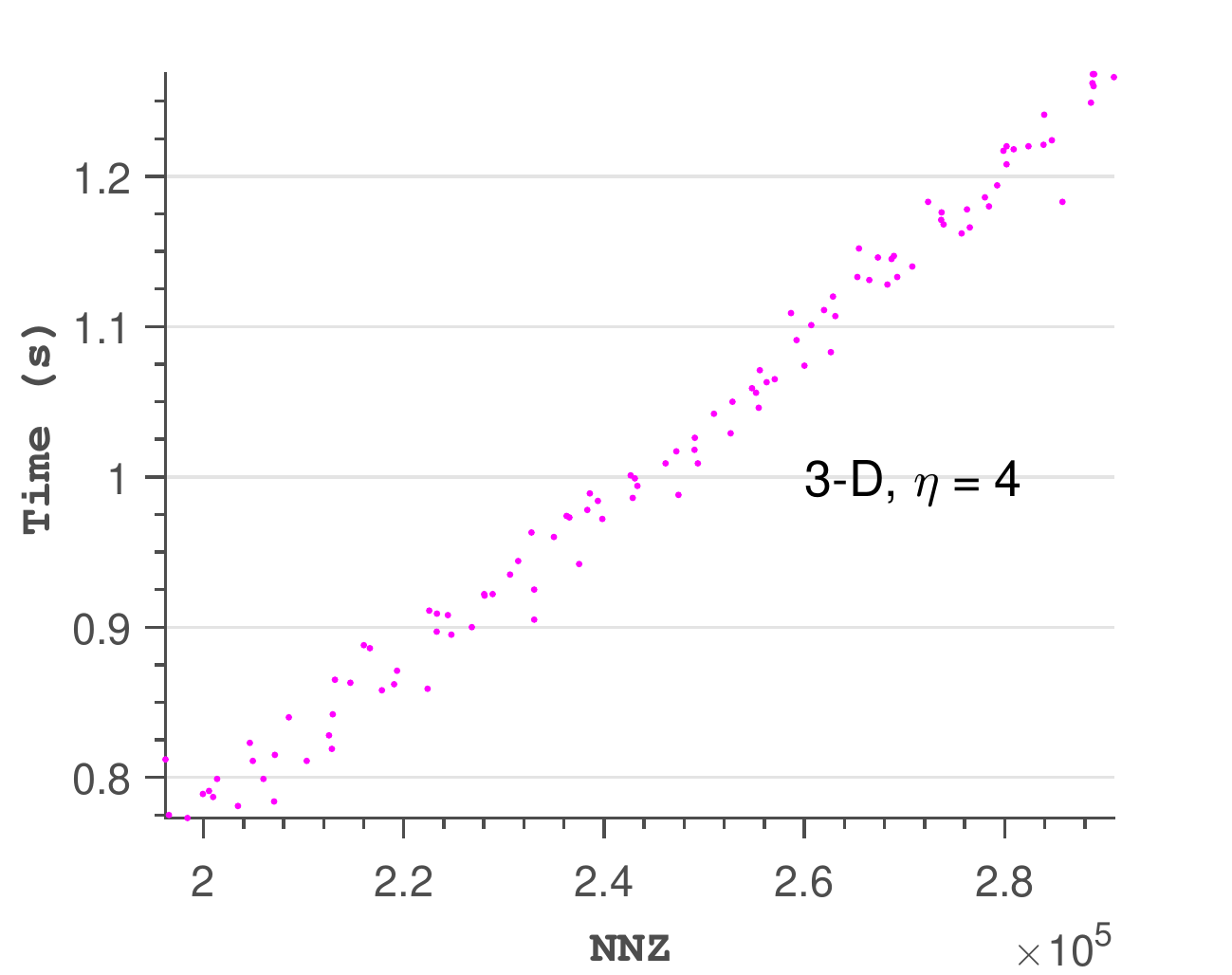}
        \caption{3-D graphs with $\eta = 4$}
        \label{fig:complexity_3d4}
    \end{subfigure}
    \caption{Illustration of the running time of \Cref{alg:pair_clustering}. (a) plots the running time against the number of vertices of random generated graphs in different dimensions and density factor $\eta$. (b) -- (d) plot the running time against the number of nonzero entries (NNZ) in the graph laplacian operators.}
    \label{fig:complexity}
\end{figure}
The first numerical example arises from solving the finite graph Laplacian system $Lx=b$, where $L$ is the Laplacian matrix of a $d$-dimensional undirected random graph $G=[V,E]$. Vertices $\bm{x}_i=(x_1^i,\cdots,x_d^i)\in V\subset\mathbb{R}^d$ are generated subject to a uniform distribution over the domain $\Omega = [0,1]^d$. Edge weights are then given by
\begin{equation*}
w_{ii}=c_i,\quad \forall i;\quad 
w_{ij}=\left\{
\begin{array}{ll}
r_{ij}^{-2}, & \mathrm{if}\ r_{ij}^2\leq \eta / n^{\frac{2}{d}}, \\
0, & \mathrm{otherwise},
\end{array}\right.\quad \forall\ i\neq j,
\end{equation*}
where $r_{ij}=\|\bm{x}_i-\bm{x}_j\|_2$, $n=\# V$ is the number of vertices and $\eta>0$ is some density factor for truncating long distance interactions. We set $c_i = 1 \ \forall i$ for the sake of well-posedness and invertibility of the graph laplacian, which gives $\Vert L^{-1} \Vert_2 = 1$. We also remark that our choice of $\eta$ ensures that the graph is locally connected and that the second smallest eigenvalue of $L$ is of order $O(1)$. This also gives $n \propto \text{Number of NonZero entries (NNZ)}$ of $L$ in this example (which is actually not required by our algorithm). The basis $\mathcal{V}\subset\mathbb{R}^n$ is given by the natural basis with respect to vertices values, and the energy decomposition $\mathcal{E}=\{E_k\}_{k=1}^m$ is collected as described in \Cref{example:Example2}, where each $E_k$ corresponds to an edge in $G$. Since $p=1$ for graph laplacian, we set $q = 1$ throughout this numerical example. 

We first verify the complexity of \Cref{alg:pair_clustering} by applying it to partition random graphs generated as described above. To be consistent, we set the prescribed accuracy $\frac{1}{\varepsilon^d} \propto n$ and the upper bound $c$ of $\delta(\mathcal{P},1)\varepsilon(\mathcal{P},1)^2$ to be 100 in all cases, which is large enough for patches to combine with each others. \Cref{fig:complexity} illustrates the nearly-linear time complexity of our algorithm with respect to the graphs' vertex number $n$, which is consistent to our complexity estimation in \Cref{subsec:complexity_paircluster}. Every dot represents the partitioning result of the given 2-D/3-D graph. In particular, the red and blue sets of dots are the partition results for 2-D graphs under the construction of $\eta = 2.5$ and $\eta = 3$ respectively, while the magenta point set represents the partitioning of the random 3-D graphs under the setting of $\eta = 4$. Similarly, \Cref{fig:complexity_2d25} - \Cref{fig:complexity_3d4} show respectively the time complexity of \Cref{alg:pair_clustering} versus the NNZ entries of $L$. Notice that for graphs having NNZ entries with order up to $10^5$, the running time is still within seconds. These demonstrate the lightweight nature of \Cref{alg:pair_clustering}.
\begin{figure}[!h]
    \begin{subfigure}[b]{0.33\textwidth}
        \centering
        \includegraphics[width=1.0\textwidth]{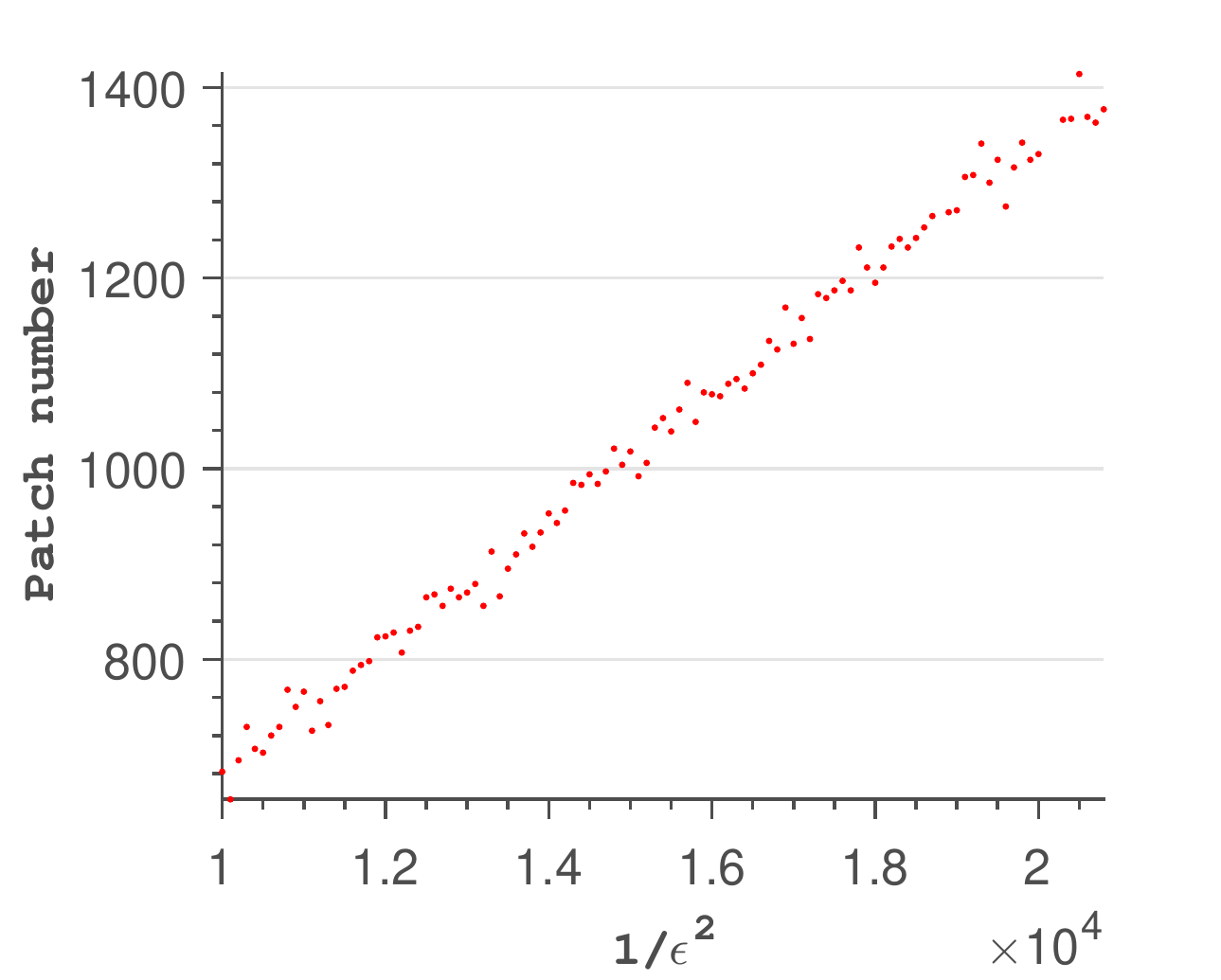}
        \caption{2-D graphs: $\# \mathcal{P} \propto \frac{1}{\epsilon^2}$}
        \label{fig:complexity2_2d25}
    \end{subfigure}
    \begin{subfigure}[b]{0.33\textwidth}
        \centering
        \includegraphics[width=1.0\textwidth]{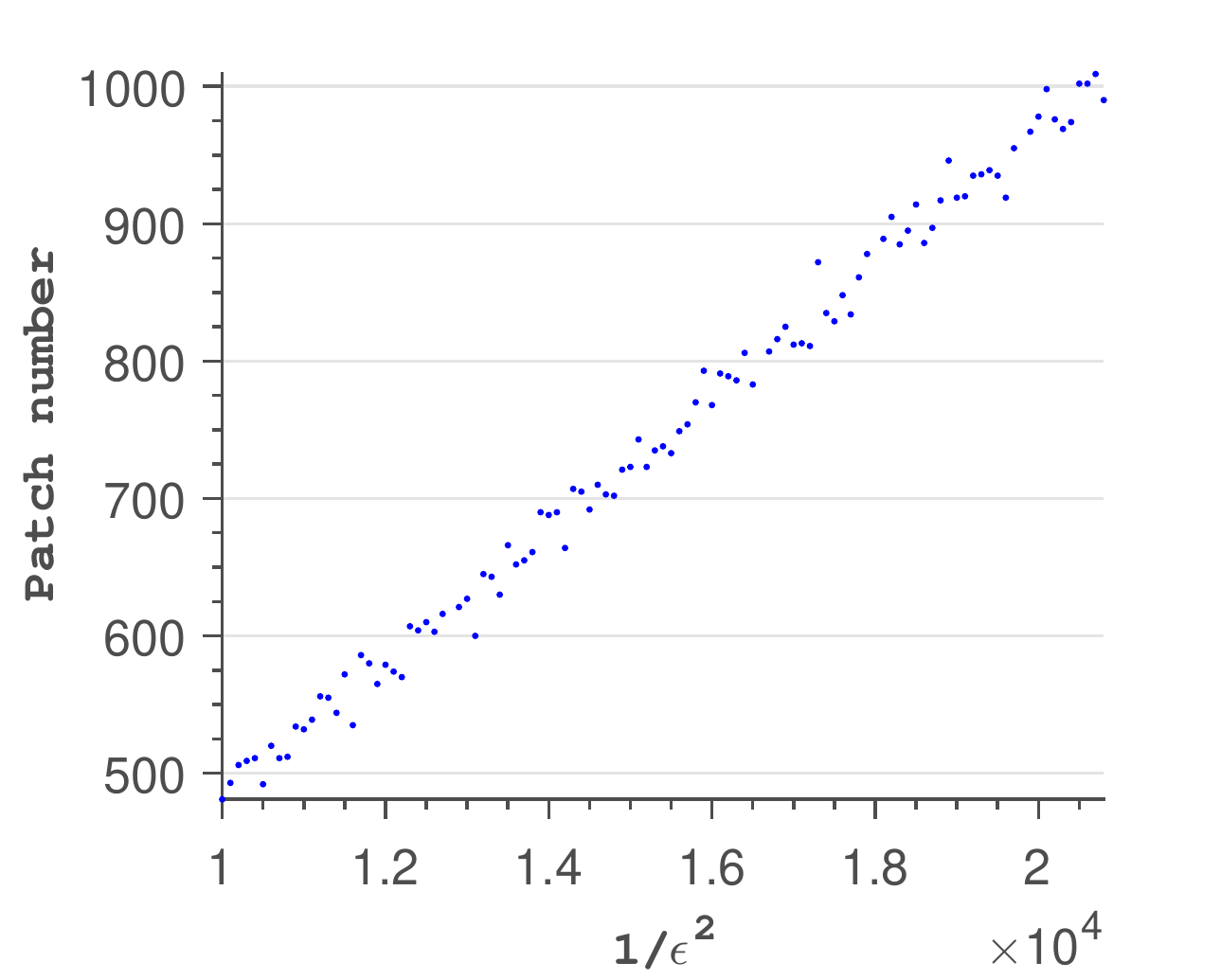}
        \caption{2-D graphs: $\# \mathcal{P} \propto \frac{1}{\epsilon^2}$}
        \label{fig:complexity2_2d3}
    \end{subfigure}
    \begin{subfigure}[b]{0.33\textwidth}
        \centering
        \includegraphics[width=1.0\textwidth]{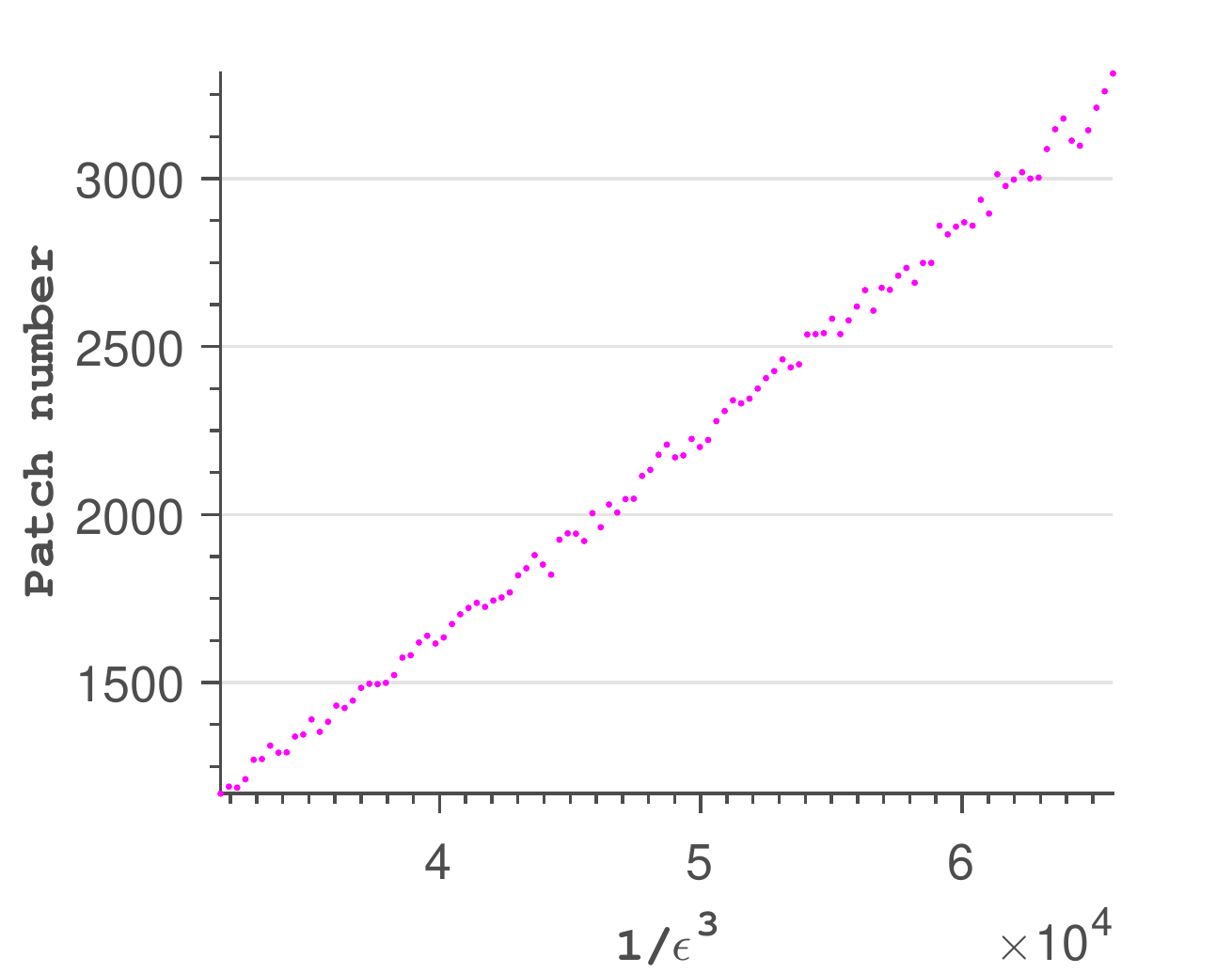}
        \caption{3-D graphs: $\# \mathcal{P} \propto \frac{1}{\epsilon^3}$}
        \label{fig:complexity2_3d4}
    \end{subfigure}
    \caption{Intrinsic dimension of the sets of graphs}
    \label{fig:complexity2}
\end{figure}

Secondly, we record the patch numbers of partition $\mathcal{P}$ obtained from \Cref{alg:pair_clustering}. We fix the domain $\Omega = [0,1] \times [0,1]$ and gradually increase the density of vertices in the graph. Therefore, we have $n \propto \frac{1}{h^d}$ and thus $\# \mathcal{P} \propto \frac{1}{\epsilon^d}$ by the observation in \Cref{eqt:patchnumber_dimension}. The relationship between $\frac{1}{\epsilon^d}$ and patch number $\# \mathcal{P}$ for the three cases is plotted in \Cref{fig:complexity2_2d25} - \Cref{fig:complexity2_3d4} respectively. \Cref{fig:complexity2_2d25} and \cref{fig:complexity2_2d3} show the linear relationship between $\# \mathcal{P}$ and $\frac{1}{\epsilon^2}$, meaning that $d=2$ in these cases. Similarly, the plot in \Cref{fig:complexity2_3d4} that discloses the dimension of the input graphs is 3-dimensional as $\# \mathcal{P} \propto \frac{1}{\epsilon^3}$. These results precisely justify the capability of our framework in capturing geometric information of the domain.

Thirdly, we focus on one particular 2-D random graph with vertex number $n=10000$ to verify the performance of our algorithm on controlling the error and well-posedness of the corresponding compressed operator. We employ the concept of the $k$-nearest neighbor (KNN) to impose local interaction. Specifically, for each vertex $\bm{x}_i$, we denote $\mathrm{NN}(\bm{x}_i,k_i)$ to be the set of the $k_i$ nearest vertices of $\bm{x}_i$. Any two vertices $\bm{x}_i$ and $\bm{x}_j$ have an edge of weight $w_{ij}=1/r_{ij}^2$ if and only if $\bm{x}_i\in \mathrm{NN}(\bm{x}_j;k_j)$ or $\bm{x}_j\in \mathrm{NN}(\bm{x}_i;k_i)$. Let $\bm{y}=(0.5,0.5)$ be the center of $\Omega$. We set $k_i=15$ if $\|\bm{x}_i-\bm{y}\|_2\leq 0.25$ and $k_i=5$ otherwise. Therefore the sub-graph inside the disk $B(\bm{y},0.25)$ has a stronger connectivity than the sub-graph outside. 

\begin{figure}[!h]
\centering
    \begin{subfigure}[b]{0.32\textwidth}
        \centering
        \includegraphics[width=1.0\textwidth]{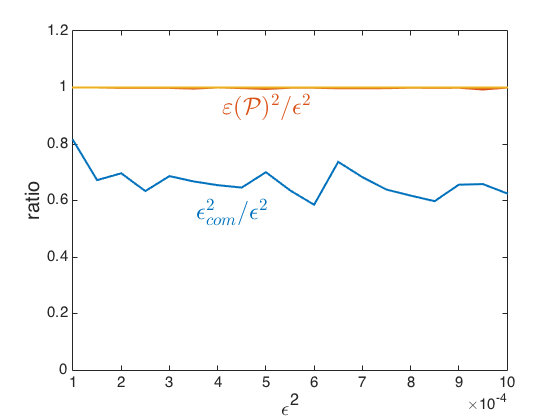}
        \caption{}
        \label{fig:errorachieve}
    \end{subfigure}
    \begin{subfigure}[b]{0.32\textwidth}
        \centering
        \includegraphics[width=1.0\textwidth]{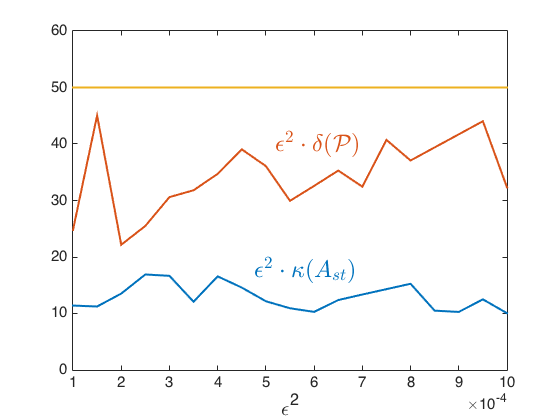}
        \caption{}
        \label{fig:conditionachieve}
    \end{subfigure}
    \begin{subfigure}[b]{0.32\textwidth}
        \centering
        \includegraphics[width=1.0\textwidth]{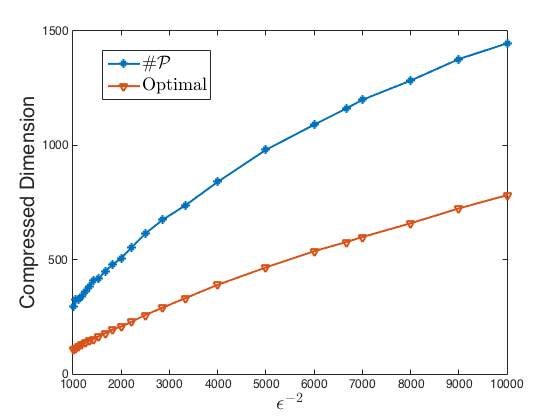}
        \caption{}
        \label{fig:dimnesionachieve}
    \end{subfigure}
    \begin{subfigure}[b]{0.32\textwidth}
        \centering
        \includegraphics[width=0.87\textwidth]{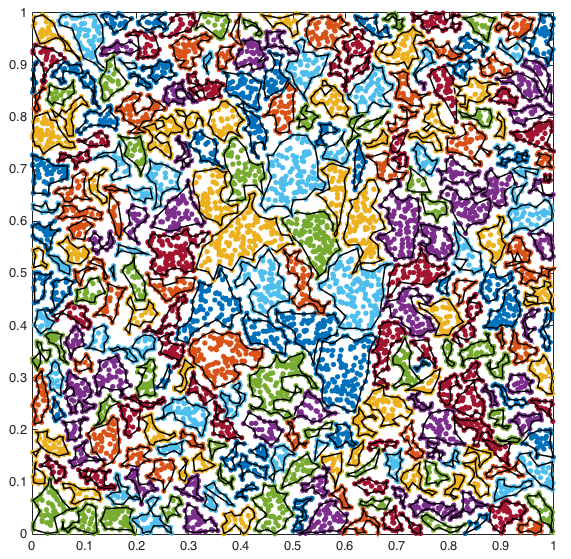}
        \caption{}
        \label{fig:patchillustration1}
    \end{subfigure}
    \begin{subfigure}[b]{0.32\textwidth}
        \centering
        \includegraphics[width=1.0\textwidth]{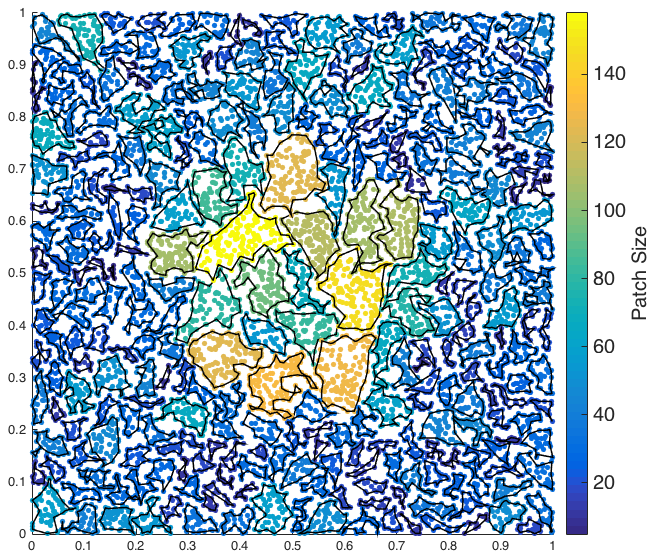}
        \caption{}
        \label{fig:patchillustration2}
    \end{subfigure}
    \caption{Error and well-posedness studies of the compressed operators}
\end{figure}

We perform \Cref{alg:operator_compression} with a fixed condition control $c=50$ and a prescribed accuracy $\epsilon^2$ varies from $0.001$ to $0.0001$. \Cref{fig:errorachieve} shows the ratios $\epsilon_{\text{com}}^2/\epsilon^2$ and $\varepsilon({\mathcal{P},1})^2/\epsilon^2$, where $\epsilon_{\text{com}}^2=\|L^{-1}-P_{\Psi}^LL^{-1}\|_2$ is the compression error. Using \Cref{alg:pair_clustering}, we achieve a nearly optimal local error control. Also notice that the global compression error ratio $\epsilon_{\text{com}}^2 / \epsilon^2$ is strictly bounded by 1 but also above $0.5$, meaning that our approach is neither playing down nor overdoing the compression. \Cref{fig:conditionachieve} shows the condition number of $A_{\text{st}}$, which is consistent to the prescribed accuracy, and is strictly bounded by $\epsilon^2 \cdot \delta(\mathcal{P},1)$ and the prescribed condition bound $c$. \Cref{fig:dimnesionachieve} plots patch number $\#\mathcal{P}$ versus $\epsilon^{-2}$ (the blue curve). Though the graph has different connectivity at different parts, it is still a 2-D graph and is locally connected in the sense of \Cref{condition:localityA}. Therefore the curve is below linear, which is consistent to estimate \cref{eqt:patchnumber_dimension} with $d=2$. As comparison, the red curve is the optimal compression dimension subject to the same prescribed accuracy given by eigenvalue decomposition. Since \Cref{alg:pair_clustering} combines patches pair-wisely, the output patch number $\#\mathcal{P}$ can be up to 2 times the optimal case. \Cref{fig:patchillustration1} shows the partition result with $\#\mathcal{P}=298$ for the case $\epsilon^2=0.001$, where the black lines outline the boundaries of patches. \Cref{fig:patchillustration2} illustrates the patch sizes of the partition. We can see that patches near the center of the domain have larger sizes than the ones near the boundary, since the graph has a higher connectivity inside the disk $B((0.5,0.5),\frac{1}{4})$.

\begin{figure}[!h]
    \begin{subfigure}[b]{0.32\textwidth}
        \centering
        \includegraphics[width=1.0\textwidth]{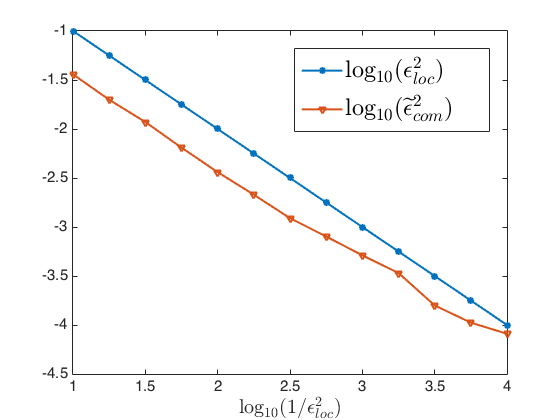}
        \caption{}
        \label{fig:erroranalysis_a}
    \end{subfigure}
    \begin{subfigure}[b]{0.32\textwidth}
        \centering
        \includegraphics[width=1.0\textwidth]{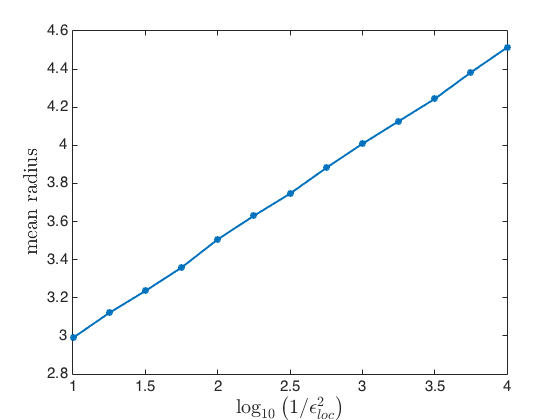}
        \caption{}
        \label{fig:erroranalysis_b}
    \end{subfigure}
    \begin{subfigure}[b]{0.32\textwidth}
        \centering
        \includegraphics[width=1.0\textwidth]{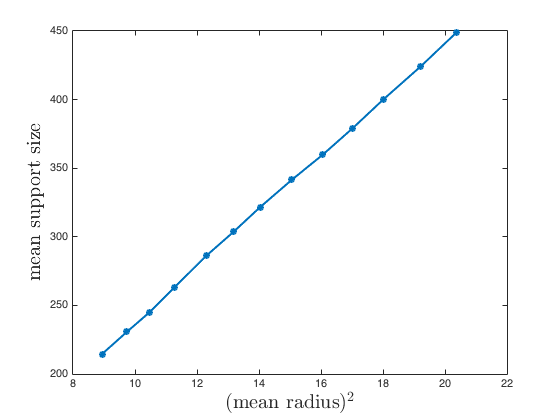}
        \caption{}
        \label{fig:erroranalysis_c}
    \end{subfigure}
    \caption{Compression error $\epsilon_{\text{com}}^2$ and the properties of the mean radius of $\widetilde{\Psi}$.}
    \label{fig:erroranalysis}
\end{figure}

We also fix the prescribed accuracy $\epsilon^2=0.0001$ to study the performance of compression with localization. In this case we have $N=\#\mathcal{P}=\dim(\Phi)=1446$. Let $\widetilde{\Psi}$ be the local approximator of $\Psi$ constructed by \Cref{alg:construct_tilde_psi} subject to localization error $\|\tilde{\psi}_i-\psi_i\|_A^2\leq\epsilon_{loc}^2$. \Cref{fig:erroranalysis} shows the compression error $\tilde{\epsilon}_{\text{com}}^2=\|L^{-1}-P_{\widetilde{\Psi}}^LL^{-1}\|_2$, mean radius and mean support size of $\widetilde{\Psi}$ with different $\epsilon_{loc}^2$ varies from $0.1$ to $0.0001$. Recall that \Cref{lemma:localization_psi} requires a localization error $\epsilon_{loc}^2=\epsilon^2 / N$ to ensure $\tilde{\epsilon}_{\text{com}}^2\leq \epsilon^2$, but \Cref{fig:erroranalysis_a} shows that $\epsilon_{loc}^2=\epsilon^2$ is adequate. \Cref{fig:erroranalysis_b} shows the linearity between mean radius and $\log\frac{1}{\epsilon_{loc}}$, which is consistent to the exponential decay of $\Psi$ proved in \Cref{thm:exponential_decay}. \Cref{fig:erroranalysis_c} shows the quadratic relation between mean support size and mean radius of $\widetilde{\Psi}$, which again reflects the geometric dimension of the graph is $2$.

\begin{figure}[!h]
    \begin{subfigure}[b]{0.32\textwidth}
        \centering
        \includegraphics[width=1.0\textwidth]{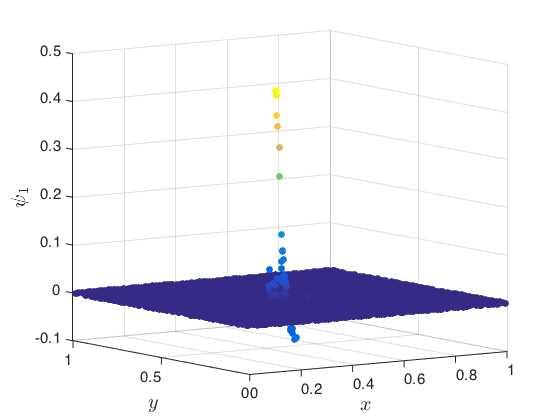}
        \caption{}
        \label{fig:psiprofile_a}
    \end{subfigure}
    \begin{subfigure}[b]{0.32\textwidth}
        \centering
        \includegraphics[width=1.0\textwidth]{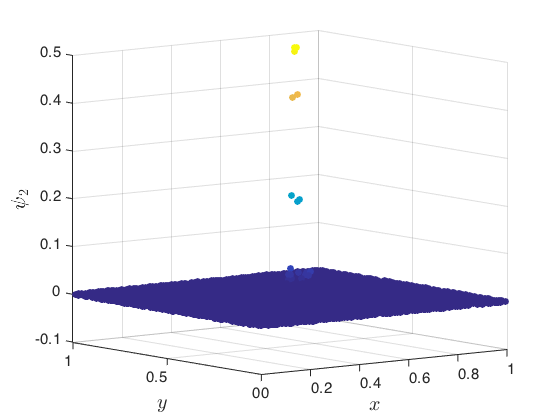}
        \caption{}
        \label{fig:psiprofile_b}
    \end{subfigure}
    \begin{subfigure}[b]{0.32\textwidth}
        \centering
        \includegraphics[width=1.0\textwidth]{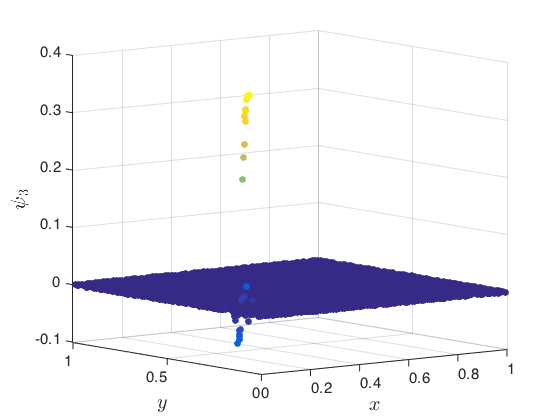}
        \caption{}
        \label{fig:psiprofile_c}
    \end{subfigure}
    
    \begin{subfigure}[b]{0.32\textwidth}
        \centering
        \includegraphics[width=1.0\textwidth]{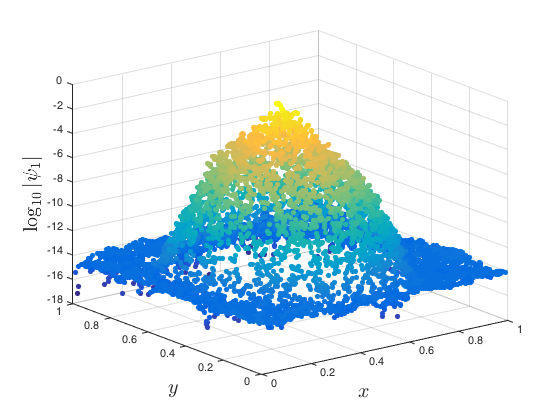}
        \caption{}
        \label{fig:psiprofile_d}
    \end{subfigure}
    \begin{subfigure}[b]{0.32\textwidth}
        \centering
        \includegraphics[width=1.0\textwidth]{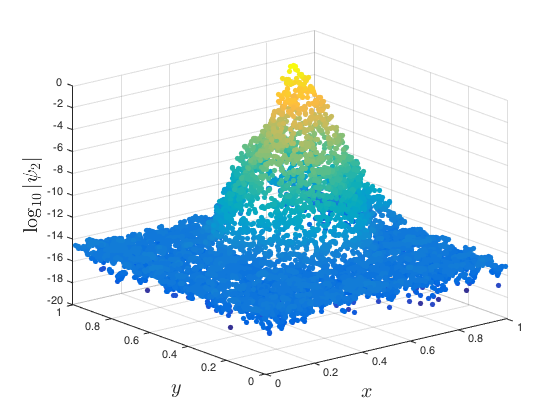}
        \caption{}
        \label{fig:psiprofile_e}
    \end{subfigure}
    \begin{subfigure}[b]{0.32\textwidth}
        \centering
        \includegraphics[width=1.0\textwidth]{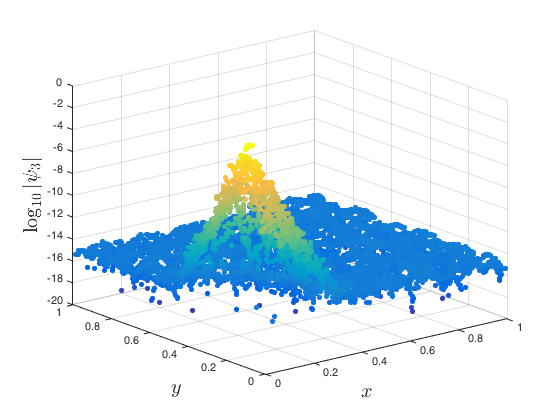}
        \caption{}
        \label{fig:psiprofile_f}
    \end{subfigure}

    \begin{subfigure}[b]{0.32\textwidth}
        \centering
        \includegraphics[width=1.0\textwidth]{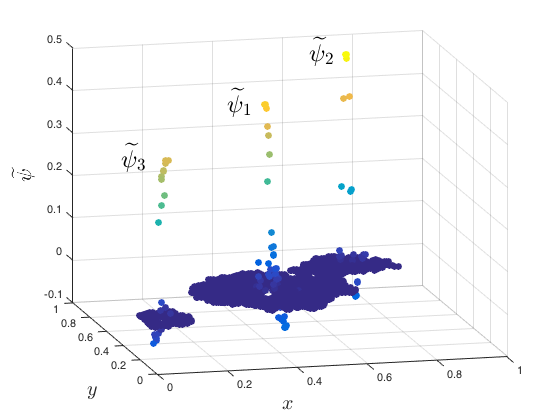}
        \caption{}
        \label{fig:psiprofile_g}
    \end{subfigure}
    \begin{subfigure}[b]{0.32\textwidth}
        \centering
        \includegraphics[width=1.0\textwidth]{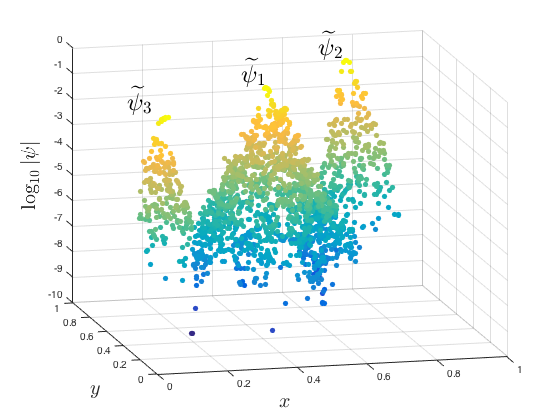}
        \caption{}
        \label{fig:psiprofile_h}
    \end{subfigure}
    \begin{subfigure}[b]{0.32\textwidth}
        \centering
        \includegraphics[width=1.0\textwidth]{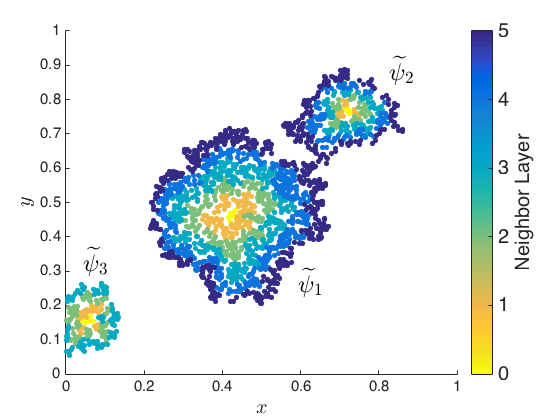}
        \caption{}
        \label{fig:psiprofile_i}
    \end{subfigure}
    \caption{Profiles of $\psi_1$, $\psi_2$, $\psi_3$, and the corresponding approximator $\tilde{\psi}_1$, $\tilde{\psi}_2$, $\tilde{\psi}_3$.}
    \label{fig:psiprofile}
\end{figure}

\begin{figure}[h!]
\centering
    \begin{subfigure}[b]{0.4\textwidth}
        \centering
        \includegraphics[width=1.0\textwidth]{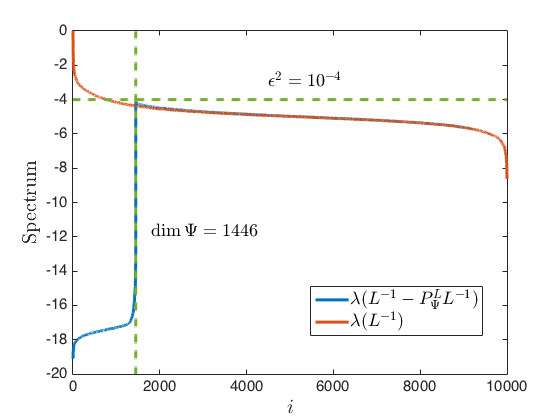}
        \caption{Spectrum of $L^{-1}$ and $L^{-1}-P_{\Psi}^LL^{-1}$}
    \end{subfigure}
    \begin{subfigure}[b]{0.4\textwidth}
        \centering
        \includegraphics[width=1.0\textwidth]{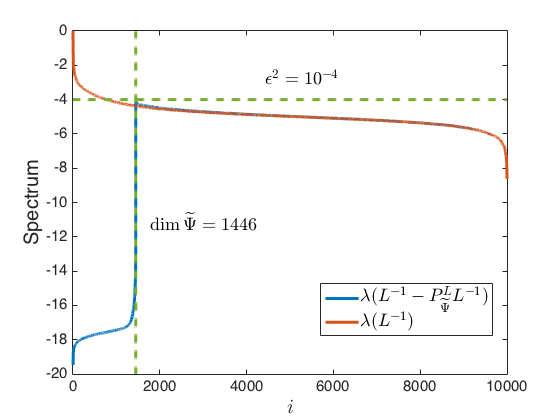}
        \caption{Spectrum of $L^{-1}$ and $L^{-1}-P_{\widetilde{\Psi}}^LL^{-1}$ }
    \end{subfigure}
    \caption{Spectrum of $L^{-1}$, $L^{-1}-P_{\Psi}^LL^{-1}$ and $L^{-1}-P_{\widetilde{\Psi}}^LL^{-1}$.}
    \label{fig:compression_spectrum}
\end{figure}

By fixing $\epsilon_{loc}^2=0.0001$, we have the mean radius of $\widetilde{\Psi} \approx 4.5$ and the mean support size $\approx 449$. We pick three functions $\psi_1,\psi_2,\psi_3$ such that $\psi_1$ is close to the center of $\Omega$, $\psi_2$ is near the boundary of connectivity change, and $\psi_3$ is close to the boundary of $\Omega$. \Cref{fig:psiprofile_a} -- \cref{fig:psiprofile_f} (first two rows of \Cref{fig:psiprofile}) show the profiles of $|\psi_i|$ and $\log_{10}|\psi_i|$ for $i=1,2,3$. Though all of them decay exponentially from their center patches to outer layers, $\psi_1$ decays slower than $\psi_3$ since the graph has a higher connectivity (i.e., larger patch sizes) near the center. \Cref{fig:psiprofile_g} and \Cref{fig:psiprofile_h} show the profiles of $|\tilde{\psi}_i|$ and $\log_{10}|\tilde{\psi}_i|$ for $i=1,2,3$. The bird's-eye view of their supports is shown in \Cref{fig:psiprofile_i}. Similarly, $\tilde{\psi}_1$ needs a larger support than $\tilde{\psi}_3$ to achieve the same accuracy, which implies that $\psi_1$ decays slower than $\psi_3$. \Cref{fig:compression_spectrum} shows the spectrum of $L^{-1}$, $L^{-1}-P_{\Psi}^LL^{-1}$ and $L^{-1}-P_{\widetilde{\Psi}}^LL^{-1}$. Notice that if we truncate the fine-scale part of $L^{-1}$ with prescribed accuracy $\epsilon^2$, then $L^{-1}-P_{\Psi}^LL^{-1}$ has rank $n-N$ and the compression error is $8.13\times 10^{-5}< \epsilon^2$. Similarly, if the local approximator $\widetilde{\psi}$ is applied (instead of $\psi$), then $L^{-1}-P_{\widetilde{\Psi}}^LL^{-1}$ also has rank $n-N$, and the compression error is also $8.13\times 10^{-5}< \epsilon^2$. This means that the same compression error is achieved by the compression $P_{\widetilde{\Psi}}^LL^{-1}$ with localization.

\begin{figure}[!h]
\centering
    \begin{subfigure}[b]{0.4\textwidth}
        \centering
        \includegraphics[width=1.0\textwidth]{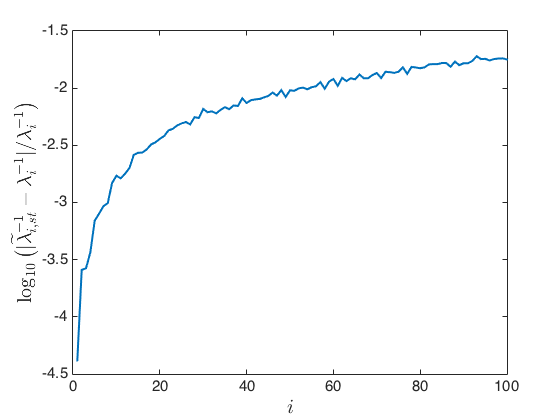}
        \caption{}
    \end{subfigure}
    \begin{subfigure}[b]{0.4\textwidth}
        \centering
        \includegraphics[width=1.0\textwidth]{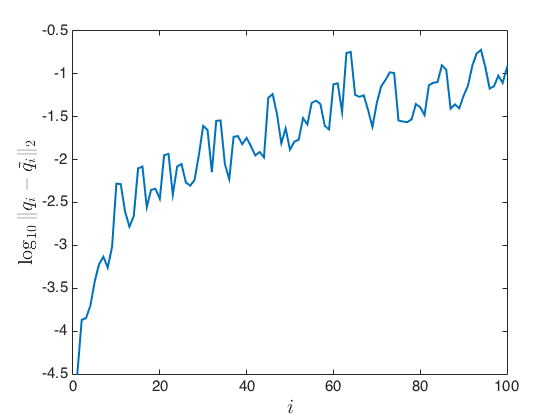}
        \caption{}
    \end{subfigure}
    \caption{Difference between the true eigenvalues/eigenvectors and the approximated ones.}
    \label{fig:eigen_difference}
\end{figure}

\begin{figure}[!h]
    \begin{subfigure}[b]{0.23\textwidth}
        \centering
        \includegraphics[width=1.0\textwidth]{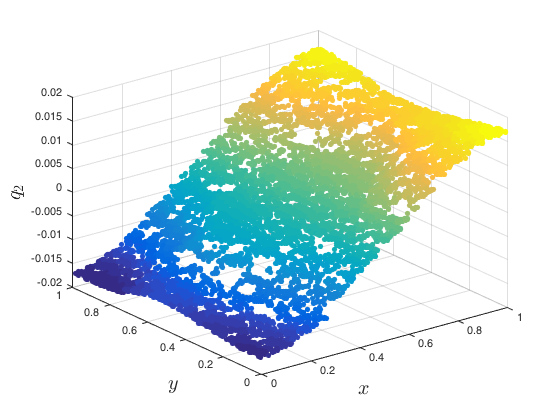}
        \caption{$q_2$}
    \end{subfigure}
    \begin{subfigure}[b]{0.23\textwidth}
        \centering
        \includegraphics[width=1.0\textwidth]{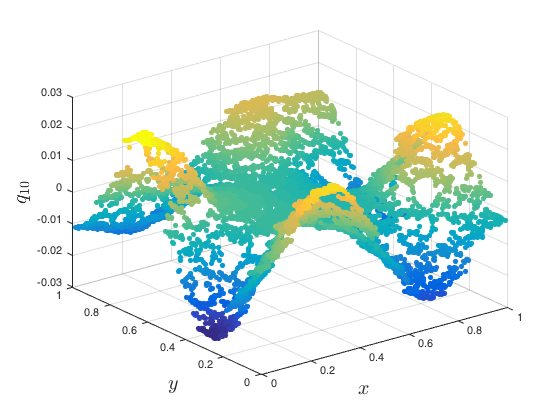}
        \caption{$q_{10}$}
    \end{subfigure}
    \begin{subfigure}[b]{0.23\textwidth}
        \centering
        \includegraphics[width=1.0\textwidth]{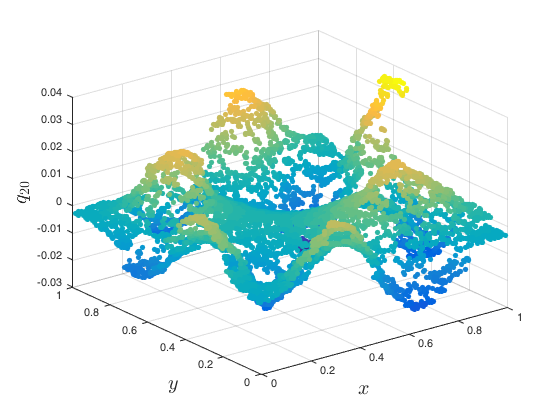}
        \caption{$q_{20}$}
    \end{subfigure}
    \begin{subfigure}[b]{0.23\textwidth}
        \centering
        \includegraphics[width=1.0\textwidth]{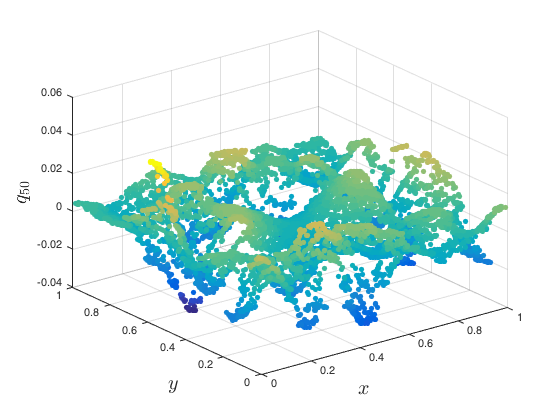}
        \caption{$q_{50}$}
    \end{subfigure}

    \begin{subfigure}[b]{0.23\textwidth}
        \centering
        \includegraphics[width=1.0\textwidth]{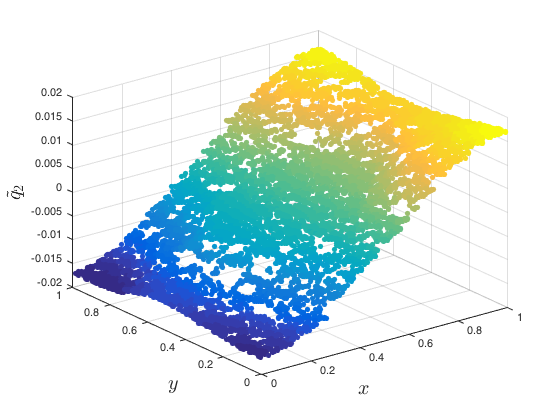}
        \caption{$\tilde{q}_{2}$}
    \end{subfigure}
    \begin{subfigure}[b]{0.23\textwidth}
        \centering
        \includegraphics[width=1.0\textwidth]{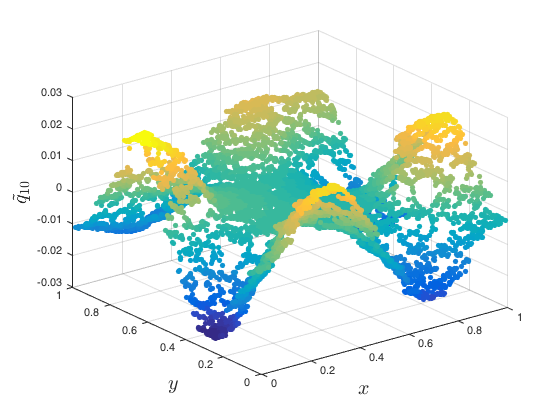}
        \caption{$\tilde{q}_{10}$}
    \end{subfigure}
    \begin{subfigure}[b]{0.23\textwidth}
        \centering
        \includegraphics[width=1.0\textwidth]{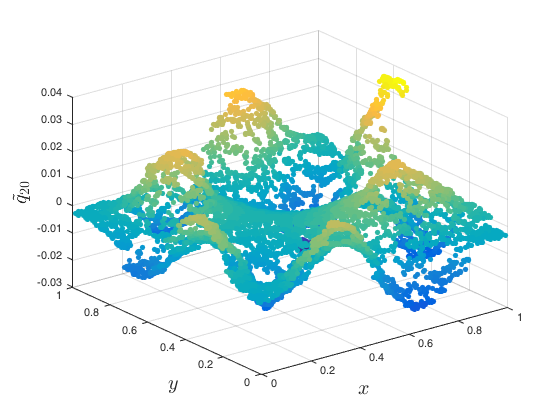}
        \caption{$\tilde{q}_{20}$}
    \end{subfigure}
    \begin{subfigure}[b]{0.23\textwidth}
        \centering
        \includegraphics[width=1.0\textwidth]{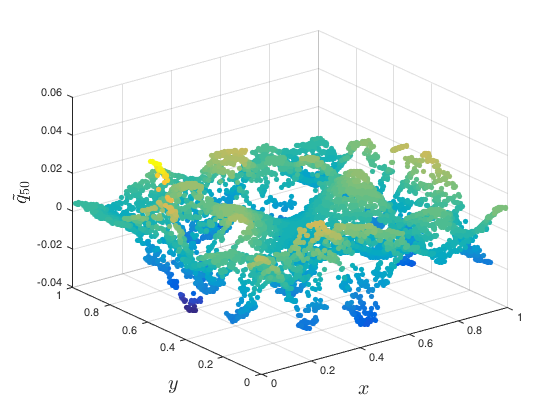}
        \caption{$\tilde{q}_{50}$}
    \end{subfigure}
    \caption{Plot of $2^{\text{nd}}$, $10^{\text{th}}$, $20^{\text{th}}$ and $50^{\text{th}}$ eigenvectors corresponding to $L^{-1}$ (first row) and $\widetilde{A}_{\text{st}}^{-1}$ (second row) respectively.}
    \label{fig:eigenvectorplot}
\end{figure}

Fourthly, we compare our results with the compression given by the PCA \cite{jolliffe2002principal}, that is $L^{-1}\approx \sum_{i=1}^{N_{\text{PCA}}}\lambda_i^{-1}q_iq_i^T$, where $\lambda_i$ is the $i_\mathrm{th}$ smallest eigenvalue of $L$ and $q_i$ is the corresponding normalized eigenvector. On the one hand, to achieve the same compression error $8.13\times 10^{-5}$, we need $N_{\text{PCA}}=893$, which is the optimal compression dimension for such accuracy. But we remark that such achievement requires solving global eigen problem and the compressed operator $\sum_{i=1}^{N_{pca}}\lambda_i^{-1}q_iq_i^T$ usually loses the original sparsity features. On the other hand, our approach has a larger number of basis functions ($N=1446$) since only local eigen information is used and patches are combined pair-wisely. But our approach gives local functions with compressed dimension just up to 2 times of the optimal dimension. Further, it turns out that we can recover the eigenvectors of $L$ corresponding to relatively small eigenvalues by solving eigenvalue problem of $A_{\text{st}}$ (or $\widetilde{A}_{\text{st}}$). \Cref{fig:eigenvectorplot} shows the $2^{\text{nd}}$, $10^{\text{th}}$, $20^{\text{th}}$ and $50^{\text{th}}$ eigenvectors corresponding to the small eigenvalues of $L^{-1}$ (first row) and $\widetilde{A}_{\text{st}}^{-1}$ (second row) respectively. Let $\tilde{\lambda}_{i,st}$ be the $i_\mathrm{th}$ smallest eigenvalue of $\widetilde{A}_{\text{st}}$, and $\tilde{\xi}_i$ be the corresponding eigenvector so that $\tilde{q}_i=\widetilde{\Psi}\tilde{\xi}_i$ has $l_2$-norm equal to 1. From the experiment, we observe that $\tilde{\lambda}_{i,st}^{-1}$ is a good approximation of $\lambda_i^{-1}$ and $\tilde{q}_i$ is a good approximation of $q_i$ for small $\lambda_i$, as shown in \Cref{fig:eigen_difference}. In other words, this procedure provides us convenience for computing the first few eigenvalues and eigenvectors of $L$, since $\widetilde{A}_{\text{st}}$ has the compressed size with a much smaller condition number ($\kappa(\widetilde{A}_{\text{st}})=1.14\times10^{5}$ vs. $\kappa(L)=4.29\times10^{8}$).


\subsection{Numerical Example 2}
\label{subsec:numerical_example2}
\begin{figure}[!h]
\centering
    \begin{subfigure}[b]{0.3\textwidth}
        \centering
        \includegraphics[width=1.0\textwidth]{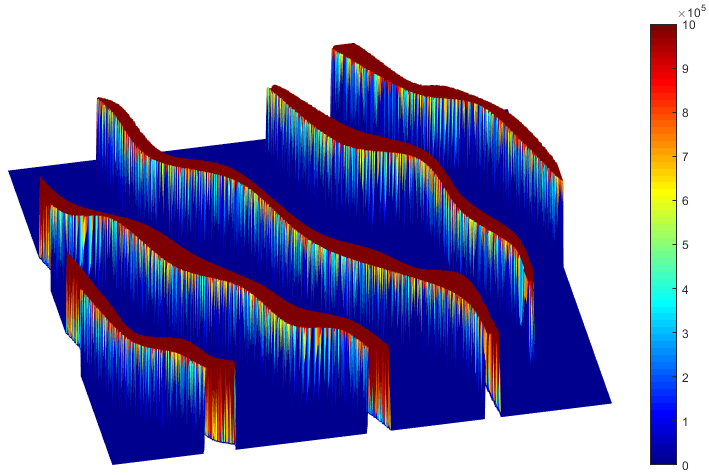}
        \caption{Distribution of $\mu(x,y)$}
        \label{fig:contrastdistribution}
    \end{subfigure}
    \begin{subfigure}[b]{0.21\textwidth}
        \centering
        \includegraphics[width=1.0\textwidth]{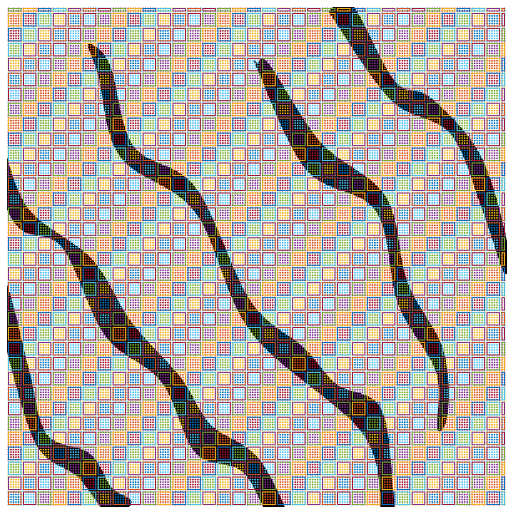}
        \caption{$200 \times 200$ regular}
        \label{fig:contrast2_regular}
    \end{subfigure}
    \begin{subfigure}[b]{0.21\textwidth}
        \centering
        \includegraphics[width=1.0\textwidth]{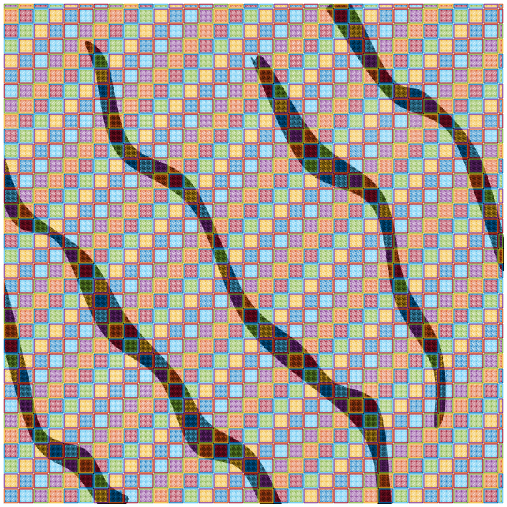}
        \caption{$400 \times 400$ regular}
        \label{fig:contrast_regular}
    \end{subfigure}
       \begin{subfigure}[b]{0.3\textwidth}
        \centering
        \includegraphics[width=1.0\textwidth]{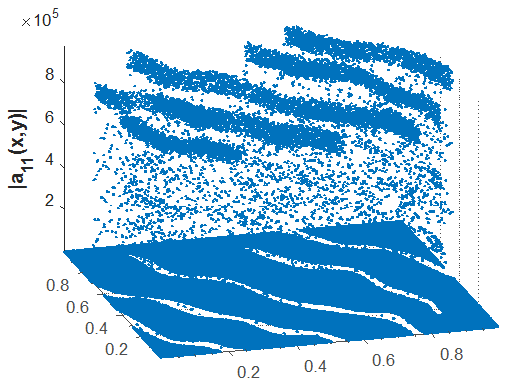}
        \caption{$|a_{11}(x,y)|$ for $(x,y) \in \Omega$}
        \label{fig:a1_distribution}
    \end{subfigure}
    \begin{subfigure}[b]{0.21\textwidth}
        \centering
        \includegraphics[width=1.0\textwidth]{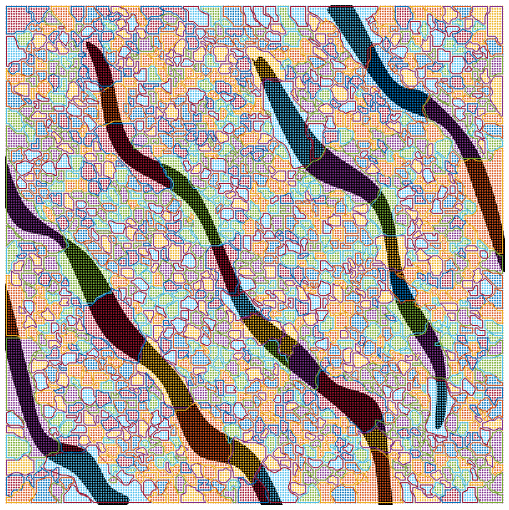}
        \caption{$200 \times 200$ proposed}
        \label{fig:contrast2_ourpartition}
    \end{subfigure}
    \begin{subfigure}[b]{0.21\textwidth}
        \centering
        \includegraphics[width=1.0\textwidth]{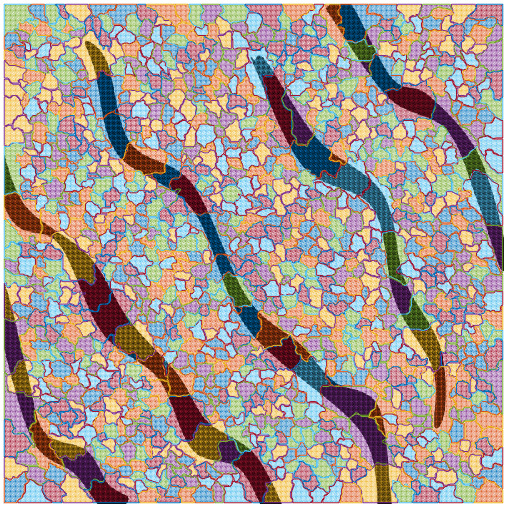}
        \caption{$400 \times 400$ proposed}
        \label{fig:contrast_ourpartition}
    \end{subfigure}
    \caption{Partitioning result: operator stemmed from FEM of Elliptic PDE with High contrasted coefficients. (a) shows the distribution of the high contrast factor in the domain $\Omega$. (b) and (c) shows the regular partition in resolution $200 \times 200$ and $400 \times 400$ respectively. (d) shows the value of $a_{11}(x,y)$ on $\Omega$. (e) and (f) show the partition results obtained from \Cref{alg:pair_clustering}.}
    \label{fig:contrast}
\end{figure}

Our second numerical example arises from using GFEM to solve the following elliptic equation with homogeneous Dirichlet boundary conditions:
\begin{equation}
-\nabla\cdot(a\cdot\nabla u)=f,\quad u\in H_0^1(\Omega),
\label{eqt:numerical2_problem}
\end{equation}
where $\Omega=[0,1]\times [0,1]$, $f\in L^2(\Omega)$, and the coefficient $a$ is a 2-by-2 matrix function of $(x,y)$ of the form 
\begin{equation}
a=\left(\begin{array}{cc}a_{11} & a_{12}\\ a_{12} & a_{22}\end{array}\right)=
\left(\begin{array}{cc} \cos\theta & \sin\theta \\ -\sin\theta & \cos\theta\end{array}\right)
\left(\begin{array}{cc} \mu\cdot e_1 & 0 \\ 0 & \mu\cdot e_2\end{array}\right)
\left(\begin{array}{cc} \cos\theta & -\sin\theta \\ \sin\theta & \cos\theta\end{array}\right).
\label{eqt:definition_coefficient}
\end{equation}
Here $\theta=\theta(x,y)\in C(\Omega)$ is the rotation (deformation) factor, $\mu=\mu(x,y)\in L^\infty(\Omega)$ is the contrast factor, and $e_i=e_i(x,y)\in L^\infty(\Omega),\ i=1,2$ are the roughness factor. For the problem to be elliptic and well-posed, we require that $\mu e_1,\mu e_2>C$ for some uniform constant $C>0$. To increase the level of difficulty in solving this elliptic PDE, we choose $e_1$ and $e_2$ to be highly oscillatory and $\mu$ varies from $O(1)$ to $O(10^6)$ (high contrasts). More precisely, $e_1,e_2$ are generated with extreme roughness as $e_i(x,y)=1+w_i(x,y),\ i=1,2$, where for each point $(x,y)$, we set $w_1(x,y), w_2(x,y)\overset{i.i.d}{\sim}\mathcal{U}([-0.1,0.1])$, a uniform distribution on $[-0.1,0.1]$. The contrast factor $\mu(x,y)$ is generated from the background permeability as shown in \Cref{fig:contrastdistribution}. $\theta$ is given by $\theta(x,y)= \pi \cdot (x + y)$. The magnitude of $|a_{11}(x,y)|$ in $\Omega$ is also plotted in \Cref{fig:a1_distribution} as a reference.

We use GFEM \cite{brenner2004finite} with a regular triangularization to form a finite system that is fine enough to capture the details of the background field. The basis $\mathcal{V}\in\mathbb{R}^n$ is the vector representation of the Galerkin nodal basis and the SPD matrix $A$ is the stiffness matrix of nodal basis with respect to the energy inner product $\int_{\Omega}(\nabla\cdot)^T a (\nabla\cdot) dxdy$. The energy decomposition $\mathcal{E}=\{E_k\}_{k=1}^m$ is the collection of all patch-wise stiffness matrices $E_k$ on every triangle $\tau_k$. Specifically, each $E_k$ has the form  
\begin{equation*}
E_k=\left( \begin{smallmatrix}
0 & & & & \\
& w_{1_k,2_k} &  w_{1_k,2_k} & w_{1_k,3_k} & \\
& w_{2_k,1_k} &  w_{2_k,2_k} & w_{2_k,3_k} & \\
& w_{3_k,1_k} &  w_{3_k,2_k} & w_{3_k,3_k} & \\
& & & & 0\\
\end{smallmatrix}\right), \quad w_{i_k,j_k}=\int_{\tau_k} (\nabla \phi_{i_k})^T a (\nabla \phi_{j_k})dxdy,\ i,j=1,2,3
\end{equation*}
where $\phi_{i_k},\ i=1,2,3 $ are the three nodal basis surrounding $\tau_k$. In this case, every finest energy element involves three functions, which generalizes the concept of graphs' edges as we mentioned before. One should also notice that for patch $\tau_k$ touching the boundary of $\Omega$, the corresponding $E_k$ reduces to involve only two or one function since nodal basis functions on boundary are not required for homogeneous Dirichlet boundary conditions. Moreover, we will choose $q = 1$ since the problem \cref{eqt:numerical2_problem} is of second order (i.e. $p = 1$).

\begin{table}[!h]
\centering
\footnotesize
\renewcommand{\arraystretch}{1.5}
    \begin{tabular}{|c|c|c|c|c|c|c|c|}
    \hline
    & Partition & $\#\mathcal{P}$ & $\varepsilon(\mathcal{P},1)^2$ & $\delta(\mathcal{P},1)$ & $\kappa(A_{\text{st}})$ & $\delta(\mathcal{P},1) \Vert A^{-1} \Vert_2$ & $\max\delta(P_j,1) \varepsilon(P_j,1)^2$ \\ \hline
\multirow{ 2}{*}{$200 \times 200$} 
& \text{Ours} 		 & 1396 & $9.9932 \times 10^{-5}$ & $3.8703 \times 10^{6}$ & $4.8023 \times 10^{3}$ & $1.1462 \times 10^{4}$ & 233.3080 \\ \cline{2-8}
& \text{Regular} 	 & 1156 & $5.7733 \times 10^{-5}$ & $2.5299 \times 10^{11}$ & $2.5325 \times 10^{8}$ & $7.4924 \times 10^8$ & 454.6966 \\ \hline
\multirow{ 2}{*}{$400 \times 400$} 
& \text{Ours} 		 & 1199 & $9.9958 \times 10^{-5}$ & $4.0467 \times 10^{6}$ & $3.1195 \times 10^{3}$ & $1.1073 \times 10^{4}$ & 277.6022 \\ \cline{2-8}
& \text{Regular} 	 & 1156 & $5.0450 \times 10^{-5}$ & $3.1410 \times 10^{11}$ & $2.2573 \times 10^{8}$ & $8.5947 \times 10^8$ & $1.0560 \times 10^3$ \\ \hline
    \end{tabular}
    \caption{Comparison with the uniform regular partitioning for $200 \times 200$ and $400 \times 400$ resolution.}
    \label{table:example5_2}
\end{table}

In this example, we compare our partitioning technique with the performance of regular partition to illustrate the adaptivity of our algorithm to the features of the coefficient function $a(x,y)$. Here we consider the cases with two different resolutions, which are the regular triangulations with $200 \times 200$ and $400 \times 400$ vertices respectively. We set $q = 1$, the prescribed accuracy $\epsilon^2 = 10^{-4}$ and the upper bound $c = 300$. We apply \Cref{alg:operator_compression} to compute the compressed operator for both cases. In the case of regular partition, we use a uniform and regular partition on $\Omega$ that achieves the prescribed accuracy (i.e., $\varepsilon(\mathcal{P},1)^2 < 10^{-4}$). Notice that in this case, the first eigenvector, $\Phi_j$ on every patch is a constant. Such choice of $\Phi$, along with the use of regular partition is equivalent to the set up in \cite{owhadi2017multigrid,hou2016sparse}. Therefore, the only modification we apply in this numerical example is the adaptive construction of the construction $\mathcal{P}$, which remarkably improves the behavior of the compressed operator $A^{-1}$. We would also like to remark that in the case without high-contrasted channels, \Cref{alg:pair_clustering} coherently gives a regular partitioning on the domain as conventional partition methods.

\Cref{table:example5_2} summarizes the partitioning results in both cases. Under regular partitioning, we have $\# \mathcal{P} = 1156$ and each $P_j$ has the size of at most $6 \times 6$ (\Cref{fig:contrast2_regular}) and $12 \times 12$ (\Cref{fig:contrast_regular}) vertices respectively. Notably, the {\bf condition factors} for both cases go up to $10^{11}$ and the corresponding true condition numbers $\kappa(A_{\text{st}})$ are having an order of $10^8$, which show that such partition will produce an ill-posed compressed operator. Using our approach, the square {\bf error factors} $\varepsilon(\mathcal{P},1)^2$ achieved in both cases are strictly bounded above by the prescribed accuracy $\epsilon^2 = 10^{-4}$ and $\delta(\mathcal{P},1)$ is of the order of $10^6$ only. Indeed, the true compression error is even smaller as the square {\bf error factor} $\varepsilon(\mathcal{P},1)^2$ is only the theoretical upper bound as required in \Cref{prop:localphi}. Furthermore, the true condition numbers $\kappa(A_{\text{st}})$ are in the order of $10^3$ (compare to $10^8$ from regular partitioning), which are again bounded by (and is much smaller than) $\delta(\mathcal{P},1) \| A^{-1} \|_2$ as observed in \Cref{thm:condition_number}. Moreover, the patch number $\# \mathcal{P}$ is comparable to the case of regular partitioning. Also notice that $\max_{P_j \in \mathcal{P}} \delta(P_j,1) \cdot \varepsilon(P_j,1)^2 < c = 300$ in both cases, which are coherent to the prescribed requirement. These results successfully illustrate the consistency between the numerical results and theoretical discoveries in the previous sections. The partition results obtained by \Cref{alg:pair_clustering} are shown in \Cref{fig:contrast2_ourpartition} and \Cref{fig:contrast_ourpartition} respectively.

\begin{figure}[!h]
\centering
    \begin{subfigure}[b]{0.32\textwidth}
        \centering
        \includegraphics[width=1.0\textwidth]{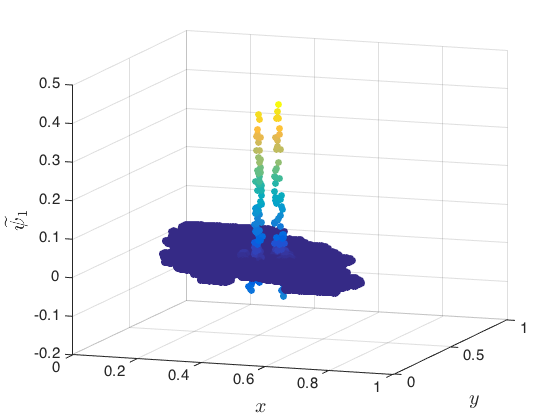}
        \caption{}
        \label{fig:psiplot_a}
    \end{subfigure}
    \begin{subfigure}[b]{0.32\textwidth}
        \centering
        \includegraphics[width=1.0\textwidth]{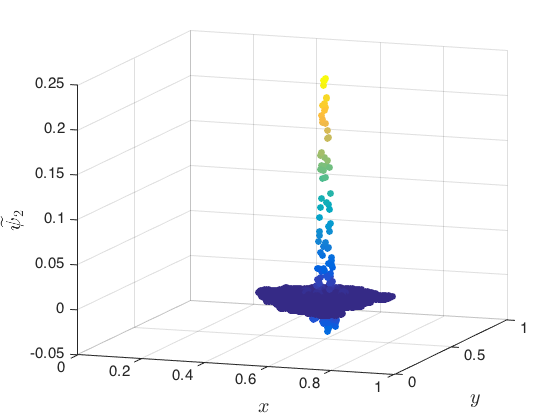}
        \caption{}
        \label{fig:psiplot_b}
    \end{subfigure}
    \begin{subfigure}[b]{0.32\textwidth}
        \centering
        \includegraphics[width=1.0\textwidth]{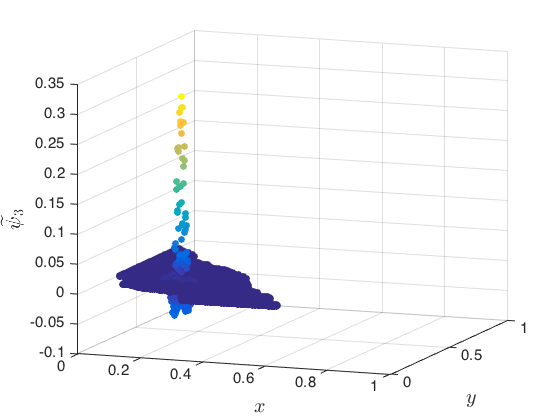}
        \caption{}
        \label{fig:psiplot_c}
    \end{subfigure}
    \begin{subfigure}[b]{0.32\textwidth}
        \centering
        \includegraphics[width=1.0\textwidth]{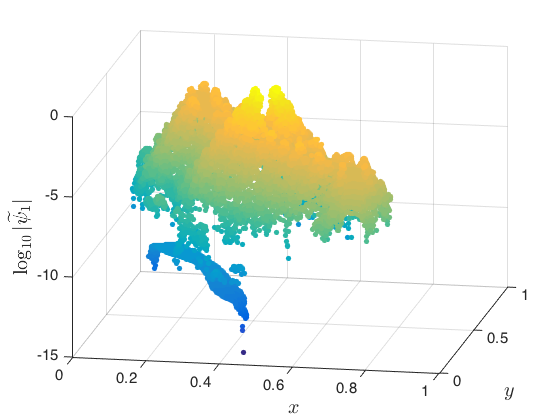}
        \caption{}
        \label{fig:psiplot_d}
    \end{subfigure}
    \begin{subfigure}[b]{0.32\textwidth}
        \centering
        \includegraphics[width=1.0\textwidth]{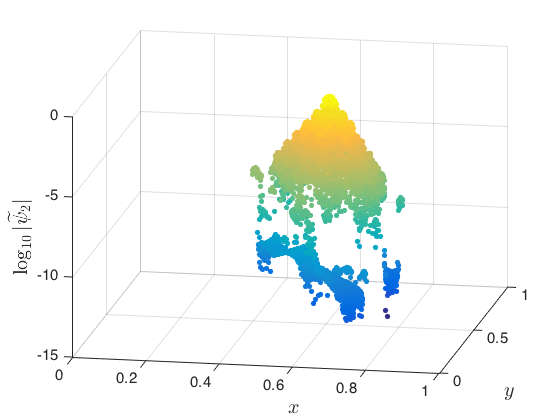}
        \caption{}
        \label{fig:psiplot_e}
    \end{subfigure}
    \begin{subfigure}[b]{0.32\textwidth}
        \centering
        \includegraphics[width=1.0\textwidth]{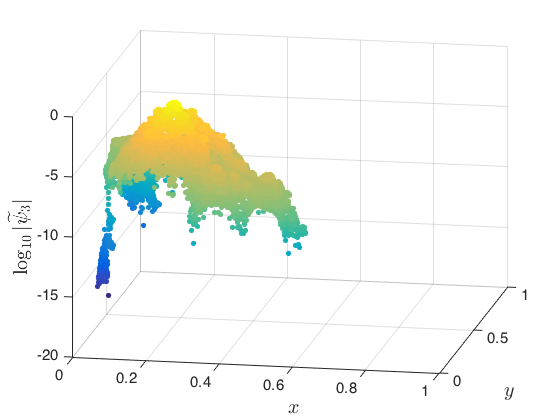}
        \caption{}
        \label{fig:psiplot_f}
    \end{subfigure}
    \begin{subfigure}[b]{0.32\textwidth}
        \centering
        \includegraphics[width=1.0\textwidth]{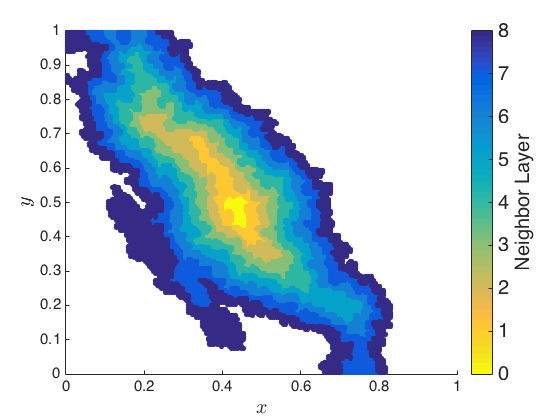}
        \caption{}
        \label{fig:psiplot_g}
    \end{subfigure}
    \begin{subfigure}[b]{0.32\textwidth}
        \centering
        \includegraphics[width=1.0\textwidth]{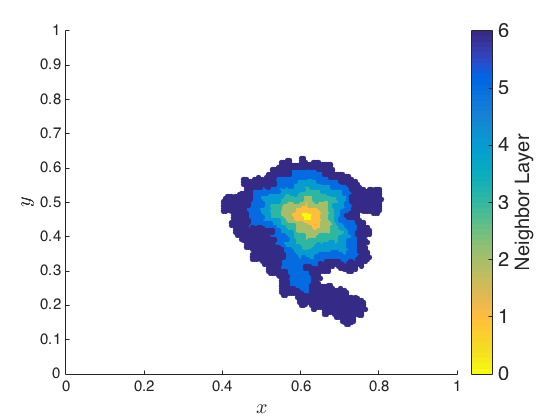}
        \caption{}
        \label{fig:psiplot_h}
    \end{subfigure}
    \begin{subfigure}[b]{0.32\textwidth}
        \centering
        \includegraphics[width=1.0\textwidth]{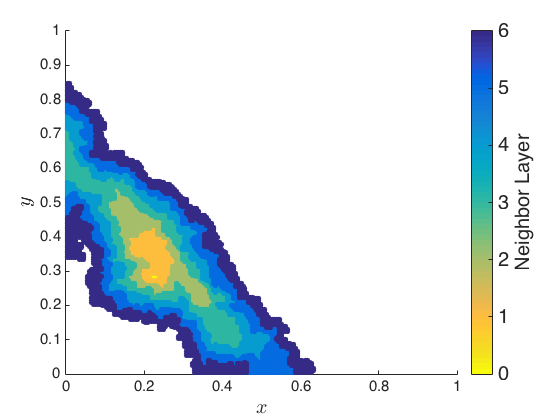}
        \caption{}
        \label{fig:psiplot_i}
    \end{subfigure}
    \caption{Samples of the localized basis.}
    \label{fig:psiplot}
\end{figure}

We remark that the reason for the huge difference in the condition number is caused by the non-adaptivity of the partitioning (i.e., regular partition) to the given operator $A$. Specifically, if some patches are fully covered in the high channel regions, the accuracy achieved by this patch is obviously very promising. However, corresponding patch-wise {\bf condition factor} will jump up to the similar order as the high contrast factor. In other words, patches which are fully covered by regions of high contrast should be avoid. As shown in \Cref{fig:contrast_ourpartition} and \Cref{fig:contrast2_ourpartition}, our proposed \Cref{alg:pair_clustering} can automatically extract the intrinsic geometric information of the operator (which is the distribution of high contrast regions) and prevent patches which are fully enclosed in the high contrast regions.

\begin{figure}[!h]
\centering
    \begin{subfigure}[b]{0.32\textwidth}
        \centering
        \includegraphics[width=1.0\textwidth]{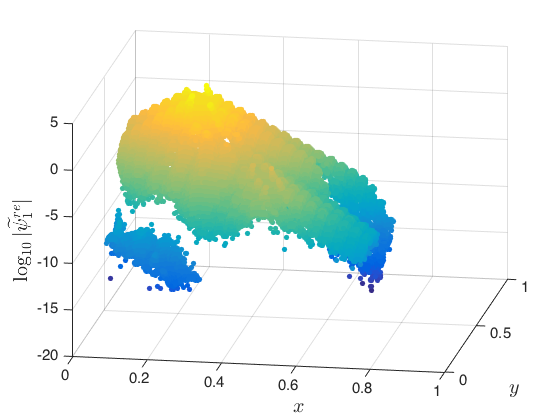}
        \caption{}
        \label{fig:repsiplot_d}
    \end{subfigure}
    \begin{subfigure}[b]{0.32\textwidth}
        \centering
        \includegraphics[width=1.0\textwidth]{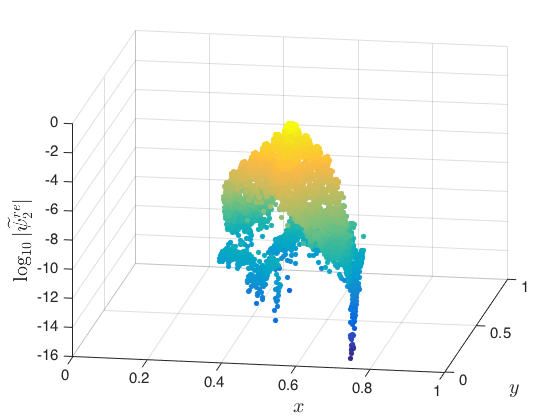}
        \caption{}
        \label{fig:repsiplot_e}
    \end{subfigure}
    \begin{subfigure}[b]{0.32\textwidth}
        \centering
        \includegraphics[width=1.0\textwidth]{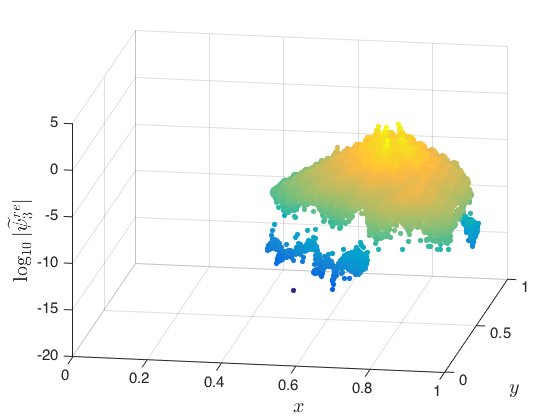}
        \caption{}
        \label{fig:repsiplot_f}
    \end{subfigure}
    \begin{subfigure}[b]{0.32\textwidth}
        \centering
        \includegraphics[width=1.0\textwidth]{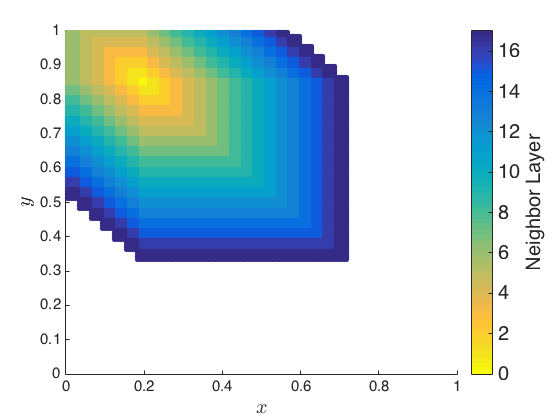}
        \caption{}
        \label{fig:repsiplot_g}
    \end{subfigure}
    \begin{subfigure}[b]{0.32\textwidth}
        \centering
        \includegraphics[width=1.0\textwidth]{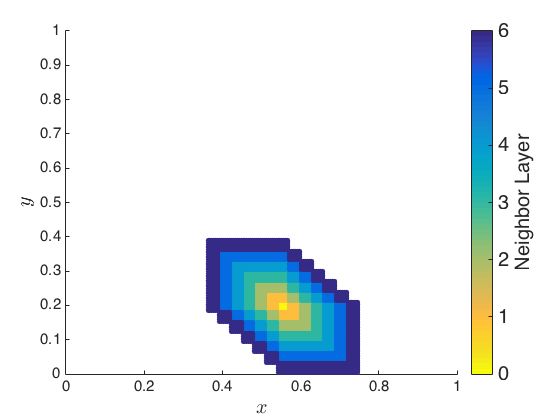}
        \caption{}
        \label{fig:repsiplot_h}
    \end{subfigure}
    \begin{subfigure}[b]{0.32\textwidth}
        \centering
        \includegraphics[width=1.0\textwidth]{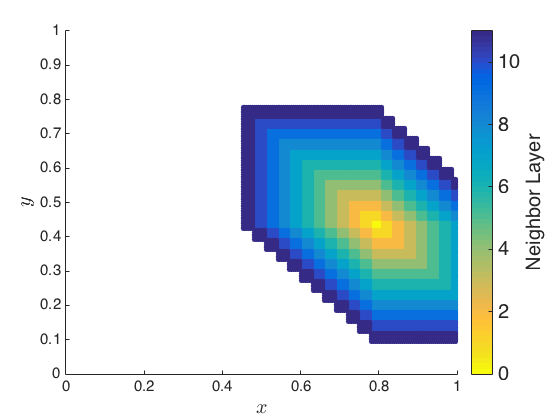}
        \caption{}
        \label{fig:repsiplot_i}
    \end{subfigure}
    \caption{Samples of the localized basis from regular partition}
    \label{fig:repsiplot}
\end{figure}

We also consider the profiles, log-profiles and supports of three localized functions $\tilde{\psi}_1,\tilde{\psi}_2,\tilde{\psi}_3$ obtained by \Cref{alg:construct_tilde_psi} with prescribed localization error $\epsilon_{loc}^2=10^{-4}$ in the $200 \times 200$ resolution case. The plots are shown in \Cref{fig:psiplot}. We can see in \Cref{fig:psiplot_a} that $\tilde{\psi}_1$ has multiple peaks (three peaks exactly), with one high contrast channel cutting through. This means that this single function characterizes the local feature of the operator. Though the exponential decaying feature is still obvious, to achieve the prescribed localization error, some local functions ($\widetilde{\psi}_1,\widetilde{\psi}_3$) have to extend along the high permeability channels and thus end up with relatively large supports. Recall that \Cref{rmk:alpha_remark} implies that the decay rate of $\psi_i$ and thus the radius of $\tilde{\psi}_i$ are invariant under local scaling (contrast scaling). But higher permeability means stronger connectivity and consequently larger patch sizes, and therefore $\tilde{\psi}_i$ extends farther in physical distance along high permeability channels. As a limitation of our approach, this long range decaying compromises the sparsity of the localized basis and the stiffness matrix, which is an issue that we plan to resolve in our future work. As comparison, \Cref{fig:repsiplot} plots the log-profile and the supports of three localized functions $\tilde{\psi}^{re}_1,\tilde{\psi}^{re}_2,\tilde{\psi}^{re}_3$ obtained similarly but with the regular partition in the $200 \times 200$ resolution case. We can see that some localized functions under regular partition also have relative large support size. However, such long distance extension is not the result of large patch sizes (since regular partition has uniform patch size), but the result of large {\bf condition factor} $\delta(\mathcal{P},1)$. Recall that to achieve an desired localization error, the radius of localized basis is also affected by $\log\delta(\mathcal{P},1)$.

\section{Multiresolution operator decomposition}
\label{Sec:Multiresolution}
One can see that what we have been doing with operator compression is essentially to truncate the microscopic/fine-scale part of $A$ while preserving the necessary macroscopic/coarse-scale part that dominates the accuracy. And in the meanwhile, the condition number of the compressed operator also drops to the level consistent to the prescribed accuracy. This consistency inspires us to perform the compression procedure hierarchically in order to separate the operator into multiple scales of resolution rather than just two.

Now consider that, instead of just compressing the inverse $A^{-1}$, we want to solve the problem $Ax=b$. Due to the sparsity of $A$, a straightforward idea is to employ iterative methods. But these methods will suffer from the large condition number of $A$. Alternatively, we would like to use the energy decomposition of $A$, and the locality (i.e. sparsity) of the energy decomposition to resolve the difficulty of large condition number. The main idea is to decompose the computation of $A^{-1}$ into hierarchical resolutions such that (i)the relative condition number in each scale/level can be well bounded, and (ii)the sub-system to be solved on each level is as sparse as the original $A$. Here we will make use of the choice of the partition $\mathcal{P}$, the basis $\Phi$ and $\Psi$ obtained in \Cref{sec:operator_compression} and \Cref{Sec:partition} to serve the purpose of multiresolution operator decomposition. In the following, We first implement a {\bf one-level decomposition}.

\subsection{One-level operator decomposition} Let $\mathcal{P}$, $\Phi$, $\Psi$ and $U$ be constructed as in \Cref{alg:operator_compression}, namely

\medskip
\begin{enumerate}[label=(\roman*)]
\item $\mathcal{P}=\{P_j\}_{j=1}^M$ is a partition of $\mathcal{V}$;
\item $\Phi=[\Phi_1,\Phi_2,\cdots,\Phi_M]$ such that every $\Phi_j\subset \Span\{P_j\}$ has dimension $q_j$;
\item $\Psi=A^{-1}\Phi(\Phi^TA^{-1}\Phi)^{-1}$;
\item $U=[U_1,U_2,\cdots,U_M]$ such that every $U_j\subset \Span\{P_j\}$ has dimension $\dim(\Span\{P_j\})-q_j$ and satisfies $\Phi_j^TU_j= {\bf 0}$ .
\end{enumerate}
\medskip

\noindent Then $[U,\Psi]$ forms a basis of $\mathbb{R}^n$, and we have
\begin{equation}
U^TA\Psi=U^T\Phi(\Phi^TA^{-1}\Phi)^{-1}=0.
\end{equation}
Thus the inverse of $A$ can be written as
\begin{equation}
\begin{split}
A^{-1}=&\ \Big(\left[\begin{array}{c} U^T\\ \Psi^T\end{array}\right]^{-1}\left[\begin{array}{c} U^T\\ \Psi^T\end{array}\right]A\left[\begin{array}{cc} U & \Psi\end{array}\right]\left[\begin{array}{cc} U & \Psi\end{array}\right]^{-1}\Big)^{-1} \\
=&\ U(U^TAU)^{-1}U^T +\Psi(\Psi^TA\Psi)^{-1}\Psi^T.
\end{split}
\label{eqt:Ainv_representation}
\end{equation}
In the following we denote $\Psi^TA\Psi=A_{\text{st}}$ and $U^TAU=B_{\text{st}}$ respectively. We also use the phrase ``solving $A^{-1}$'' to mean ``solving $A^{-1} b$ for any $b$". From \Cref{eqt:Ainv_representation}, we observe that solving $A^{-1}$ is equivalent to solving $A_{\text{st}}^{-1}$ and $B_{\text{st}}^{-1}$ separately. For $B_{\text{st}}$, notice that since the space/basis $U$ is constructed locally with respect to each patch $P_j$, $B_{\text{st}}$ will inherit the sparsity characteristic from $A$ if $A$ is local/sparse. Thus it will be efficient to solve $B_{\text{st}}^{-1}$ using iterative type methods if the condition number of $B_{\text{st}}$ is bounded. In the following, we introduce \Cref{lemma:Bst_eigenvalue} which provides an upper bounded of the $B_{\text{st}}$ that ensures the efficiency of solving $B_{\text{st}}^{-1}$. The proof of the lemma imitates the proof from Theorem 10.9 of \cite{owhadi2017universal}, where the required condition \cref{eqt:phi_constraint} corresponds to Equation (2.3) in \cite{owhadi2017universal}.

\begin{lemma}
If $\Phi$ satisfies the condition \cref{eqt:phi_constraint} with constant $\epsilon$, then 
\begin{equation}
\lambda_{\max}(B_{\text{st}})\leq\lambda_{\max}(A)\cdot\lambda_{\max}(U^T U),\qquad \lambda_{\min}(B_{\text{st}})\geq \frac{1}{\epsilon^2}\cdot\lambda_{\min}(U^T U),
\end{equation}
and thus 
\begin{equation}
\kappa(B_{\text{st}})\leq \epsilon^2\cdot\lambda_{\max}(A)\cdot\kappa(U^TU).
\end{equation}
\label{lemma:Bst_eigenvalue}
\end{lemma}

\begin{proof} For $\lambda_{\max}(B_{\text{st}})$, we have
\begin{equation}
\lambda_{\max}(B_{\text{st}}) = \|B_{\text{st}}\|_2 = \|U^TAU\|_2\leq \|A\|_2\|U\|_2^2 = \|A\|_2\|U^TU\|_2 =\lambda_{\max}(A)\lambda_{\max}(U^T U).
\end{equation}
For $\lambda_{\min}(B_{\text{st}})$, since $\Phi$ satisfies the condition \cref{eqt:phi_constraint} with constant $\epsilon$  and $\Phi^TU = {\bf 0}$, we have
\[\|x\|_2^2\leq \frac{1}{\lambda_{\min}(U^TU)} x^TU^TUx\leq \frac{ \epsilon^2 }{\lambda_{\min}(U^TU)} x^TU^TAUx=\frac{ \epsilon^2 }{\lambda_{\min}(U^TU)} x^TB_{\text{st}}x,\]
thus $\lambda_{\min}(B_{\text{st}})\geq \frac{1}{\epsilon^2}\lambda_{\min}(U^T U)$.
\end{proof}

We shall have some discussion on the bound $\epsilon^2\cdot\lambda_{\max}(A)\cdot\kappa(U^TU)$ separately into two parts, namely (i) $\epsilon^2\cdot\lambda_{\max}(A)$; and (ii) $\kappa(U^TU)$. Notice that $U^TU$ is actually block-diagonal with blocks $U_j^TU_j$, therefore 
\begin{equation}
\kappa(U^TU) = \frac{\lambda_{\max}(U^TU)}{\lambda_{\min}(U^TU)} = \frac{\max_{1\leq j\leq M}{\lambda_{\max}(U_j^T U_j)}}{\min_{1\leq j\leq M}{\lambda_{\min}(U_j^T U_j)}}.
\end{equation}
In other words, we can bound $\kappa(U^TU)$ well by choosing proper $U_j$ for each $P_j$. For instance, if we allow any kind of local computation on $P_j$, we may simply extend $\Phi_j$ to an orthonormal basis of $\Span\{P_j\}$ to get $U_j$ by using QR factorization \cite{saad2011numerical}. In this case, we have $\kappa(U^TU)=1$. 

For the part $\epsilon^2\lambda_{\max}(A)$, recall that we construct $\Phi$ based on a partition $\mathcal{P}$ and an integer $q$ so that $\Phi$ satisfies condition \cref{eqt:phi_constraint} with constant $\varepsilon(\mathcal{P},q)$, thus the 
posterior bound of $\kappa(B_{\text{st}})$ is $\varepsilon(\mathcal{P},q)^2\lambda_{\max}(A)$(when $\kappa(U^TU)=1$). Recall that $\kappa(A_{\text{st}})$ is bounded by $\delta(\mathcal{P},q)\|A^{-1}\|_2$, therefore $\kappa(A_{\text{st}})\kappa(B_{\text{st}})\leq \varepsilon(\mathcal{P},q)^2\delta(\mathcal{P},q)\kappa(A)$. That is, $A_{\text{st}}$ and $B_{\text{st}}$ divide the burden of the large condition number of $A$ with an amplification factor $\varepsilon(\mathcal{P},q)^2\delta(\mathcal{P},q)$. We call $\kappa(\mathcal{P},q)\triangleq\varepsilon(\mathcal{P},q)^2\delta(\mathcal{P},q)$ the $q_\mathrm{th}$-order condition number of the partition $\mathcal{P}$. This explains why we attempt to bound $\kappa(\mathcal{P},q)$ in the construction of the partition $\mathcal{P}$.

Ideally, we hope the one-level operator decomposition gives $\kappa(A_{\text{st}})\approx\kappa(B_{\text{st}})$, so that the two parts equally share the burden in parallel. But such result may not be good enough when $\kappa(A_{\text{st}})$ and $\kappa(B_{\text{st}})$ are still large. To fully decompose the large condition number of $A$, a simple idea is to recursively apply the one-level decomposition. That is, we first set a small enough $\epsilon$ to sufficiently bound $\kappa(B_{\text{st}})$; then if $\kappa(A_{\text{st}})$ is still large, we apply the decomposition to $A_{\text{st}}^{-1}$ again to further decompose $\kappa(A_{\text{st}})$. However, the decomposition of $A^{-1}$ is based on the construction of $\mathcal{P}$ and $\Phi$, namely on the underlying energy decomposition $\mathcal{E}=\{E_k\}_{k=1}^m$ of $A$. Hence, we have to construct the corresponding {\bf energy decomposition} of $A_{\text{st}}$ before we implement the same operator decomposition on $A_{\text{st}}^{-1}$.

\subsection{Inherited energy decomposition} Let $\mathcal{E}=\{E_k\}_{k=1}^m$ be the energy decomposition of $A$, then the {\bf inherited energy decomposition} of $A_{\text{st}}=\Psi^TA\Psi$ with respect to $\mathcal{E}$ is simply given by $\mathcal{E}^\Psi=\{E_k^\Psi\}_{k=1}^m$ where 
\begin{equation}
E_k^\Psi=\Psi^TE_k\Psi,\quad k=1,2,\cdots,m.
\end{equation}
Notice that this inherited energy decomposition of $A_{\text{st}}$ with respect to $\mathcal{E}$ has the same number of energy elements as $\mathcal{E}$, which is not preferred and actually redundant in practice. Therefore we shall consider to reduce the energy decomposition of $A_{\text{st}}$. Indeed we will use $\widetilde{\Psi}$ instead of $\Psi$ in practice, where each $\tilde{\psi}_i$ is some local approximator of $\psi_i$ (obtained by \Cref{construction:psi}). Specifically, we will actually deal with $\widetilde{A}_{\text{st}}=\widetilde{\Psi}^TA\widetilde{\Psi}$ and thus we shall consider to find a proper condensed energy decomposition of $\widetilde{A}_{\text{st}}$.

If we see $\widetilde{A}_{\text{st}}$ as a matrix with respect to the reduced space $\mathbb{R}^N$, then for any vector $\bm{x} = \{x_1,\ldots, x_N\} \in \mathbb{R}^N$, the connection $\bm{x} \sim E^{\widetilde{\Psi}}$ between $\bm{x}$ and some $E^{\widetilde{\Psi}}=\widetilde{\Psi}^TE\widetilde{\Psi}$ comes from the connection between $E$ and those $\tilde{\psi}_i$ corresponding to nonzero $x_i$, and such connections are the key to constructing a partition $\mathcal{P}_{\text{st}}$ for $\widetilde{A}_{\text{st}}$. Recall that the support of each $\tilde{\psi}_i$($S_k(P_{j_i})$ for some $k$) is a union of patches, there is no need to distinguish among energy elements interior to the same patch when we deal with the connections between these elements and the basis $\widetilde{\Psi}$. Therefore we introduce the {\bf reduced inherited energy decomposition} of $\widetilde{A}_{\text{st}}=\widetilde{\Psi}^TA\widetilde{\Psi}$ as follows:

\begin{definition}[Reduced inherited energy decomposition]
With respect to the underlying energy decomposition $\mathcal{E}$ of $A$, the partition $\mathcal{P}$ and the corresponding $\widetilde{\Psi}$, the {\bf reduced inherited energy decomposition} of $\widetilde{A}_{\text{st}}=\widetilde{\Psi}^TA\widetilde{\Psi}$ is given by $\mathcal{E}^{\widetilde{\Psi}}_{re}=\{\underline{A}_{P_j}^{\widetilde{\Psi}}\}_{j=1}^M\cup\{E^{\widetilde{\Psi}}:E\in \mathcal{E}_{\mathcal{P}}^c\}$ with 
\begin{align}
\underline{A}_{P_j}^{\widetilde{\Psi}}&=\widetilde{\Psi}^T\underline{A}_{P_j}\widetilde{\Psi},\quad j=1,2,\cdots,M, \text{ and } \\
E^{\widetilde{\Psi}}&=\widetilde{\Psi}^TE\widetilde{\Psi},\quad \forall E\in \mathcal{E}_{\mathcal{P}}^c,
\label{eqt:reduced_energy}
\end{align}
where $\mathcal{E}_{\mathcal{P}}^c=\mathcal{E}\backslash \mathcal{E}_{\mathcal{P}}$ with $\mathcal{E}_{\mathcal{P}}=\{E\in \mathcal{E}:\exists P_j \in \mathcal{P}\ s.t.\ E\in P_j\}$.
\label{def:reduced_energy}
\end{definition}

Once we have the underlying energy decomposition of $A_{\text{st}}$ (or $\widetilde{A}_{\text{st}}$), we can repeat the procedure to decompose $A_{\text{st}}^{-1}$ (or $\widetilde{A}_{\text{st}}^{-1}$) in $\mathbb{R}^N$ as what we have done to $A^{-1}$ in $\mathbb{R}^n$. We will introduce the {\it multi-level decomposition} of $A^{-1}$ in the following subsection.

\subsection{Multiresolution operator decomposition}
\label{subsec:multi_operatordecomp}

Let $A^{(0)}=A$, and we construct $A^{(k)},B^{(k)}$ recursively from $A^{(0)}$. More precisely, let $\mathcal{E}^{(k-1)}$ be the underlying energy decomposition of $A^{(k-1)}$, and $\mathcal{P}^{(k)}$, $\Phi^{(k)}$, $\Psi^{(k)}$ and $U^{(k)}$ be constructed corresponding to $A^{(k-1)}$ and $\mathcal{E}^{(k-1)}$ in space $\mathbb{R}^{N^{(k-1)}}$, where $N^{(k-1)}$ is the dimension of $A^{(k-1)}$. We use one-level operator decomposition to decompose $(A^{(k-1)})^{-1}$ as  
\[
(A^{(k-1)})^{-1} = U^{(k)}\big((U^{(k)})^TA^{(k-1)}U^{(k)}\big)^{-1}(U^{(k)})^T + \Psi^{(k)}\big((\Psi^{(k)})^TA^{(k-1)}\Psi^{(k)}\big)^{-1}(\Psi^{(k)})^T,
\]
and then define $A^{(k)}=(\Psi^{(k)})^TA^{(k-1)}\Psi^{(k)}$, $B^{(k)}=(U^{(k)})^TA^{(k-1)}U^{(k)}$, and $\mathcal{E}^{(k)}=(\mathcal{E}^{(k-1)})^{\Psi^{(k)}}_{re}$ as in \Cref{def:reduced_energy}. Moreover, if we write 
\begin{subequations}
\begin{align}
\bm{\Phi}^{(1)}&=\Phi^{(1)},\quad \bm{\Phi}^{(k)}=\Phi^{(1)}\Phi^{(2)}\cdots\Phi^{(k-1)}\Phi^{(k)},\quad k\geq 1, \label{eqt:mutliresolution_phi&U&psi_a}\\
\mathcal{U}^{(1)}&=U^{(1)},\quad \mathcal{U}^{(k)}=\Psi^{(1)}\Psi^{(2)}\cdots\Psi^{(k-1)}U^{(k)},\quad k\geq1, \label{eqt:mutliresolution_phi&U&psi_b}\\
\bm{\Psi}^{(1)}&=\Psi^{(1)},\quad \bm{\Psi}^{(k)}=\Psi^{(1)}\Psi^{(2)}\cdots\Psi^{(k-1)}\Psi^{(k)},\quad k\geq 1, \label{eqt:mutliresolution_phi&U&psi_c}
\end{align}
\label{eqt:mutliresolution_phi&U&psi}
\end{subequations}
then one can prove by induction that for $k\geq1$,
\[A^{(k)}=(\bm{\Psi}^{(k)})^TA\bm{\Psi}^{(k)}=\big((\bm{\Phi}^{(k)})^TA^{-1}\bm{\Phi}^{(k)}\big)^{-1},\quad B^{(k)}=(\mathcal{U}^{(k)})^TA\mathcal{U}^{(k)},\]
\[(\bm{\Phi}^{(k)})^T\bm{\Phi}^{(k)}=(\bm{\Phi}^{(k)})^T\bm{\Psi}^{(k)}=I_{N^{(k)}},\quad \bm{\Psi}^{(k)}=A^{-1}\bm{\Phi}^{(k)}\big((\bm{\Phi}^{(k)})^TA^{-1}\bm{\Phi}^{(k)}\big)^{-1},\]
and for any integer $K$, 
\begin{equation}
A^{-1} = (A^{(0)})^{-1} = \sum_{k=1}^{K}\mathcal{U}^{(k)}\big((\mathcal{U}^{(k)})^TA\mathcal{U}^{(k)}\big)^{-1}(\mathcal{U}^{(k)})^T+\bm{\Psi}^{(K)}\big((\bm{\Psi}^{(K)})^TA\bm{\Psi}^{(K)}\big)^{-1}(\bm{\Psi}^{(K)})^T.
\label{eqt:Ainverse_overalldecomposition}
\end{equation}

\begin{remark}
$ $ \linebreak
\vspace{-3mm}
\begin{itemize}
\item[-] One shall notice that the partition $\mathcal{P}^{(k)}$ on each level $k$ is not a partition of the whole space $\mathbb{R}^n$, but a partition of the reduced space $\mathbb{R}^{N^{(k-1)}}$, and $\Phi^{(k)},\Psi^{(k)},U^{(k)}$ are all constructed corresponding to this $\mathcal{P}^{(k)}$ in the same reduced space. Intuitively, if the average patch sizes (basis number in a patch) for partition $\mathcal{P}^{(k)}$ is $s^{(k)}$, then we have $N^{(k)}= \frac{q^{(k)}}{s^{(k)}}N^{(k-1)}$, where $q^{(k)}$ is the integer for constructing $\Phi^{(k)}$.
\item[-] Generally, methods of multiresolution type use nested partitions/meshes that are generated only based on the computational domain \cite{bathe1976numerical,trottenberg2000multigrid}. But here the nested partitions are replaced by level-wisely constructed ones which are adaptive to $A^{(k)}$ on each level and require no a priori knowledge of the computational domain/space.
\item[-] In the Gamblet setting introduced in \cite{owhadi2017multigrid}, equations \cref{eqt:mutliresolution_phi&U&psi} together with \cref{eqt:Ainverse_overalldecomposition} can be viewed as the Gamblet Transform.
\end{itemize}
\end{remark}

The multiresolution operator decomposition up to a level $K$ is essentially equivalent to a decomposition of the whole space $\mathbb{R}^n$ \cite{owhadi2017multigrid} as 
\[\mathbb{R}^n=\mathcal{U}^{(1)}\oplus\mathcal{U}^{(2)}\oplus\cdots\oplus\mathcal{U}^{(K)}\oplus \bm{\Psi}^{(K)},\]
where again we also use $\mathcal{U}^{(k)}$( or $\bm{\Psi}^{(k)}$) to denote the subspace spanned by the basis $\mathcal{U}^{(k)}$( or $\bm{\Psi}^{(k)}$). Due to the $A$-orthogonality between these subspaces, using this decomposition to solve $A^{-1}$ is equivalent to solving $A^{-1}$ in each subspace separately (or more precisely solving $(B^{(k)})^{-1},k=1,\cdots,K$, or $(A^{(K)})^{-1}$), and by doing so we decompose the large condition number of $A$ into bounded pieces as the following corollary states.

\begin{corollary} If on each level $\Phi^{(k)}$ is given by \Cref{construction:phi} with integer $q^{(k)}$, then for $k\geq1$ we have
\[\lambda_{\max}(A^{(k)})\leq \delta(\mathcal{P}^{(k)},q^{(k)}),\qquad \lambda_{\min}(A^{(k)})\geq \lambda_{\min}(A),\]
\[\lambda_{\max}(B^{(k)})\leq \delta(\mathcal{P}^{(k-1)},q^{(k-1)})\lambda_{\max}\big((U^{(k)})^TU^{(k)}\big),\]
\[\lambda_{\min}(B^{(k)})\geq\frac{1}{\varepsilon(\mathcal{P}^{(k)},q^{(k)})^2}\lambda_{\min}\big((U^{(k)})^TU^{(k)}\big),\]
and thus 
\[\kappa(A^{(k)})\leq \delta(\mathcal{P}^{(k)},q^{(k)})\|A^{-1}\|_2,\]
\[\kappa(B^{(k)})\leq \varepsilon(\mathcal{P}^{(k)},q^{(k)})^2\delta(\mathcal{P}^{(k-1)},q^{(k-1)})\kappa\big((U^{(k)})^TU^{(k)}\big).\]
For consistency, we write $\delta(\mathcal{P}^{(0)},q^{(0)})=\lambda_{\max}(A^{(0)})=\lambda_{\max}(A)$.
\label{cor:Bconditonnumber}
\end{corollary}

\begin{proof}
These results follow directly from \Cref{thm:condition_number} and \Cref{lemma:Bst_eigenvalue}.
\end{proof}

\begin{remark}
The fact $\lambda_{\min}(B^{(k)})\gtrsim\frac{1}{\varepsilon(\mathcal{P}^{(k)},q^{(k)})^2}$ implies that the level $k$ is a level with resolution of scale no greater than $\varepsilon(\mathcal{P}^{(k)},q^{(k)})$, namely the space $\mathcal{U}^{(k)}$ is a subspace of the whole space $\mathbb{R}^n$ of scale finer than $\varepsilon(\mathcal{P}^{(k)},q^{(k)})$ with respect to $A$. This is essentially what multiresolution means in this decomposition.
\end{remark}

Now we have a multiresolution decomposition of $A^{-1}$, the applying of $A^{-1}$ (namely solving linear system $Ax=b)$) can break into the applying of $(B^{(k)})^{-1}$ on each level and the applying of $(A^{(K)})^{-1}$ on the bottom level. In what follows, we always assume $\kappa((U^{(k)})^TU^{(k)})=1$. Then the efficiency of the multiresolution decomposition in resolving the difficulty of large condition number of $A$ lies in the effort to bound each $\varepsilon(\mathcal{P}^{(k)},q^{(k)})^2\delta(\mathcal{P}^{(k-1)},q^{(k-1)})$ so that $B^{(k)}$ has a controlled spectrum width and can be efficiently solved using the CG type method. Define $\kappa(\mathcal{P}^{(k)},q^{(k)})=\varepsilon(\mathcal{P}^{(k)},q^{(k)})^2\delta(\mathcal{P}^{(k)},q^{(k)})$ and $\gamma^{(k)}=\frac{\varepsilon(\mathcal{P}^{(k)},q^{(k)})}{\varepsilon(\mathcal{P}^{(k-1)},q^{(k-1)})}$, then we can write
\begin{equation}\kappa(B^{(k)})\leq \varepsilon(\mathcal{P}^{(k)},q^{(k)})^2\delta(\mathcal{P}^{(k-1)},q^{(k-1)})=(\gamma^{(k)})^2\kappa(\mathcal{P}^{(k-1)},q^{(k-1)}).
\end{equation}
The partition condition number $\kappa(\mathcal{P}^{(k)},q^{(k))})$ is a level-wise information only concerning the partition $\mathcal{P}^{(k)}$. Similar to what we do in \Cref{alg:pair_clustering}, we will impose a uniform bound $c$ in the partitioning process so that $\kappa(\mathcal{P}^{(k)},q^{(k))})\leq c$ on every level. The ratio $\gamma^{(k)}$ reflects the scale gap between level $k-1$ and $k$, which is why it should measure the condition number(spectrum width) of $B^{(k)}$. However, it turns out that the choice of $\gamma^{(k)}$ is not arbitrary, and it will be subject to a restriction derived out of concern of sparsity.

So far the $A^{(k)}$ and $B^{(k+1)}$ are dense for $k\geq1$ since the basis $\Psi^{(k)}$ are global. It would be pointless to bound the condition number of $B^{(k)}$ if we cannot take advantage of the locality/sparsity of $A$. So in practice, the multiresolution operator decomposition is performed with localization on each level to ensure locality/sparsity. Thus we have the modified multiresolution operator decomposition in the following subsection.

\subsection{Multiresolution operator decomposition with localization} 
Let $\widetilde{A}^{(0)}=A$, and we construct $\widetilde{A}^{(k)},\widetilde{B}^{(k)}$ recursively from $\widetilde{A}^{(0)}$. More precisely, let $\widetilde{\mathcal{E}}^{(k-1)}$ be the underlying energy decomposition of $\widetilde{A}^{(k-1)}$, and $\mathcal{P}^{(k)}$, $\Phi^{(k)}$, $\Psi^{(k)}$ and $U^{(k)}$ be constructed corresponding to $\widetilde{A}^{(k-1)}$ and $\widetilde{\mathcal{E}}^{(k-1)}$ in space $\mathbb{R}^{N^{(k-1)}}$. We decompose $(\widetilde{A}^{(k-1)})^{-1}$ as 
\begin{equation}
(\widetilde{A}^{(k-1)})^{-1} = U^{(k)}\big((U^{(k)})^T\widetilde{A}^{(k-1)}U^{(k)}\big)^{-1}(U^{(k)})^T + \Psi^{(k)}\big((\Psi^{(k)})^T\widetilde{A}^{(k-1)}\Psi^{(k)}\big)^{-1}(\Psi^{(k)})^T.
\end{equation}
Let $\widetilde{\Psi}^{(k)}$ be a local approximator of $\Psi^{(k)}$. Then we define 
\begin{equation}
\widetilde{A}^{(k)}=(\widetilde{\Psi}^{(k)})^T\widetilde{A}^{(k-1)}\widetilde{\Psi}^{(k)},\qquad \widetilde{B}^{(k)}=(U^{(k)})^T\widetilde{A}^{(k-1)}U^{(k)},
\label{eqt:AB_construct_tilde}
\end{equation}
and $\widetilde{\mathcal{E}}^{(k)}=(\widetilde{\mathcal{E}}^{(k-1)})^{\widetilde{\Psi}^{(k)}}_{re}$ as in \Cref{def:reduced_energy}.\\

Similar to \Cref{cor:Bconditonnumber}, we have the following estimates on the condition numbers of $A^{(k)}$ and $B^{(k)}$:

\begin{corollary} If on each level $\Phi^{(k)}$ is given by \Cref{construction:phi} with integer $q^{(k)}$, and $\widetilde{\Psi}^{(k)}$ is a local approximator of $\Psi^{(k)}$ subject to localization error $\|\widetilde{\psi}_i^{(k)}-\psi_i^{(k)}\|_{A^{(k-1)}}\leq\frac{\epsilon}{\sqrt{N^{(k)}}}$, then for $k\geq1$ we have
\[\lambda_{\max}(\widetilde{A}^{(k)})\leq \Big(1+\frac{\epsilon}{\sqrt{\delta(\mathcal{P}^{(k)},q^{(k)})}}\Big)^2\delta(\mathcal{P}^{(k)},q^{(k)}),\qquad \lambda_{\min}(\widetilde{A}^{(k)})\geq \lambda_{\min}(A),\]
\[\lambda_{\max}(\widetilde{B}^{(k)})\leq \Big(1+\frac{\epsilon}{\sqrt{\delta(\mathcal{P}^{k-1)},q^{(k-1)})}}\Big)^2\delta(\mathcal{P}^{(k-1)},q^{(k-1)})\lambda_{\max}\big((U^{(k)})^TU^{(k)}\big),\]
\[\lambda_{\min}(\widetilde{B}^{(k)})\geq\frac{1}{\varepsilon(\mathcal{P}^{(k)},q^{(k)})^2}\lambda_{\min}\big((U^{(k)})^TU^{(k)}\big),\]
and thus 
\[\kappa(\widetilde{A}^{(k)})\leq \Big(1+\frac{\epsilon}{\sqrt{\delta(\mathcal{P}^{(k)},q^{(k)})}}\Big)^2\delta(\mathcal{P}^{(k)},q^{(k)})\|A^{-1}\|_2,\]
\[\kappa(\widetilde{B}^{(k)})\leq \Big(1+\frac{\epsilon}{\sqrt{\delta(\mathcal{P}^{(k-1)},q^{(k-1)})}}\Big)^2\varepsilon(\mathcal{P}^{(k)},q^{(k)})^2\delta(\mathcal{P}^{(k-1)},q^{(k-1)})\kappa\big((U^{(k)})^TU^{(k)}\big).\]
For consistency, we write $\delta(\mathcal{P}^{(0)},q^{(0)})=\lambda_{\max}(\widetilde{A}^{(0)})=\lambda_{\max}(A)$.
\label{cor:Bconditonnumber2}
\end{corollary}

\begin{proof}
These results follow directly from the proof of \Cref{thm:condition_number}, \Cref{cor:localpsi_conditionnum} and \Cref{lemma:Bst_eigenvalue}.
\end{proof}

One can indeed prove that $\delta(\mathcal{P}^{(k)},q^{(k)})\geq \frac{1}{\varepsilon(\mathcal{P}^{(k)},q^{(k)})^2}$, and thus $\frac{\epsilon}{\sqrt{\delta(\mathcal{P}^{(k-1)},q^{(k-1)})}}\leq \epsilon \varepsilon(\mathcal{P}^{(k)},q^{(k)})$ is a small number. Therefore \Cref{cor:Bconditonnumber2} states that the multiresolution decomposition with localization has estimates on condition numbers of the same order as in \Cref{cor:Bconditonnumber}, i.e. $\kappa(\widetilde{A}^{(k)})\leq O(\delta(\mathcal{P}^{(k)},q^{(k)})\|A^{-1}\|_2)$ and $\kappa(\widetilde{B}^{(k)})\leq O(\varepsilon(\mathcal{P}^{(k)},q^{(k)})^2\delta(\mathcal{P}^{(k-1)},q^{(k-1)}))$. Having this in hand, we proceed to discuss the desired sparsity of $\widetilde{A}^{(k)}$ and $\widetilde{B}^{(k)}$.\\

\textbf{Locality Preservation:} Similar to the locality discussion of $\widetilde{A}_{st}$ in \Cref{Sec:partition}, under the locality condition \cref{condition:localityA}, we have the following recursive estimate on the number of nonzero entries of each $A^{(k)}$ as
\begin{equation}
nnz(\widetilde{A}^{(k)})= O(nnz(\widetilde{A}^{(k-1)})\cdot (q^{(k)})^2\cdot \frac{1}{s^{(k)}}\cdot (r^{(k)})^d),
\label{eqt:Aktilde_sparsity}
\end{equation}
where $s^{(k)}$ is the average patch size of $\mathcal{P}^{(k)}$, and $r^{(k)}$ is the decay radius of $\widetilde{\Psi}^{(k)}$. Also, noticing that $\widetilde{B}^{(k)}=(U^{(k)})^T\widetilde{A}^{(k-1)}U^{(k)}$ and that the basis $U^{(k)}$ are local vectors of support size $s^{(k)}$, we have 
\begin{equation}nnz(\widetilde{B}^{(k)})= O( nnz(\widetilde{A}^{(k-1)})\cdot s^{(k)}).
\end{equation}
In fact, the basis $U^{(k)}$ can be computed from $\Phi^{(k)}$ using the implicit QR factorization\cite{saad2011numerical}, and thus the matrix multiplication with respect to $U^{(k)}$ can be done by using the Householder vectors in time linear to $q^{(k)}\cdot N^{(k)}$. Therefore, when we evaluate $\widetilde{B}^{(k)}=(U^{(k)})^T\widetilde{A}^{(k-1)}U^{(k)}$ (in iterative method), only the NNZ of $\widetilde{A}^{(k-1)}$ matters. In brief, we need to preserve the locality of $A^{(k)}$ down through all the levels to ensure the efficiency of the multiresolution decomposition with localization. But the accumulation of the factor $\frac{(q^{(k)})^2 (r^{(k)})^d}{s^{(k)}}$, if not well controlled, will compromise the sparsity inherited from $\widetilde{A}^{(0)}=A$. Therefore a necessary condition for the decomposition to keep sparsity is 
\[o(s^{(k)})\geq (q^{(k)})^2 (r^{(k)})^d,\quad k\geq 1,\]
under which we have the sparsity estimate $nnz(\widetilde{A}^{(k)})=O(nnz(A))$. In particular, when we impose the localization error $\|\tilde{\psi}_i^{(k)}-\psi_{i}^{(k)}\|_{A^{(k-1)}}\leq\frac{\epsilon}{\sqrt{N^{(k)}}}$ on each level $k$ for some uniform $\epsilon$, we have $r^{(k)}=O(\log\frac{1}{\epsilon}+\log N^{(k)} +\log\delta(\mathcal{P}^{(k)},q^{(k)}))$ according the discussions in \Cref{sec:choose_psi}. Then the sparsity condition becomes 
\begin{equation}
o(s^{(k)})\geq (q^{(k)})^2 \big(\log\frac{1}{\epsilon}+\log N^{(k)} +\log\delta(\mathcal{P}^{(k)},q^{(k)})\big)^d,\quad k\geq 1.
\label{con:sparsitypreservation}
\end{equation}
This lower bound of the patch size $s^{(k)}$ means that we need to compress enough dimensions from higher level to lower level in order to preserve sparsity due to the outreaching support of the localized basis $\widetilde{\Psi}^{(k)}$.

In practice, we will choose some $\epsilon$ smaller than the top level scale $\varepsilon(\mathcal{P}^{(1)},q^{(1)})$ and a uniform integer $q$. By imposing uniform condition bound $\kappa(\mathcal{P}^{(k)},q^{(k)})\leq c$ we have $\delta(\mathcal{P}^{(k)},q^{(k)})\leq\frac{c}{\varepsilon(\mathcal{P}^{(k)},q^{(k)})^2}\leq \frac{c}{\epsilon^2}$. Therefore a safe uniform criterion for patch size $s^{(k)}$ is 
\begin{equation}
O(s^{(k)})=s=q^2(\log\frac{1}{\epsilon}+\log n)^{d+l},
\label{eqt:patchsize}
\end{equation}
for some small $l>0$, which asymptotically, when $n$ goes large and $\epsilon$ goes small, will ensure $nnz(A^{(k)})=O(nnz(A))$ down through the decomposition. Since the decomposition should stop when $N^{(K)}$, the dimension of $A^{(K)}$ is small enough, namely when $n=O((s/q)^K)$, \Cref{eqt:patchsize} also gives us an estimate of the total level number as
\begin{equation}
K=O(\log_{s/q}n)=O\big(\frac{\log n}{\log(q(\log\frac{1}{\epsilon}+\log n)^{d+l})}\big)=O\big(\frac{\log n}{\log(\log\frac{1}{\epsilon}+\log n)}\big).
\label{eqt:levelnumber}
\end{equation}

\subsection*{Choice of scale ratio $\gamma$} Recall that the partition $\mathcal{P}^{(k)}$ is a partition of basis in the space $\mathbb{R}^{N^{(k-1)}}$. By tracing back to the top level, we can also see it as a partition in the original space $\mathbb{R}^n$. Denoting $R^{(k)}$ to be the average radius (with respect to adjacency defined by $A$) and $S^{(k)}$ to be the average patch size of the patches (with respect to $\mathbb{R}^n$) in $\mathcal{P}^{(k)}$, we have $S^{(k)}=O((R^{(k)})^d)$ under the locality condition \cref{condition:localityA}, and an intuitive geometry estimate gives $\frac{S^{(k)}}{S^{(k-1)}}=s^{(k)}$. As a consequence, under the local energy decomposition condition \cref{condition:scalesproperty} of order $(q,p)$, we have the following estimate
\begin{equation}
\gamma^{(k)}=\frac{\varepsilon(\mathcal{P}^{(k)},q)}{\varepsilon(\mathcal{P}^{(k-1)},q)}=O(\Big(\frac{R^{(k)}}{R^{(k-1)}}\Big)^p)=O(\Big(\frac{S^{(k)}}{S^{(k-1)}}\Big)^\frac{p}{d})=O((s^{(k)})^\frac{p}{d}).
\label{eqt:scaleestimate}
\end{equation}
Such estimate arises naturally in a lot of PDE problems, especially when the smallest eigenvalues of local operators have clear dependence on the domain size, the dimension of the space and the order of the equation \cite{hou2016sparse}. Under the sparsity condition \cref{con:sparsitypreservation} and considering $q$ as a constant, we require
\begin{equation}
o(\gamma^{(k)})\geq (\log\frac{1}{\epsilon}+\log N^{(k)} +\log\delta(\mathcal{P}^{(k)},q^{(k)}))^p,
\end{equation}
to ensure the sparsity of the decomposition, and similarly a safe, uniform choice of the scale ratio $\gamma^{(k)}$ is 
\begin{equation}
\gamma^{(k)}=\gamma=(\log\frac{1}{\epsilon}+\log n)^{p+l},
\label{eqt:scaleratio}
\end{equation}
for some small $l>0$. Such choice provides a uniform bound on the condition number of $\widetilde{B}^{(k)}$ as
\begin{equation}
\kappa(\widetilde{B}^{(k)})\leq O(\varepsilon(\mathcal{P}^{(k)},q^{(k)})^2\delta(\mathcal{P}^{(k-1)},q^{(k-1)})) \leq O((\log\frac{1}{\epsilon}+\log n)^{p+l})
\label{eqt:tildeBconditionnumber}
\end{equation}
when a uniform condition bound $\kappa(\mathcal{P}^{(k)},q^{(k)})\leq c$ is imposed by algorithm. Notice that the ratio $\gamma^{(k)}$ is only defined for $k\geq2$, thus the estimate \cref{eqt:tildeBconditionnumber} is valid for $k\geq2$. For consistency, we choose $\varepsilon(\mathcal{P}^{(1)},q^{(1)})^2=O(\frac{(\log\frac{1}{\epsilon}+\log n)^{p+l}}{\|A\|_2})$ so that \cref{eqt:tildeBconditionnumber} is also valid for $k=1$.

\begin{remark} By estimate \cref{eqt:tildeBconditionnumber}, the bound on $\kappa(\widetilde{B}^{(k)})$ will go to infinity when $n$ goes to infinity. In our construction of the multiresolution decomposition for resolving large condition number of $A$, we can not asymptotically have an absolute constant bound for $\kappa(\widetilde{B}^{(k)})$ on all levels, due to the required preservation of sparsity. This difficulty comes from the inductive nature of the algorithm that the posterior estimate of the sparsity of $\widetilde{A}^{(k)}$ is based on the sparsity of $\widetilde{A}^{(k-1)}$, as shown in \cref{eqt:Aktilde_sparsity}. However, in \cite{owhadi2017universal}, the existence of nested measurement function $\Phi$ is assumed a-priori before the construction of the multiresolution structure, and thus the sparsity of $\widetilde{A}^{(k)}$ can be inherited directly from $\widetilde{A}^{(0)}$, which avoids the accumulation of the factor $\frac{(q^{(k)})^2 (r^{(k)})^d}{s^{(k)}}$ through levels. As a result, the sparsity of $\widetilde{A}^{(k)}$ does not contradict the uniform bound of $\kappa(\widetilde{B}^{(k)})$.\\
\end{remark}

\subsection*{Error estimate} Using multiresolution operator decomposition with localization to solve $A^{-1}$, the error on each level comes from two main sources: (i) the localization error between $\Psi^{(k)}\big((\Psi^{(k)})^T\widetilde{A}^{(k-1)}\Psi^{(k)}\big)^{-1}(\Psi^{(k)})^T$ and $\widetilde{\Psi}^{(k)}\big((\widetilde{\Psi}^{(k)})^T\widetilde{A}^{(k-1)}\widetilde{\Psi}^{(k)}\big)^{-1}(\widetilde{\Psi}^{(k)})^T$; (ii) the error caused by solving $(A^{(k)})^{-1}=\big((\widetilde{\Psi}^{(k)})^T\widetilde{A}^{(k-1)}\widetilde{\Psi}^{(k)}\big)^{-1}$ (or $(\widetilde{B}^{(k)})^{-1}$) with iterative type methods. To come up with an estimate of the total error, we perform a standard analysis of error accumulation in an inductive manner.

\begin{thm} 
Given an integer $K$, let $\mathrm{Inv}(A)$ denote the solver for $A^{-1}$ using $K$-level's multiresolution operator decomposition with localization. Assume that 
\begin{enumerate}[label=(\roman*)]
\item each $(\widetilde{B}^{(k)})^{-1}$ can be solved efficiently subject to a uniform relative error bound $\mathrm{err}_B$ in the sense that the solver $\mathrm{Inv}(\widetilde{B}^{(k)})$ (as a linear operator) satisfies 
\begin{equation}
\|(\widetilde{B}^{(k)})^{-1}b-\mathrm{Inv}(\widetilde{B}^{(k)})b\|_{\widetilde{B}^{(k)}}\leq \mathrm{err}_B \| b \|_{(\widetilde{B}^{(k)})^{-1}},\quad \forall b\in \mathbb{R}^{N^{(k)}},\quad \forall 1\leq k\leq K;
\end{equation}
\label{itm:assumption1}
\item at level $K$, $(\widetilde{A}^{(K)})^{-1}$ can be solved efficiently subject to a relative error $\mathrm{err}^{(K)}_A$ in the sense that the solver $\mathrm{Inv}(\widetilde{A}^{(K)})$ satisfies 
\begin{equation}
\|(\widetilde{A}^{(K)})^{-1}b-\mathrm{Inv}(\widetilde{A}^{(K)})b\|_{\widetilde{A}^{(K)}}\leq \mathrm{err}^{(K)}_A\|b\|_{(\widetilde{A}^{(K)})^{-1}},\quad \forall b\in \mathbb{R}^{N^{(K)}}.
\end{equation}
\label{itm:assumption2}
\item each $\widetilde{\Psi}^{(k)}$ satisfies the localization approximation property 
\begin{equation}\|\widetilde{\psi}^{(k)}_i-\psi^{(k)}_i\|_{\widetilde{A}^{(k-1)}}\leq\frac{\mathrm{err}_{loc}}{2\sqrt{N^{(k)}\|A^{-1}\|_2}},\quad 1\leq i\leq N^{(k)},
\end{equation}
with a uniform constant $\mathrm{err}_{loc}$.
\label{itm:assumption3}
\end{enumerate}
Then we have
\[\|A^{-1}b-\mathrm{Inv}(A)b\|_A\leq \mathrm{err}_{total}\|b\|_{A^{-1}},\quad \forall b\in \mathbb{R}^n, \]
and in $\|\cdot\|_2$,
\[\|A^{-1}b-\mathrm{Inv}(A)b\|_2\leq \mathrm{err}_{total}\|A^{-1}\|_2\|b\|_2,\quad \forall b\in \mathbb{R}^n.\]
where 
\[\mathrm{err}_{total}=K(\mathrm{err}_B+\mathrm{err}_{loc})+\mathrm{err}^{(K)}_A.\]
\label{theorem:MMD_levelerrorbound}
\end{thm}
\begin{remark}
Assumption \ref{itm:assumption1} is reasonable since each $\widetilde{B}^{(k)}$ inherit the sparsity from $A$, and its condition numbers can be well bounded in order $O((\log\frac{1}{\mathrm{err}_{loc}}+\log n)^{p+l})$. Assumption \ref{itm:assumption2} is reasonable since each $\widetilde{A}^{(K)}$ is of small dimension when $K$ is as large as in \Cref{eqt:levelnumber}. Assumption \ref{itm:assumption3} is reasonable due to the exponential decay property of each $\Psi^{(k)}$. Indeed, to ensure locality of reduced energy decomposition, the localization error control can be relaxed in practice. Such relaxed error can be fixed by doing compensation computation as we will see in \Cref{numericalexample_mmd}.
\end{remark}

\begin{figure}[!h]
\centering
\includegraphics[width=1.0\textwidth]{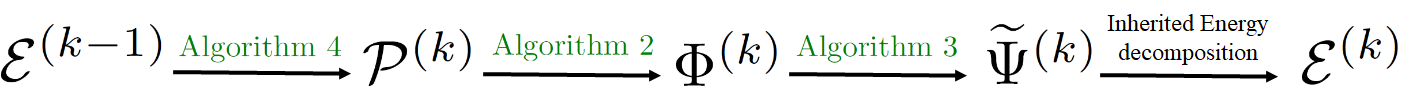}
\caption{Process flowchart of \Cref{alg:MMD_localized}.}
\label{fig:flowchart2}
\end{figure}

\subsection{Algorithm} We now summarize the procedure of {\bf Multiresolution Matrix Decomposition (MMD) with Localization} as \Cref{alg:MMD_localized}, and the use of MMD to solve linear system as \Cref{alg:MMD_solver}. Also, \Cref{fig:flowchart2} shows the flowchart of \Cref{alg:MMD_localized}.

\begin{algorithm}[!h]
\caption{\it Multiresolution Matrix Decomposition (MMD) with Localization}
\label{alg:MMD_localized}
\begin{algorithmic}[1]
\REQUIRE{SPD matrix $A=\widetilde{A}^{(0)}$, energy decomposition $\mathcal{E}=\widetilde{\mathcal{E}}^{(0)}$, underlying basis $\mathcal{V}$, localization constant $\epsilon$, level number $K$, $q^{(k)}$, error factor bound $\varepsilon^{(k)}$ and condition bound $c^{(k)}$ for each level.}
\ENSURE{$\widetilde{A}^{(K)}$, $\widetilde{\Psi}^{(k)}, U^{(k)}$, and $\widetilde{B}^{(k)}$.}
\FOR{$k=1:K$}
\STATE{Construct $\mathcal{P}^{(k)},\Phi^{(k)},U^{(k)},\widetilde{\Psi}^{(k)}$ with \Cref{alg:operator_compression}, with respect to $\widetilde{A}^{(k-1)},\widetilde{\mathcal{E}}^{(k-1)}$, and subject to $q^{(k)},\varepsilon^{(k)},c^{(k)}$ and localization error $\frac{\epsilon}{\sqrt{N^{(k)}}}$;} \label{line2:MMD_localized}
\STATE{Compute $\widetilde{A}^{(k)}$ and $\widetilde{B}^{(k)}$ by \Cref{eqt:AB_construct_tilde};} \label{line3:MMD_localized}
\STATE{Compute reduced energy $\widetilde{\mathcal{E}}^{(k)}$ by \Cref{eqt:reduced_energy};} \label{line4:MMD_localized}
\STATE{output/store $\widetilde{\Psi}^{(k)},U^{(k)},\widetilde{B}^{(k)}$;}
\ENDFOR
\STATE{output/store $\widetilde{A}^{(K)}$.}
\end{algorithmic}
\end{algorithm}

\begin{algorithm}[!h]
\caption{\it Solving linear system with MMD with Localization}
\label{alg:MMD_solver}
\begin{algorithmic}[1]
\REQUIRE{$\widetilde{\Psi}^{(k)},U^{(k)},\widetilde{B}^{(k)}$ for $k=1,2,\cdots,K$, $\widetilde{A}^{(K)}$, load vector $b=b^{(0)}$, prescribed relative accuracy $\epsilon$}
\ENSURE{Approximated solution $x^{(0)}$.}
\FOR{$k=1:K$} \label{line:MMD_solver_line1}
\STATE{$z^{(k)}=(U^{(k)})^T b^{(k-1)}$;}
\STATE{Solve $\widetilde{B}^{(k)}y^{(k)}=z^{(k)}$ up to relative error $\epsilon$;} \label{line:MMD_solver_line3}
\STATE{$b^{(k)}=(\widetilde{\Psi}^{(k)})^Tb^{(k-1)}$;}
\ENDFOR
\STATE{Solve $\widetilde{A}^{(K)}x^{(K)}=b^{(K)}$ up to relative error $\epsilon$;}
\FOR{$k=K:1$}
\STATE{$x^{(k-1)}=U^{(k)}y^{(k)}+\widetilde{\Psi}^{(k)}x^{(k)}$;}
\ENDFOR
\end{algorithmic}
\end{algorithm}

\begin{remark}
$ $ \linebreak
\vspace{-3mm}
\begin{itemize}
\item[-] Once the MMD is obtained, the first for-loop (Line~\ref{line:MMD_solver_line1}) in \Cref{alg:MMD_solver} can be performed in parallel, which makes it much more efficient than non-parallelizable iterative methods.
\item[-] Once the whole decomposition structure is completed, we can a posterior omit the level-wise energy decompositions and partitions. Then if we see our level-wisely constructed $\Phi$ as a nested sequence \cref{eqt:mutliresolution_phi&U&psi_a}, our decomposition is structurally equivalent to the result obtained in \cite{owhadi2017universal}, where the existence of such nested $\Phi$ is a priori assumed. Therefore, the required properties of the nested sequence in Condition 2.3 of \cite{owhadi2017universal} are similar to the assumption in \Cref{theorem:MMD_levelerrorbound}.
\end{itemize}
\end{remark}

\subsection*{Complexity of \Cref{alg:MMD_localized}} Assume that locality conditions \cref{condition:localityA},\cref{condition:scalesproperty},\cref{condition:localregularity} are trun with constant $d,p,q,c$. Then all $q^{(k)}$ and $c^{(k)}$ are chosen uniformly over levels to be $q, c$ respectively. $\varepsilon^{(k)}$ is chosen subject to scale ratio choice \cref{eqt:scaleratio}, $(\varepsilon^{(1)})^2=\frac{(\log\frac{1}{\epsilon}+\log n)^{p+l}}{\|A\|_2}$ for some small $l>0$, and $\epsilon$ is chosen so that $\epsilon\leq \varepsilon^{(1)}$. Due to the condition bound $c$, we have condition number estimate \cref{eqt:tildeBconditionnumber}. Then the complexity of Line~\ref{line2:MMD_localized} can be modified from \cref{eqt:operatorcompression_complexity} as 
\[O\big(d\cdot s^2\cdot\log s\cdot n\big)+O(\log s\cdot n\cdot\log n)+O\big(q\cdot n \cdot (\log\frac{1}{\epsilon}+\log n)^{p+l}\cdot (\log\frac{1}{\epsilon}+\log n)^{d+1}\big),\]
where $s=O((\log\frac{1}{\epsilon}+\log n)^{d(1+l/p)})$ according to estimate \cref{eqt:scaleestimate}. The complexity of Line~\ref{line3:MMD_localized} and \ref{line4:MMD_localized}(sparse matrices multiplication) together can be bounded by 
\[O(n\cdot (\log\frac{1}{\epsilon}+\log n)^{3d})\]
due to the locality of $\widetilde{\Psi}^{(k)}$ and the inherited locality of $\widetilde{A}^{(k-1)}$ and $\widetilde{B}^{(k-1)}$. Therefore the complexity on each level can be bounded by
\begin{equation}
O\big(d\cdot s^2\cdot\log s\cdot n\big)+O(\log s\cdot n\cdot\log n)+O\big(q\cdot n\cdot (\log\frac{1}{\epsilon}+\log n)^{3d+p}\big),
\end{equation}
where we have assumed that $d\geq 1\geq l$.
Then the total complexity of \Cref{alg:MMD_localized} is then the level number $K$ times \cref{eqt:MMD_solver_layercomplexity}. By \Cref{eqt:levelnumber}, we have $K=O\big(\frac{\log n}{\log(\log\frac{1}{\epsilon}+\log n)}\big)\leq O(\log n)$ and $K\log s=O(\log n)$. Thus the total complexity of \Cref{alg:MMD_localized} is
\begin{align}
& O\big(d\cdot s^2\cdot\log n\cdot n\big)+O(n\cdot(\log n)^2)+O\big(K\cdot q\cdot n\cdot (\log\frac{1}{\epsilon}+\log n)^{3d+p}\big)\nonumber \\
\leq &\ O(m\cdot \log n\cdot (\log\frac{1}{\epsilon}+\log n)^{3d+p} ).
\label{eqt:MMD_solver_layercomplexity}
\end{align}
where $m=O(d\cdot n)$ is the number of nonzero entries of $A$.

\subsection*{Complexity of \Cref{alg:MMD_solver}} Assume that the relative accuracy $\epsilon$ is the same as the $\epsilon$ in \Cref{alg:MMD_localized}. Recall that the number of nonzero entries of each $\widetilde{A}^{(k)}$ is bounded by by $O(nnz(A))=O(m)$, and the condition number each $\widetilde{B}^{(k)}$ can be bounded by $O((\log\frac{1}{\epsilon}+\log n)^{p+l})$, then the complexity of solving linear system in Line~\ref{line:MMD_solver_line3} using a CG type method is bounded by 
\[O(m\cdot (\log\frac{1}{\epsilon}+\log n)^{p+l}\cdot\log\frac{1}{\epsilon}).\]
Therefore if we use a CG type method to solve all inverse problems involved in \Cref{alg:MMD_solver}, based on the MMD with localization given by \Cref{alg:MMD_localized}, the running time of \Cref{alg:MMD_solver} subject to level-wise relative accuracy $\epsilon$ is 
\[O(K\cdot m \cdot (\log\frac{1}{\epsilon}+\log n)^{p+l}\cdot \log\frac{1}{\epsilon})\leq O(m \cdot (\log\frac{1}{\epsilon}+\log n)^{p+l}\cdot \log\frac{1}{\epsilon}\cdot \log n).\]
However, by \Cref{theorem:MMD_levelerrorbound}, the total accuracy is $\epsilon_{total}=O(K\epsilon)$. Thus the complexity of \Cref{alg:MMD_solver} subject to a total relative accuracy $\epsilon_{total}$ is 
\begin{equation}
O(m \cdot (\log\frac{1}{\epsilon_{total}}+\log n)^{p+l}\cdot (\log\frac{1}{\epsilon_{total}}+\log\log n)\cdot \log n).
\end{equation}

\subsection{Multilevel operator compression:} Also we can consider the multiresolution matrix decomposition from the perspective of operator compression. For any $K$, by omitting the finer scale subspaces $\mathcal{U}^{(k)}$, $k=1,2,\cdots,K$, we get an effective approximator of $A^{-1}$ as
\begin{equation}
A^{-1}\approx\bm{\Psi}^{(K)}\big((\bm{\Psi}^{(K)})^TA\bm{\Psi}^{(K)}\big)^{-1}(\bm{\Psi}^{(K)})^T=P_{\bm{\Psi}^{(K)}}^AA^{-1}.
\label{eqt:invA_compression}
\end{equation}
Intuitively, this approximation lies above the scale of $\varepsilon(\mathcal{P}^{(K)},q^{(K)})$, and therefore should have a corresponding dominant compression error. However we should again notice that the composite basis $\bm{\Phi}^{(k)}$ is not given a priori and directly in $\mathbb{R}^n$, but constructed level by level using the information of $A^{(k)}$ on each level, and recall that the {\bf error factor} $\varepsilon(\mathcal{P}^{(k)},q^{(k)})$ is computed with respect to the reduced space $\mathbb{R}^{N^{(k)}}$, not to the whole space $\mathbb{R}^n$. Thus the total error of compression \cref{eqt:invA_compression} is accumulated over all levels finer than level $K$. To quantify such compression error, we introduce the following theorem:
\begin{thm}
\label{thm:error_accumulate}
Assume that on each level $\Phi^{(k)}$ is given by \Cref{construction:phi} with integer $q^{(k)}$. Then we have 
\begin{equation}
\|x-P_{\bm{\Phi}^{(K)}}x\|_2\leq   \Big(\sum_{k=0}^{K}\varepsilon(\mathcal{P}^{(k)},q^{(k)})^2\Big)^\frac{1}{2}\|x\|_A,\quad \forall x\in \mathbb{R}^n,
\label{eqt:error_accumulate}
\end{equation}
and thus for any $x\in \mathbb{R}^n$ and $b=Ax$, we have
\[\|x-P_{\bm{\Psi}^{(K)}}^Ax\|_A\leq \Big(\sum_{k=1}^{K}\varepsilon(\mathcal{P}^{(k)},q^{(k)})^2\Big)^\frac{1}{2}\|b\|_2,  \]
\[\|x-P_{\bm{\Psi}^{(K)}}^Ax\|_2\leq \Big(\sum_{k=1}^{K}\varepsilon(\mathcal{P}^{(k)},q^{(k)})^2\Big)\|b\|_2,\]
\[\|A^{-1}-P_{\bm{\Psi}^{(K)}}^AA^{-1}\|_2\leq \Big(\sum_{k=1}^{K}\varepsilon(\mathcal{P}^{(k)},q^{(k)})^2\Big).\]
\end{thm}
\begin{remark}
$ $ \linebreak
\vspace{-3mm}
\begin{itemize}
\item[-] Though the compression error is in a cumulative form, if we assume that $\varepsilon(\mathcal{P}^{(k)},q^{(k)})$ increases with $k$ at a certain ratio $\frac{\varepsilon(\mathcal{P}^{(k)},q^{(k)})}{\varepsilon(\mathcal{P}^{(k-1)},q^{(k-1)})}\geq\gamma$ for some $\gamma > 1$, then it's easy to see that 
\[\sum_{k=1}^{K}\varepsilon(\mathcal{P}^{(k)},q^{(k)})^2\leq \frac{\gamma^2}{\gamma^2-1}\varepsilon(\mathcal{P}^{(K)},q^{(K)})^2,\]
which is an error of scale $\varepsilon(\mathcal{P}^{(K)},q^{(K)})^2$ as we expected.
\item[-] Again one shall be aware of the difference between the one-level compression with {\bf error factor} $\varepsilon(\mathcal{P}^{(K)},q^{(K)})$ and the multi-level compression in \Cref{subsec:multi_operatordecomp}. A one-level compression with {\bf error factor} $\varepsilon(\mathcal{P}^{(K)},q^{(K)})$ requires to construct $\mathcal{P}^{(K)},\bm{\Phi}^{(k)}$ and so on directly with respect to $A$ in the whole space $\mathbb{R}^n$, which involves solving eigenvalue problems on considerably large patches in $\mathcal{P}^{(K)}$ when $\varepsilon(\mathcal{P}^{(K)},q^{(K)})$ is a coarse scale. But the multi-level compression in \Cref{subsec:multi_operatordecomp} is computed hierarchically with bounded compression ratio between levels, and thus only involves eigenvalue problems on patches of well bounded size($s=O(\log\frac{1}{\epsilon}+\log n)^{d+l}$) in each reduced space $\mathbb{R}^{N^{(k)}}$, and is thus more tractable in practice.
\item[-] One can also analyze the compression error when localization of each $\Psi^{(k)}$ is considered. The analysis would be similar to the one in \Cref{thm:error_accumulate}. 
\end{itemize}
\end{remark}

\subsection{Numerical Example for Multiresultion Matrix Decomposition (MMD)}
\label{numericalexample_mmd}
\begin{figure}[!h]
\centering
    \begin{subfigure}[b]{0.35\textwidth}
        \centering
        \includegraphics[width=1.0\textwidth]{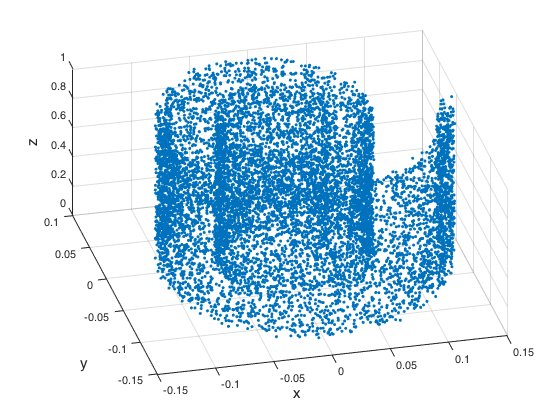}
        \caption{Front view}
        \label{fig:MMD_ptcloud1}
    \end{subfigure}
    \begin{subfigure}[b]{0.35\textwidth}
        \centering
        \includegraphics[width=1.0\textwidth]{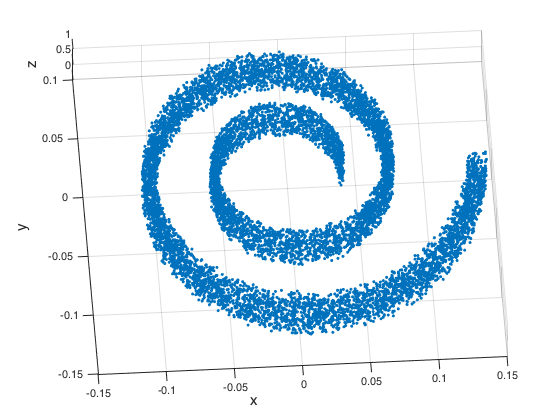}
        \caption{Top view}
        \label{fig:MMD_ptcloud2}
    \end{subfigure}
    \caption{``Roll surface'' constructed by \Cref{eqt:form_rollsurface}.}
\end{figure}
Our third numerical example shows the effectiveness of using multiresolution matrix decomposition (MMD) with localization to solve a graph Laplacian system. Again we use the same setup in the Example 1 in \Cref{subsec:numerical_example1}, with density factor $\eta = 2$. But this time the vertices of the graph are randomly distributed around a 2 dimensional roll surface of area 1 in $\mathbb{R}^3$. The distribution is a combination of a uniform distribution over the surface and an up to $10\%$ random displacement off the surface. More precisely, the two-dimensional roll is characterized as
\[(x(t),y(t),z)=(\rho(t)\cos(\theta(t)),\rho(t)\sin(\theta(t)),z),\quad t\in[0,1],z\in[0,1],\]
where 
\[\theta(t)=\frac{1}{a}\log\Big(1+t\big(e^{4\pi a}-1\big)\Big),\quad \rho(t)=\frac{a}{\sqrt{1+a^2}}\big(t+\frac{1}{e^{4\pi a}-1}\big), \]
and so $\sqrt{(\rho'(t))^2+(\rho(t)\theta'(t))^2}=1$. Each vertex $(x_i,y_i,z_i)$ is generated by
\begin{equation}
(x_i,y_i,z_i)=(\eta_i\rho(t_i)\cos(\theta(t_i)),\eta_i\rho(t_i)\sin(\theta(t_i)),z_i),\quad i=1,2,\cdots,n,
\label{eqt:form_rollsurface}
\end{equation}
where $t_i\overset{i.i.d}{\sim} \mathcal{U}[0,1]$, $z_i\overset{i.i.d}{\sim} \mathcal{U}[0,1]$, and $\eta_i\overset{i.i.d}{\sim} \mathcal{U}[0.9,1.1]$. In this example we take $n=10000$ and $a=0.1$. \Cref{fig:MMD_ptcloud1} and \cref{fig:MMD_ptcloud2} show the point cloud of all vertices. This explicit expression, however, is considered as a hidden geometric information, and is not employed in our partitioning algorithm.
\begin{table}[!h]
\centering
\footnotesize
\renewcommand{\arraystretch}{1.5}
    \begin{tabular}{|c|c|c|c|c|}
    \hline
    Level & Size & $\#$Nonzeros & Condition Number &Complexity\\ \hline
    $L$ & $10000\times10000$ & $128018\triangleq m$ & $1.93\times 10^7$ & $2.47\times 10^{12}$ \\ \hline
    $\widetilde{B}^{(1)}$ & $1898\times1898$ & $9812\approx0.07m$ & $1.72\times 10^2$ & $1.69\times 10^{6}$\\
    $\widetilde{B}^{(2)}$ & $6639\times6639$ & $391499\approx3.06m$ & $1.80\times 10^1$ & $7.04\times 10^{6}$\\
    $\widetilde{B}^{(3)}$ & $1244\times1244$ & $417156\approx3.26m$ & $7.27\times 10^1$ & $3.03\times 10^{7}$\\
    $\widetilde{B}^{(4)}$ & $186\times186$ & $34596\approx0.27m$ & $4.47\times 10^1$ & $1.55\times 10^{6}$\\
    $\widetilde{A}^{(4)}$ & $33\times33$ & $1025\approx0.008m$ & $2.83\times 10^3$ & $2.90\times 10^{6}$\\
    total & - & $854088\approx 6.67m$ & - & $4.34\times 10^{7}$ \\
    \hline
    \end{tabular}
    \caption{Complexity results of 4-level MMD with localization using \Cref{alg:MMD_localized}.}
    \label{table:MMD_example1}
\end{table}
\begin{figure}[!h]
\centering
    \begin{subfigure}[b]{0.35\textwidth}
        \centering
        \includegraphics[width=1.0\textwidth]{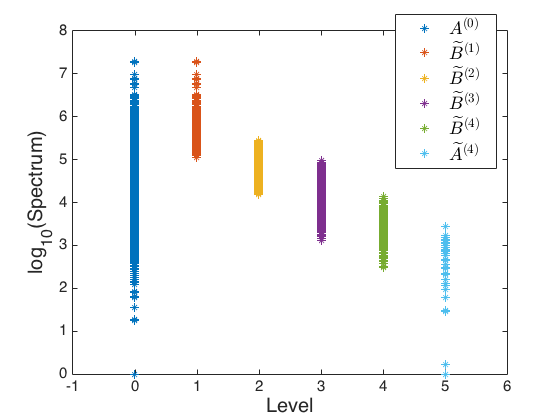}
        \caption{Spectrum}
    \end{subfigure}
    \begin{subfigure}[b]{0.35\textwidth}
        \centering
        \includegraphics[width=1.0\textwidth]{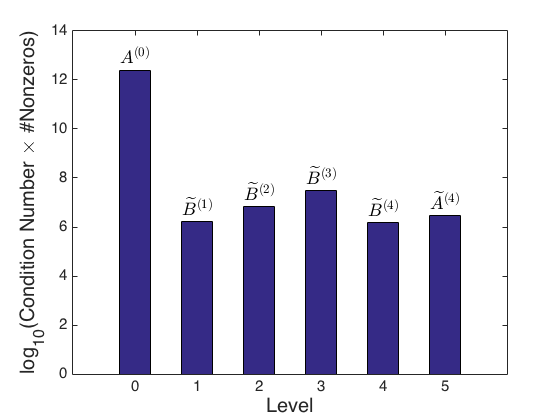}
        \caption{Complexity}
    \end{subfigure}
    \caption{Spectrum and complexity of each layer in 4-level MMD obtained from \Cref{alg:MMD_localized}.}
    \label{fig:MMD_case1}
\end{figure}

The Laplacian $L=A_0$ and the energy decomposition $\mathcal{E}$ are given as in \Cref{example:Example2}, and we apply a 4-level multiresolution matrix decomposition with localization using \Cref{alg:MMD_localized} to decompose the problem of solving $L^{-1}$. In this particular case we have $\lambda_{\max}(L)=1.93\times 10^7$ and $\lambda_{\min}(L)=1$. Again for graph laplacian, we choose $q=1$. On each level $k$, the partition is constructed subject to $\varepsilon(\mathcal{P}^{(k)},1)^2=10^{k-6}$ (i.e. $\{\varepsilon(\mathcal{P}^{(k)},1)^2\}_{k=1}^4=\{0.00001,0.0001,0.001,0.01\}$) and $\delta(\mathcal{P}^{(k)},1)\varepsilon(\mathcal{P}^{(k)},1)^2\leq50$. The compliment space $U^{(k)}$ are extended from $\Phi^{(k)}$ using a patch-wise QR factorization. So according to \Cref{cor:Bconditonnumber2}, each $\kappa(\widetilde{B}^{(k)}),\ k=2,3,4$ is expected to be bounded by $\delta(\mathcal{P}^{(k-1)},1)\varepsilon(\mathcal{P}^{(k)},1)^2 = \delta(\mathcal{P}^{(k-1)},1) \varepsilon(\mathcal{P}^{(k-1)},1)^2 \frac{\varepsilon(\mathcal{P}^{(k)},1)^2}{\varepsilon(\mathcal{P}^{(k-1)},1)^2} \leq 500$, $\kappa(\widetilde{B}^{(1)})$ is expected bounded by $\lambda_{\max}(L)\varepsilon(\mathcal{P}^{(0)},1)^2 = 1.93 \times 10^2$, and $\kappa(\widetilde{A}^{(4)})$ is expected to be bounded by $\delta(\mathcal{P}^{(3)},1)\lambda_{\min}(L)^{-1} = \delta(\mathcal{P}^{(3)},1) \varepsilon(\mathcal{P}^{(3)},1)^2 \frac{\lambda_{\min}(L)^{-1}}{\varepsilon(\mathcal{P}^{(3)},1)^2} \leq 5000$. Since we will use a conjugate gradient (CG) type method to compare the effectiveness of solving $L^{-1}$ directly and using the 4-level decomposition, the complexities of both approaches are proportional to the product of the number of nonzero (NNZ) entries and the condition number of the matrix concerned, given a fixed prescribed relative accuracy \cite{saad2011numerical}. Therefore we define the complexity of a matrix as the product of its NNZ entries and its condition number. Though here we use the sparsity of $B^{(k)} = (U^{(k)})^T A^{(k-1)} U^{(k)}$, in practice only the sparsity of $A^{(k-1)}$ matters and the matrix multiplication with respect to $U^{(k)}$ can be done by using the Householder vectors from the implicit QR factorization \cite{saad2011numerical}. The results not only satisfy the theoretical prediction, but also turn out to be much better than expected as shown in \Cref{table:MMD_example1} and \Cref{fig:MMD_case1}.
\begin{table}[!h]
\centering
\footnotesize
\renewcommand{\arraystretch}{1.5}
    \begin{tabular}{|c|c|c|c|c|}
    \hline
    Level & Size & $\#$Nonzero & Condition Number &Complexity\\ \hline
    $L$ & $10000\times10000$ & $128018\triangleq m$ & $1.93\times 10^7$ & $2.47\times 10^{12}$ \\ \hline
    $\widetilde{B}^{(1)}$ & $1898\times1898$ & $9812\approx0.07m$ & $1.72\times 10^2$ & $1.69\times 10^{6}$\\
    $\widetilde{B}^{(2)}$ & $6639\times6639$ & $391499\approx3.06m$ & $1.80\times 10^1$ & $7.04\times 10^{6}$\\
    $\widetilde{B}^{(3)}$ & $1014\times1014$ & $237212\approx1.85m$ & $2.56\times 10^1$ & $6.09\times 10^{6}$\\
    $\widetilde{B}^{(4)}$ & $313\times313$ & $86323\approx0.67m$ & $1.62\times 10^1$ & $1.40\times 10^{6}$\\
    $\widetilde{B}^{(5)}$ & $114\times114$ & $12996\approx0.10m$ & $5.54\times 10^1$ & $7.20\times 10^{5}$\\
    $\widetilde{A}^{(5)}$ & $22\times22$ & $442\approx0.004m$ & $1.60\times 10^3$ & $7.07\times 10^{5}$\\
    total & - & $738284\approx 5.77m$ & - & $1.76\times 10^{7}$ \\
    \hline
    \end{tabular}
    \caption{Complexity results of 5-level MMD with localization using \Cref{alg:MMD_localized}.}
    \label{table:MMD_example2}
\end{table}
\begin{figure}[!h]
\centering
    \begin{subfigure}[b]{0.35\textwidth}
        \centering
        \includegraphics[width=1.0\textwidth]{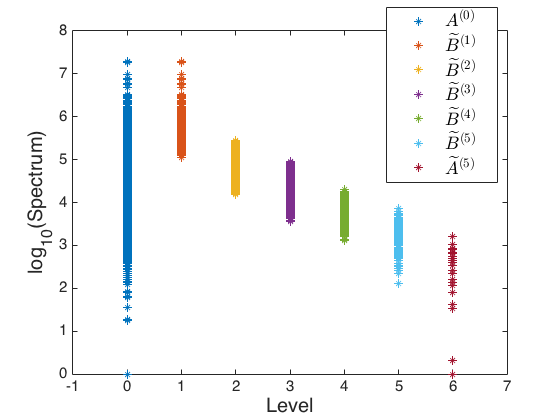}
        \caption{}
    \end{subfigure}
    \begin{subfigure}[b]{0.35\textwidth}
        \centering
        \includegraphics[width=1.0\textwidth]{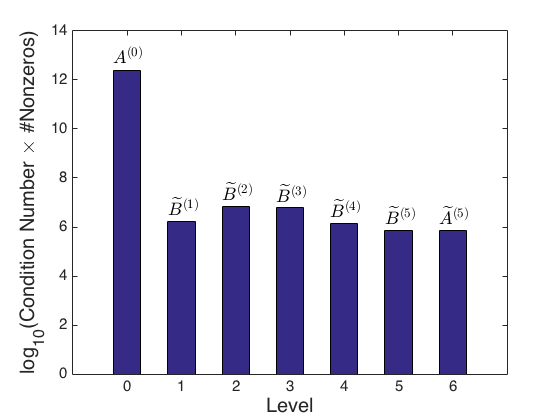}
        \caption{}
    \end{subfigure}
    \caption{Spectrum and complexity of each layer in 5-level MMD obtained from \Cref{alg:MMD_localized}.}
    \label{fig:MMD_case2}
\end{figure}
We now verify the performance of the 4-level decomposition by solving two particular systems $Lu^*=b$, where
\begin{equation*}
\begin{split}
\text{Case 1: }& \ u_i^*=(x_i^2+y_i^2+z_i^2)^\frac{1}{2},\ i=1,2,\cdots,n, b=Lu^*; \\
\text{Case 2: }& \ u_i^*=x_i+y_i+\sin(z_i),\ i=1,2,\cdots,n, b=Lu^*.
\end{split}
\end{equation*}
\noindent For both cases, we set the prescribed relative accuracy to be $\epsilon=10^{-5}$ such that $\|\hat{u}-u^*\|_L\leq \epsilon\|b\|_2$. By \Cref{theorem:MMD_levelerrorbound}, this accuracy can be achieved by imposing a corresponding accuracy control on each level's linear system relative error (i.e., $\mathrm{err}^{(K)}_A$ and $\mathrm{err}_B$ ) and localization error ($\mathrm{err}_{loc}$). In practice, instead of imposing a hard error control, we relax the localization error to $\varepsilon(\mathcal{P}^{(k)},1)$ (instead of $\epsilon)$ in order to ensure sparsity, which is actually how we obtain the 4-level decomposition with localization. Such relaxed localization error can be fixed by doing a compensation correction at level-0, which takes the output of \Cref{alg:MMD_solver} as an initialization to solve $Lu=b$. As shown in the gray columns ``$\#$ Iteration" in \Cref{table:MMD_case1} and \Cref{table:MMD_case2}, the number of iterations in the compensation calculation (which is the computation at Level 0) are less than 25 in both cases, which indicate that the localization error on each level is still small even when relaxed.

In particular, we use a preconditioned CG method to solve any involved linear systems $Au=b$ in \Cref{alg:MMD_solver}. The precondition matrix $D$ is chosen as the diagonal part of $A$ and we take $\bm{0}$ (all-zeros vector) as initials if no preconditioning vector is provided. The main computational cost of a single use of the PCG method is measured by the product of the number of iterations and the NNZ entries of the matrix involved (See the gray column ``Main Cost'' in both \Cref{table:MMD_case1} and \Cref{table:MMD_case2}).

In both cases, we can see that the total computational costs of our approach are obviously reduced compared to the direct use of PCG. Moreover, since the downward level-wise computation can be done in parallel, the effective computational costs of our approach are even less, which is the sum of the maximal cost among all levels and the cost of the compensation correction (See ``Parallel" row in \Cref{table:MMD_case1} and \Cref{table:MMD_case2} respectively).

Further, from the results we can see that the costs among all levels in both cases are mainly concentrated on level 2, namely, the inverting of $\widetilde{B}^{(3)}$. This observation implies that the structural/geometrical details of $L$ have more proportion on the scale corresponding to level 2 than on other scales, which is consistent to the fact that $\widetilde{B}^{(3)}$ on level 2 has the largest complexity of all. Though in practice we don't have the information in \Cref{table:MMD_case1} and \cref{table:MMD_case2}, we may observe the dominance of time complexity on level 2 after numbers of calls of our solver. As a natural improvement, we can simply further decompose the problem at level 2. More precisely, to relieve the dominance of level 2, we add one extra scale of $0.0003$ between the scales of $0.0001$ and $0.001$. Consequently, we obtain a similar 5-level decomposition with $\{\varepsilon(\mathcal{P}^{(k)},1)^2\}_{k=1}^5=\{0.00001,0.0001,0.0003,0.001,0.01\}$ (See \Cref{table:MMD_example2} and \Cref{fig:MMD_case2}). From the ``Main Cost 5-level'' columns of \Cref{table:MMD_case1} and \Cref{table:MMD_case2}), we can see that this simple improvement of further decomposition based on the feedback of computational results does reduce the computational cost.
\begin{table}[!h]
\centering
\footnotesize
{\tabulinesep=1.2mm \begin{tabu}{| c | c c | c c | c c | c c |}
\hline 
& \multicolumn{2}{c|}{Scale} & \multicolumn{2}{c|}{Matrix} & \multicolumn{2}{c|}{\# Iteration} & \multicolumn{2}{c|}{Main Cost} \\
\cline{2-9}
& \multicolumn{2}{c|}{$1 / \lambda_{\max}(L) \sim 1 / \lambda_{\min}(L)$} & \multicolumn{2}{c|}{$L$} & \multicolumn{2}{c|}{572} & \multicolumn{2}{c|}{$7.32 \times 10^7$} \\
\cline{1-9}
Level & 4-level & 5-level & 4-level & 5-level & 4-level & 5-level & 4-level & 5-level \\
\hline
0
&\cellcolor{Gray} $10^{-5} \sim 1 / \lambda_{\max}(L)$ & $10^{-5} \sim 1 / \lambda_{\max}(L)$
&\cellcolor{Gray} $\widetilde{B}^{(1)}$ & $\widetilde{B}^{(1)}$ 
&\cellcolor{Gray} 22 & 22
&\cellcolor{Gray} $2.15\times 10^5$ & $2.15\times 10^5$ \\
\hline
1
&\cellcolor{Gray} $10^{-5} \sim 10^{-4}$ & $10^{-5} \sim 10^{-4}$
&\cellcolor{Gray} $\widetilde{B}^{(2)}$ & $\widetilde{B}^{(2)}$
&\cellcolor{Gray} 11 & 11
&\cellcolor{Gray} $4.31\times 10^6$ & $4.31\times 10^6$ \\
\hline
\multirow{2}{*}{2}
&\cellcolor{Gray} $10^{-4} \sim 10^{-3}$ & $10^{-4} \sim 3\times 10^{-3}$
&\cellcolor{Gray} $\widetilde{B}^{(3)}$ & $\widetilde{B}^{(3)}$
&\cellcolor{Gray} 25 & 14
&\cellcolor{Gray} $1.04\times 10^7$ & $3.32\times 10^6$ \\

&\cellcolor{Gray}  & $3 \times 10^{-3} \sim 10^{-3}$
&\cellcolor{Gray}  & $\widetilde{B}^{(4)}$
&\cellcolor{Gray}  & 14
&\cellcolor{Gray}  & $1.21\times 10^6$ \\
\hline
3
&\cellcolor{Gray} $10^{-3} \sim 10^{-2}$ & $10^{-3} \sim 10^{-2}$
&\cellcolor{Gray} $\widetilde{B}^{(4)}$ & $\widetilde{B}^{(5)}$ 
&\cellcolor{Gray} 23 & 22
&\cellcolor{Gray} $7.96\times 10^5$ & $2.86\times 10^5$ \\
\hline
4
&\cellcolor{Gray} $10^{-2} \sim 1 / \lambda_{\min}(L)$ & $10^{-2} \sim 1 / \lambda_{\min}(L)$
&\cellcolor{Gray} $\widetilde{A}^{(4)}$ & $\widetilde{A}^{(5)}$ 
&\cellcolor{Gray} 31 & 22
&\cellcolor{Gray} $3.18\times 10^4$ & $9.72\times 10^3$ \\
\hline
0
&\cellcolor{Gray} - & -
&\cellcolor{Gray} L & L
&\cellcolor{Gray} 25 & 30
&\cellcolor{Gray} $3.20\times 10^6$ & $3.84\times 10^6$ \\
\hline
Total
&\cellcolor{Gray} - & -
&\cellcolor{Gray} - & -
&\cellcolor{Gray} 137 & 135
&\cellcolor{Gray} $1.90\times 10^7$ & $1.31\times 10^7$ \\
\hline
Parallel
&\cellcolor{Gray} - & -
&\cellcolor{Gray} - & -
&\cellcolor{Gray} - & -
&\cellcolor{Gray} $1.36\times 10^7$ & $8.15\times 10^6$ \\
\hline
\end{tabu}}
\caption{Case 1: Difference between 4-level and 5-level decomposition.}
\label{table:MMD_case1}
\end{table}
\begin{table}[!h]
\centering
\footnotesize
{\tabulinesep=1.2mm \begin{tabu}{| c | c c | c c | c c | c c |}
\hline 
& \multicolumn{2}{c|}{Scale} & \multicolumn{2}{c|}{Matrix} & \multicolumn{2}{c|}{\# Iteration} & \multicolumn{2}{c|}{Main Cost} \\
\cline{2-9}
& \multicolumn{2}{c|}{$1 / \lambda_{\max}(L) \sim 1 / \lambda_{\min}(L)$} & \multicolumn{2}{c|}{$L$} & \multicolumn{2}{c|}{586} & \multicolumn{2}{c|}{$7.50 \times 10^7$} \\
\cline{1-9}
Level & 4-level & 5-level & 4-level & 5-level & 4-level & 5-level & 4-level & 5-level \\
\hline
0
&\cellcolor{Gray} $10^{-5} \sim 1 / \lambda_{\max}(L)$ & $10^{-5} \sim 1 / \lambda_{\max}(L)$
&\cellcolor{Gray} $\widetilde{B}^{(1)}$ & $\widetilde{B}^{(1)}$ 
&\cellcolor{Gray} 24 & 24
&\cellcolor{Gray} $2.35\times 10^5$ & $2.35\times 10^5$ \\
\hline
1
&\cellcolor{Gray} $10^{-5} \sim 10^{-4}$ & $10^{-5} \sim 10^{-4}$
&\cellcolor{Gray} $\widetilde{B}^{(2)}$ & $\widetilde{B}^{(2)}$
&\cellcolor{Gray} 11 & 11
&\cellcolor{Gray} $4.31\times 10^6$ & $4.31\times 10^6$ \\
\hline
\multirow{2}{*}{2}
&\cellcolor{Gray} $10^{-4} \sim 10^{-3}$ & $10^{-4} \sim 3\times 10^{-3}$
&\cellcolor{Gray} $\widetilde{B}^{(3)}$ & $\widetilde{B}^{(3)}$
&\cellcolor{Gray} 25 & 14
&\cellcolor{Gray} $1.04\times 10^7$ & $3.32\times 10^6$ \\
&\cellcolor{Gray}  & $3 \times 10^{-3} \sim 10^{-3}$
&\cellcolor{Gray}  & $\widetilde{B}^{(4)}$
&\cellcolor{Gray}  & 14
&\cellcolor{Gray}  & $1.21\times 10^6$ \\
\hline
3
&\cellcolor{Gray} $10^{-3} \sim 10^{-2}$ & $10^{-3} \sim 10^{-2}$
&\cellcolor{Gray} $\widetilde{B}^{(4)}$ & $\widetilde{B}^{(5)}$ 
&\cellcolor{Gray} 23 & 23
&\cellcolor{Gray} $7.96\times 10^5$ & $2.99\times 10^5$ \\
\hline
4
&\cellcolor{Gray} $10^{-2} \sim 1 / \lambda_{\min}(L)$ & $10^{-2} \sim 1 / \lambda_{\min}(L)$
&\cellcolor{Gray} $\widetilde{A}^{(4)}$ & $\widetilde{A}^{(5)}$ 
&\cellcolor{Gray} 30 & 22
&\cellcolor{Gray} $3.08\times 10^4$ & $9.72\times 10^3$ \\
\hline
0
&\cellcolor{Gray} - & -
&\cellcolor{Gray} L & L
&\cellcolor{Gray} 16 & 23
&\cellcolor{Gray} $2.05\times 10^6$ & $2.94\times 10^6$ \\
\hline
Total
&\cellcolor{Gray} - & -
&\cellcolor{Gray} - & -
&\cellcolor{Gray} 129 & 131
&\cellcolor{Gray} $1.78\times 10^7$ & $1.23\times 10^7$ \\
\hline
Parallel
&\cellcolor{Gray} - & -
&\cellcolor{Gray} - & -
&\cellcolor{Gray} - & -
&\cellcolor{Gray} $1.24\times 10^7$ & $7.25\times 10^6$ \\
\hline
\end{tabu}}
\caption{Case 2: Difference between 4-level and 5-level decomposition.}
\label{table:MMD_case2}
\end{table}

\section{Conclusion and Future Works}
\label{sec:conclusion}

\subsection{Summary}
In this work, we introduce the notion of energy decomposition for symmetric positive definite matrix $A$ using the representation of energy elements. These energy elements help extracting the hidden geometric information of the operator, which serves the purpose of finding an appropriate partitioning of the basis. Specifically, we introduce the closed and interior energies which tightly bound the restricted energy/restricted matrix in terms of positive definiteness. These bounds further lead into the introduction of two important local measurements, the {\bf error factor} and the {\bf condition factor}, which can be calculated efficiently by solving a local and partial eigen problem. Using these local measurements, we propose a nearly-linear time algorithm to obtain an appropriate basis partitioning for compressing the operator with prescribed accuracy and bounded condition number. Extending the idea of operator compression into hierarchical formulation, we also propose a nearly-linear time solver for general symmetric positive definite matrix. The main idea is to decompose the operator into multiple scales of resolution such that the relative condition number in each scale can be bounded. Experimental results are reported to demonstrate the efficacy of our proposed algorithms.

\subsection{Future works}
This groundwork introduces the idea of {\bf energy decomposition} and its applications in operator compression and solving SPD linear system. We believe that the energy framework may prompt further research. Particularly, we discover further possible improvement of our algorithms during the development stage.


Due to the pairing characteristic of \Cref{alg:pair_clustering}, we are quite affirmative that our partitioning algorithm is not optimal. Instead, the clustering problem could be reformulated into some local optimization problem such that the construction of the partition $\mathcal{P}$ can be more robust. Secondly, our current implementation is a combination of Matlab and C++ coding and no parallel computing is included. Therefore, one of our future works is to develop an optimal coding such that more comparison experiments with state-of-the-art algorithms can be conducted.

In application-wise planning, we observe in \Cref{subsec:numerical_example2} that the compression scheme not only satisfies the prescribed accuracy and well-posedness requirement of the operator, but also preserves the geometric information of the spectrum. More specifically, we recall that the eigenvectors corresponding to the first few eigenvectors of the compressed operator $A_{\text{st}}$ give accurate approximation of the corresponding eigenvectors of the original inverse operator $A_{-1}$ (See \Cref{fig:eigenvectorplot}). Inspired by this observation, we place confidence in modifying our compression scheme into an efficient solver for partial eigen problem. Secondly, as demonstrated in the high contrast problem in \Cref{subsec:numerical_example2}, we believe that our energy decomposition framework can be specifically modified to suit the purpose of solving elliptic PDEs with high contrast coefficients. Based on our energy decomposition, more in-depth analysis and improvement could be made to show that such framework is one of the possible candidates to solve the elliptic type problem with highly varying coefficients. Thirdly, regarding the partitioning procedure and the locality of the basis in our algorithms, the localized MMD solver can be further improved to fit into the needs of frequent updating of the solver. For example, in graph laplacian system, our MMD solver can be updated dynamically if new vertices/new edges are added to the given graphs. This dynamic update greatly reduces the time for the regeneration of the solver, especially when the updates size is small.

\section*{Acknowledgment}
This research is supported in part by NSF Grant No. DMS-1318377 and DMS-1613861. We would also like to thank Prof. Houman Owhadi for the suggestions and discussions throughout the development of this work.

\section*{Appendix: Proofs of Theorems}
\begin{proof}[Proof of \Cref{lemma:phi_error}]
$ $\\
\vspace{-3mm}
\begin{enumerate}
\item Let $y=A^{-1}(P_\Phi b)\in \Psi$, then 
\begin{align*}
\|x-y\|_A^2=&\ (x-y,A(x-y))= (x-y-P_\Phi(x-y),b- P_\Phi b)\\
& \leq \|x-y-P_\Phi(x-y)\|_2\|b-P_\Phi b\|_2\leq \epsilon\|x-y\|_A\|b\|_2,
\end{align*}
and thus we have
\[\|x-P_\Psi^A x\|_A\leq \|x-y\|_A\leq \epsilon \|b\|_2. \]
\item Let $z=A^{-1}(x-P_\Psi^A x)$, then 
\begin{align*}
\|x-P_\Psi^A x\|_2^2=&\ (x- P_\Psi^A x,x-P_\Psi^A x) = (x- P_\Psi^A x,A z) \\
& = (x- P_\Psi^A x,z-P_\Psi^A z)_A \leq \|x- P_\Psi^A x\|_A\|z-P_\Psi^A z\|_A \leq \epsilon\|x- P_\Psi^A x\|_A\|Az\|_2,
\end{align*}
and thus
\[\|x-P_\Psi^A x\|_2\leq \epsilon\|x- P_\Psi^A x\|_A\leq \epsilon^2\|b\|_2.\]
\item Immediate result of 2.
\end{enumerate}
\end{proof}

\begin{proof}[Proof of \Cref{lemma:int_spectrum}]
$ $\\
\vspace{-3mm}

Let $\Phi^k\subset \Span\{\mathcal{S}\}$ denote the eigenspace of $\underline{A}_\mathcal{S}$(as a operator restricted to $\Span\{\mathcal{S}\}$) corresponding to interior eigenvalues $\lambda_1(\mathcal{S})\leq \lambda_2(\mathcal{S})\leq \cdots\leq\lambda_k(\mathcal{S})$. On the one hand, for all $x\in \Span\{\mathcal{S}\}$, we have
\begin{align*}
\|x\|_{\underline{A}_\mathcal{S}}^2\geq\|x-P_{\Phi^{q(\epsilon)}}x\|_{\underline{A}_\mathcal{S}}^2 
&= (x-P_{\Phi^{q(\epsilon)}}x)^T\underline{A}_\mathcal{S}(x-P_{\Phi^{q(\epsilon)}}x) \geq \lambda_{q(\epsilon)+1}(\mathcal{S})\|x-P_{\Phi^{q(\epsilon)}}x\|_2^2, \\
\Longrightarrow \|x-P_{\Phi^{q(\epsilon)}}x\|_2 &\leq \frac{1}{\sqrt{\lambda_{q(\epsilon)+1}(\mathcal{S})}}\|x\|_{\underline{A}_\mathcal{S}} \leq \epsilon \|x\|_{\underline{A}_\mathcal{S}}.
\end{align*} 
Thus $\Phi^{q(\epsilon)}\in \mathcal{G}(\epsilon)$, $q(\epsilon)\geq p(\epsilon)$. On the other hand, assume that the minimum $p(\epsilon)$ is achieved by some space $\widetilde{\Theta}$, then one can check that 
\begin{align*}\lambda_{p(\epsilon)+1}=&\ \max_{\begin{subarray}{c}\Theta \subset \Span\{\mathcal{S}\}\\ \dim\Theta= p(\epsilon)\end{subarray}}\min_{x\in \Span\{\mathcal{S}\}}\frac{\|x-P_\Theta x\|_{\underline{A}_\mathcal{S}}^2}{\|x-P_\Theta x\|^2_2}\\ 
\geq&\ \min_{x\in \Span\{\mathcal{S}\}}\frac{\|x-P_{\widetilde{\Theta}} x\|_{\underline{A}_\mathcal{S}}^2}{\|x-P_{\widetilde{\Theta}} x\|^2_2}\ \\
=&\ \min_{x\in \mathrm{span}\{\mathcal{S}\}}\frac{\|x-P_{\widetilde{\Theta}} x\|_{\underline{A}_\mathcal{S}}^2}{\|x-P_{\widetilde{\Theta}} x-P_{\widetilde{\Theta}}(x-P_{\widetilde{\Theta}} x)\|^2_2}\qquad \Big(\text{since}\ P_{\widetilde{\Theta}}(x-P_{\widetilde{\Theta}} x)=0\Big)\\
=&\ \min_{\begin{subarray}{c} y=x-P_{\widetilde{\Theta}} x\\ x\in \mathrm{span}\{\mathcal{S}\}\end{subarray}}\frac{\|y\|_{\underline{A}_\mathcal{S}}^2}{\|y-P_{\widetilde{\Theta}}y\|^2_2}\\
\geq&\ \min_{ y\in \mathrm{span}\{\mathcal{S}\}}\frac{\|y\|_{\underline{A}_\mathcal{S}}^2}{\|y-P_{\widetilde{\Theta}}y\|^2_2} \geq\ \frac{1}{\epsilon^2},
\end{align*}
which implies $p(\epsilon)\geq q(\epsilon)$ by the definition of $q(\epsilon)$. Finally we have $p(\epsilon)= q(\epsilon)$.
\end{proof}

\begin{proof}[Proof of \Cref{lemma:localization_psi}]
$ $\\
\vspace{-3mm}

We only need to prove property 1, properties 2 and 3 follow by using the same argument as in \Cref{lemma:phi_error}. Recall that we have
\[\|x-P_\Psi^A x\|_A=\|x-\Psi c\|_A\leq\epsilon\|b\|_2,\]
with $c=A_{\text{st}}^{-1}\Psi^TAx$. Let $y_1=\Psi c=\sum_{i=1}^N c_i\psi_i$, $y_2=\widetilde{\Psi}c=\sum_{i=1}^N c_i\tilde{\psi}_i$. Then we have
\[\|y_1-y_2\|_A=\left\Vert\sum_{i=1}^N c_i(\psi_i-\tilde{\psi}_i)\right\Vert_A\leq \sum_{i=1}^N |c_i|\|\psi_i-\tilde{\psi}_i\|_A\leq \frac{C\epsilon \sqrt{N}}{\sqrt{N}}\big(\sum_{i=1}^N c_i^2 \big)^\frac{1}{2}=C\epsilon \sqrt{c^Tc}. \] 
Notice that 
\[c^Tc=x^TA\Psi A_{\text{st}}^{-2}\Psi^TAx\leq \|A^\frac{1}{2}\Psi A_{\text{st}}^{-2}\Psi^TA^\frac{1}{2}\|_2\|x\|_A^2, \]
\[\|A^\frac{1}{2}\Psi A_{\text{st}}^{-2}\Psi^TA^\frac{1}{2}\|_2=\|A_{\text{st}}^{-1}\Psi^TA\Psi A_{\text{st}}^{-1}\|_2=\|A_{\text{st}}^{-1}\|_2,\]
\[\|x\|_A^2=b^TA^{-1}b\leq \|A^{-1}\|_2\|b\|_2^2,\]
therefore we get
\[\|y_1-y_2\|_A\leq C\epsilon\sqrt{\|A_{\text{st}}^{-1}\|_2\|A^{-1}\|_2} \|b\|_2\leq  C\epsilon\|A^{-1}\|_2\|b\|_2.\]
Then we have 
\[\|x-y_2\|_A\leq\|x-y_1\|_A+\|y_1-y_2\|_A\leq\epsilon\|b\|_2+C\epsilon\|A^{-1}\|_2\|b\|_2=(1+C\|A^{-1}\|_2)\epsilon\|b\|_2.\]
Since $y_2\in \Span\{\widetilde{\Psi}\}$, we obtain
\[\|x-P_{\widetilde{\Psi}}^A x\|_A\leq\|x-y_2\|_A\leq (1+C\|A^{-1}\|_2)\epsilon\|b\|_2.\]
\end{proof}

\begin{proof}[Proof of \Cref{thm:psi_delta}]
$ $\\
\vspace{-3mm} 

$\|\psi_i\|_A\leq \|\psi_i^k\|_A\leq\|\psi_i^0\|_A$ has been proved in the construction of local approximators. We only need to prove $\|\psi_i^0\|_A\leq \sqrt{\delta(P_{j_i})}$. Recall that $\psi_i^0$ is defined as 
\begin{equation*}
\begin{array}{rcl}
 \psi_i^0 &=& \underset{{x\in \Span\{P_{j_i}\}}}{\argmin} \|x\|_A,\\
\text{subject to } \quad \varphi_{i'}^Tx &=& \delta_{i',i},\quad \forall\ i' = 1,\ldots, N.
\end{array}
\end{equation*}

Without loss of generality, we can assume that $\varphi_i$ is the first column of $\Phi_{j_i}$. And notice that $\|x\|_A=\|x\|_{A_{P_{j_i}}}$, therefore the optimization formation can be rewritten as 
\begin{equation}
\begin{array}{rcl}
 \psi_i^0 &=& \underset{{x\in \Span\{P_{j_i}\}}}{\argmin} \|x\|_{A_{P_{j_i}}},\\
\text{subject to } \quad \Phi_{j_i}^Tx &=& z_i,
\end{array}
\label{eqt:optimization_initial}
\end{equation}
where $z_i=(1,0,0,\cdots,0)^T\in \mathbb{R}^{q_{j_i}}$. This optimization problem can be uniquely and explicitly solved as 
\begin{equation}
\psi_i^0=A_{P_{j_i}}^{-1}\Phi_{j_i}(\Phi_{j_i}^TA_{P_{j_i}}^{-1}\Phi_{j_i})^{-1}z_i,
\end{equation}
where again $A_{P_{j_i}}^{-1}$ denotes the inverse of $A_{P_{j_i}}$ as an operator restricted to $\Span\{P_{j_i}\}$. And thus we have
\begin{equation}
\|\psi_i^0\|_A^2=\|\psi_i^0\|_{A_{P_{j_i}}}^2=z_i^T(\Phi_{j_i}^TA_{P_{j_i}}^{-1}\Phi_{j_i})^{-1}z_i.
\end{equation}
Notice that 
\begin{equation*}
\overline{A}_{P_{j_i}}\succeq A_{P_{j_i}} \Rightarrow\ A_{P_{j_i}}^{-1}\succeq \overline{A}_{P_{j_i}}^{-1} \Rightarrow\ \Phi_{j_i}^TA_{P_{j_i}}^{-1}\Phi_{j_i}\succeq \Phi_{j_i}^T\overline{A}_{P_{j_i}}^{-1}\Phi_{j_i}\ \Rightarrow\ (\Phi_{j_i}^T\overline{A}_{P_{j_i}}^{-1}\Phi_{j_i})^{-1}\succeq (\Phi_{j_i}^TA_{P_{j_i}}^{-1}\Phi_{j_i})^{-1},
\end{equation*}
since $\overline{A}_{P_{j_i}}$ and $A_{P_{j_i}}$ are both symmetric, positive definite as operators restricted to $\Span\{P_{j_i}\}$. Therefore by \Cref{eqt:delta&Ast2} we have
\begin{equation}
\|\psi_i^0\|_A^2\leq z_i^T(\Phi_{j_i}^T\overline{A}_{P_{j_i}}^{-1}\Phi_{j_i})^{-1}z_i\leq \|(\Phi_{j_i}^T\overline{A}_{P_{j_i}}^{-1}\Phi_{j_i})^{-1}\|_2\|z_i\|_2^2=\delta(P_{j_i}).
\end{equation}
\end{proof}

\begin{proof}[Proof of \Cref{cor:expdecay}]
$ $\\
\vspace{-3mm}

This proof basically follows the idea in \cite{owhadi2017multigrid}. For any $\psi_i$, recall that $\psi_i^k\in \Span(S_k(P_{j_i}))$, thus $\underline{A}_{S_k^c(P_{j_i})}\psi_i^k=\bm{0}$, and
\[\|\psi_i\|_{\underline{A}_{S_k^c(P_{j_i})}}^2=\|\psi_i-\psi_i^k\|_{\underline{A}_{S_k^c(P_{j_i})}}^2\leq \|\psi_i-\psi_i^k\|_A^2\leq \Big(\frac{\alpha(\mathcal{P})-1}{\alpha(\mathcal{P})}\Big)^k\delta(P_{j_i}).\]
For any two $\psi_i,\psi_{i'}$, since $(\psi_i-\psi_i^k)^T\Phi=\bm{0}$, and $A\psi_{i'}\in \Phi$, we have $(\psi_i-\psi_i^k)^TA\psi_{i'}=0$ for all $k$. Also notice that $(\psi_i^{k})^TA\psi_{i'}^k=0$ for $k<\frac{k_{ii'}}{2}$, since $S_k(P_{j_i})\cap S_k(P_{j_{i'}})=\emptyset$ when $2k<k_{ii'}$. Therefore taking $k=\lceil\frac{k_{ii'}}{2}\rceil-1$ we have
\begin{align*}
|\psi_i^TA\psi_{i'}|=|(\psi_i^k)^TA(\psi_{i'}-\psi_{i'}^k)|\leq&\ \|\psi_i^k\|_A\|\psi_{i'}-\psi_{i'}^k\|_A\\
\leq&\ \Big(\frac{\alpha(\mathcal{P})-1}{\alpha(\mathcal{P})}\Big)^\frac{k}{2}\delta(\mathcal{P})\leq \Big(\frac{\alpha(\mathcal{P})-1}{\alpha(\mathcal{P})}\Big)^{\frac{k_{ii'}}{4}-\frac{1}{2}}\delta(\mathcal{P}).
\end{align*}
\end{proof}

\begin{proof}[Proof of \Cref{thm:condition_number}]
$ $\\
\vspace{-3mm}

Thanks to \Cref{eqt:stiffness}, and since $\Phi^T\Phi=I_M$, we have
\begin{equation}
\|A_{\text{st}}^{-1}\|_2=\|\Phi^TA^{-1}\Phi\|_2\leq \|A^{-1}\|_2\quad \Longrightarrow\quad \lambda_{\min}(A_{\text{st}})\geq \lambda_{\min}(A).
\end{equation}
Thanks to \Cref{eqt:delta&Ast2}, and since $A\preceq \sum_{j=1}^M \overline{A}_{P_j}$, we have
\begin{equation}
\delta(\mathcal{P})I_M \succeq (\Phi^T(\sum_{j=1}^M \overline{A}_{P_j})^{-1}\Phi)^{-1} \succeq (\Phi^TA^{-1}\Phi)^{-1}=A_{\text{st}} \Longrightarrow \lambda_{\max}(A_{\text{st}}) \leq \delta(\mathcal{P}).
\end{equation}
\end{proof}

\begin{proof}[Proof of \Cref{cor:localpsi_conditionnum}] 
$ $\\
\vspace{-3mm}

Notice that $P_{\widetilde{\Psi}}^A= \widetilde{\Psi} \widetilde{A}_{\text{st}}^{-1}\widetilde{\Psi}^TA$ is a projection with respect to the energy inner product, we have
\[A\widetilde{\Psi}\widetilde{A}_{\text{st}}^{-1}\widetilde{\Psi}^TA\preceq A\quad \Longrightarrow\quad \widetilde{\Psi}\widetilde{A}_{\text{st}}^{-1}\widetilde{\Psi}^T\preceq A^{-1},\]
and since $\Phi^T\widetilde{\Psi}=I_{N}$, we have
\[\|\widetilde{A}_{\text{st}}^{-1}\|_2\leq \|\Phi^TA^{-1}\Phi\|_2\leq \|A^{-1}\|_2\quad \Longrightarrow\quad \lambda_{\min}(\widetilde{A}_{\text{st}})\geq \lambda_{\min}(A).\]
For any $c\in \mathbb{R}^N$, using a similar argument in \Cref{lemma:localization_psi}, we have 
\[\|\Psi c-\widetilde{\Psi} c\|_A\leq C\epsilon\|c\|_2,\]
then we get
\[c^T\widetilde{A}_{\text{st}}c=\|\widetilde{\Psi} c\|_A^2\leq \big(\|\Psi c-\widetilde{\Psi} c\|_A+\|\Psi c\|_A\big)^2\leq(C\epsilon+\sqrt{\lambda_{\max}(A_{\text{st}})})^2\|c\|_2^2,\]
and using $\lambda_{\max}(A_{\text{st}})\leq \delta(\mathcal{P})$, we have
\[\lambda_{\max}(\widetilde{A}_{\text{st}})\leq \big(1+\frac{C\epsilon}{\sqrt{\delta(\mathcal{P})}}\big)^2\delta(\mathcal{P}).\]
\end{proof}

\begin{proof}[Proof of \Cref{theorem:MMD_levelerrorbound}]
$ $\\
\vspace{-3mm}

First by assumption \ref{itm:assumption2}, we have 
\[\|(\widetilde{A}^{(K)})^{-1}b-\mathrm{Inv}(\widetilde{A}^{(k)})b\|_{\widetilde{A}^{(K)}}\leq \mathrm{err}^{(K)}_A\|b\|_{(\widetilde{A}^{(K)})^{-1}},\quad \forall b\in \mathbb{R}^{N^{(K)}}.\]
To perform induction, we assume that at level $k$, $(\widetilde{A}^{(k)})^{-1}$ can be solved subject to a relative error $\mathrm{err}^{(k)}_A$ in the sense that the solver $\mathrm{Inv}(\widetilde{A}^{(k)})$ satisfies 
\[\|(\widetilde{A}^{(k)})^{-1}b-\mathrm{Inv}(\widetilde{A}^{(k)})b\|_{\widetilde{A}^{(k)}}\leq \mathrm{err}^{(k)}_A\|b\|_{(\widetilde{A}^{(k)})^{-1}},\quad \forall b\in \mathbb{R}^{N^{(k)}}.\]
Recall that $(\widetilde{A}^{(k-1)})^{-1}$ and the solver $\mathrm{Inv}(\widetilde{A}^{(k-1)})$ are given by 
\begin{equation}
(\widetilde{A}^{(k-1)})^{-1}=U^{(k)}\big((U^{(k)})^T\widetilde{A}^{(k-1)}U^{(k)}\big)^{-1}(U^{(k)})^T + \Psi^{(k)}\big((\Psi^{(k)})^T\widetilde{A}^{(k-1)}\Psi^{(k)}\big)^{-1}(\Psi^{(k)})^T,
\end{equation}
\begin{equation}
\mathrm{Inv}(\widetilde{A}^{(k-1)})=U^{(k)}\mathrm{Inv}(\widetilde{B}^{(k)})(U^{(k)})^T + \widetilde{\Psi}^{(k)}\mathrm{Inv}(\widetilde{A}^{(k)})(\widetilde{\Psi}^{(k)})^T,
\end{equation}
then for any $b\in \mathbb{R}^{N^{(k-1)}}$, we have
\begin{align*}
&\ \|(\widetilde{A}^{(k-1)})^{-1}b-\mathrm{Inv}(\widetilde{A}^{(k)})b\|_{\widetilde{A}^{(k-1)}}\\
\leq &\ \|U^{(k)}\big((U^{(k)})^T\widetilde{A}^{(k-1)}U^{(k)}\big)^{-1}(U^{(k)})^Tb-U^{(k)}\mathrm{Inv}(\widetilde{B}^{(k)})(U^{(k)})^Tb\|_{\widetilde{A}^{(k-1)}}\\
&\ + \|\Psi^{(k)}\big((\Psi^{(k)})^T\widetilde{A}^{(k-1)}\Psi^{(k)}\big)^{-1}(\Psi^{(k)})^Tb-\widetilde{\Psi}^{(k)}\big((\widetilde{\Psi}^{(k)})^T\widetilde{A}^{(k-1)}\widetilde{\Psi}^{(k)}\big)^{-1}(\widetilde{\Psi}^{(k)})^Tb\|_{\widetilde{A}^{(k-1)}}\\
&\ + \|\widetilde{\Psi}^{(k)}\big((\widetilde{\Psi}^{(k)})^T\widetilde{A}^{(k-1)}\widetilde{\Psi}^{(k)}\big)^{-1}(\widetilde{\Psi}^{(k)})^Tb-\widetilde{\Psi}^{(k)}\mathrm{Inv}(\widetilde{A}^{(k)})(\widetilde{\Psi}^{(k)})^Tb\|_{\widetilde{A}^{(k-1)}}\\
=&\ I_1+I_2+I_3.
\end{align*}
Recall that $\widetilde{A}^{(k)}=(\widetilde{\Psi}^{(k)})^T\widetilde{A}^{(k-1)}\widetilde{\Psi}^{(k)}$, $\widetilde{B}^{(k)}=(U^{(k)})^T\widetilde{A}^{(k-1)}U^{(k)}$, then by assumption \ref{itm:assumption1} we have
\begin{align*}
I_1=&\ \|(\widetilde{B}^{(k)})^{-1}(U^{(k)})^Tb-\mathrm{Inv}(\widetilde{B}^{(k)})(U^{(k)})^Tb\|_{\widetilde{B}^{(k)}}\\
\leq&\ \mathrm{err}_B \Big(b^TU^{(k)}(\widetilde{B}^{(k)})^{-1}(U^{(k)})^Tb\Big)^\frac{1}{2}\\
\leq&\ \mathrm{err}_B \|(\widetilde{A}^{(k-1)})^\frac{1}{2}U^{(k)}((U^{(k)})^T\widetilde{A}^{(k-1)}U^{(k)})^{-1}(U^{(k)})^T(\widetilde{A}^{(k-1)})^\frac{1}{2}\|_2\|b\|_{(\widetilde{A}^{(k-1)})^{-1}}\\
\leq&\ \mathrm{err}_B\|b\|_{(\widetilde{A}^{(k-1)})^{-1}}.
\end{align*}
Similarly by the assumption of induction, we have
\begin{equation*}
I_3 = \|(\widetilde{A}^{(k)})^{-1}(\widetilde{\Psi}^{(k)})^Tb-\mathrm{Inv}(\widetilde{A}^{(k)})(\widetilde{\Psi}^{(k)})^Tb\|_{\widetilde{A}^{(k)}} \leq \mathrm{err}^{(k)}_A\|b\|_{(\widetilde{A}^{(k-1)})^{-1}}.
\end{equation*}
Let $x=(\widetilde{A}^{(k-1)})^{-1}b$, then we get 
\begin{equation*}
I_2=\|P_{\Psi^{(k)}}^{\widetilde{A}^{(k-1)}}x-P_{\widetilde{\Psi}^{(k)}}^{\widetilde{A}^{(k-1)}}x\|_{\widetilde{A}^{(k-1)}}\leq \|P_{\Psi^{(k)}}^{\widetilde{A}^{(k-1)}}x-P_{\widetilde{\Psi}^{(k)}}^{\widetilde{A}^{(k-1)}}P_{\Psi^{(k)}}^{\widetilde{A}^{(k-1)}}x\|_{\widetilde{A}^{(k-1)}}+\|P_{\widetilde{\Psi}^{(k)}}^{\widetilde{A}^{(k-1)}}P_{U^{(k)}}^{\widetilde{A}^{(k-1)}}x\|_{\widetilde{A}^{(k-1)}}.
\end{equation*}
Using a similar argument in \Cref{lemma:localization_psi} we can actually prove by assumption \ref{itm:assumption3} that 
\[\|P_{\Psi^{(k)}}^{\widetilde{A}^{(k-1)}}x-P_{\widetilde{\Psi}^{(k)}}^{\widetilde{A}^{(k-1)}}P_{\Psi^{(k)}}^{\widetilde{A}^{(k-1)}}x\|_{\widetilde{A}^{(k-1)}}\leq \frac{1}{2}\mathrm{err}_{loc}\|P_{\Psi^{(k)}}^{\widetilde{A}^{(k-1)}}x\|_{\widetilde{A}^{(k-1)}} \leq \frac{1}{2}\mathrm{err}_{loc}\|b\|_{(\widetilde{A}^{(k-1)})^{-1}},\]
and 
\begin{align*}
\|P_{\widetilde{\Psi}^{(k)}}^{\widetilde{A}^{(k-1)}}P_{U^{(k)}}^{\widetilde{A}^{(k-1)}}x\|_{\widetilde{A}^{(k-1)}}\leq&\  \|((\widetilde{\Psi}^{(k)})^T\widetilde{A}^{(k-1)}\widetilde{\Psi}^{(k)})^{-1}\|_2^\frac{1}{2}\|(\widetilde{\Psi}^{(k)})^T\widetilde{A}^{(k-1)}P_{U^{(k)}}^{\widetilde{A}^{(k-1)}}x\|_2\\
\leq &\  \|A^{-1}\|_2^\frac{1}{2}\|(\widetilde{\Psi}^{(k)}-\Psi^{(k)})^T\widetilde{A}^{(k-1)}P_{U^{(k)}}^{\widetilde{A}^{(k-1)}}x\|_2\\
\leq &\ \|A^{-1}\|_2^\frac{1}{2}\|(\widetilde{\Psi}^{(k)}-\Psi^{(k)})^T\widetilde{A}^{(k-1)}(\widetilde{\Psi}^{(k)}-\Psi^{(k)})\|_2^\frac{1}{2}\|P_{U^{(k)}}^{\widetilde{A}^{(k-1)}}x\|_{\widetilde{A}^{(k-1)}}\\
\leq &\ \frac{1}{2}\mathrm{err}_{loc}\|b\|_{(\widetilde{A}^{(k-1)})^{-1}},
\end{align*}
thus $I_2\leq \mathrm{err}_{loc}\|b\|_{(\widetilde{A}^{(k-1)})^{-1}}$. Finally we have
\[\|(\widetilde{A}^{(k-1)})^{-1}b-\mathrm{Inv}(\widetilde{A}^{(k-1)})b\|_{\widetilde{A}^{(k-1)}}\leq(\mathrm{err}_B+\mathrm{err}_{loc}+\mathrm{err}^{(k)}_A)\|b\|_{(\widetilde{A}^{(k-1)})^{-1}},\]
that is we have
\[\mathrm{err}^{(k-1)}_A=(\mathrm{err}_B+\mathrm{err}_{loc}+\mathrm{err}^{(k)}_A).\]
Then by induction the relative total error using $K$-level's decomposition with localization for solving $A^{-1}$ is 
\[\mathrm{err}_{total}=\mathrm{err}^{(0)}_A=K(\mathrm{err}_B+\mathrm{err}_{loc})+\mathrm{err}^{(K)}_A.\]
\end{proof}

\begin{proof}[Proof of \Cref{thm:error_accumulate}]
$ $\\
\vspace{-3mm}

Again by \Cref{lemma:phi_error}, we only need to prove \Cref{eqt:error_accumulate}. For consistency, we write $\bm{\Phi}^{(0)}=I_n$, and correspondingly $\bm{\Psi}^{(0)}=I_n$, $P_{\bm{\Phi}^{(0)}}=I_n$, $P_{\bm{\Psi}^{(0)}}^{A}=I_n$. Using \Cref{eqt:mutliresolution_phi&U&psi}, it is easy to check that for any $x\in \mathbb{R}^n$ and any $k_1\leq k_2\leq k_3$, 
\[(P_{\bm{\Phi}^{(k_1)}}x-P_{\bm{\Phi}^{(k_2)}}x)^T(P_{\bm{\Phi}^{(k_2)}}x-P_{\bm{\Phi}^{(k_3)}}x)=0,\]
thus we have
\[\|x-P_{\bm{\Phi}^{(K)}}x\|_2^2=\|\sum_{k=1}^K(P_{\bm{\Phi}^{(k-1)}}x-P_{\bm{\Phi}^{(k)}}x)\|_2^2=\sum_{k=1}^K\|P_{\bm{\Phi}^{(k-1)}}x-P_{\bm{\Phi}^{(k)}}x\|_2^2.\]
Notice that 
\[P_{\bm{\Phi}^{(k-1)}}x-P_{\bm{\Phi}^{(k)}}x=\bm{\Phi}^{(k-1)}(\bm{\Phi}^{(k-1)})^Tx-\bm{\Phi}^{(k-1)}\Phi^{(k)}(\Phi^{(k)})^T(\bm{\Phi}^{(k-1)})^Tx,\]
thus by the construction of $\Phi^{(k)}$( or $\bm{\Phi}^{(k)}$), we have
\begin{align*}
\|P_{\bm{\Phi}^{(k-1)}}x-P_{\bm{\Phi}^{(k)}}x\|_2^2=&\ \|(\bm{\Phi}^{(k-1)})^Tx-\Phi^{(k)}(\Phi^{(k)})^T(\bm{\Phi}^{(k-1)})^Tx\|_2^2\\
\leq&\ \varepsilon(\mathcal{P}^{(k)},q^{(k)})^2\|(\bm{\Phi}^{(k-1)})^Tx\|_{A^{(k-1)}}^2 \\
=&\ \varepsilon(\mathcal{P}^{(k)},q^{(k)})^2 x^T\bm{\Phi}^{(k-1)}\big((\bm{\Phi}^{(k-1)})^TA^{-1}\bm{\Phi}^{(k-1)}\big)^{-1}(\bm{\Phi}^{(k-1)})^Tx\\
=&\ \varepsilon(\mathcal{P}^{(k)},q^{(k)})^2\|P_{\bm{\Psi}^{(k-1)}}^Ax\|_A^2\\
\leq&\ \varepsilon(\mathcal{P}^{(k)},q^{(k)})^2\|x\|_A^2.
\end{align*}
We have used the fact that 
\begin{align*}
\|P_{\bm{\Psi}^{(k)}}^Ax\|_A^2=&\ x^TA\bm{\Psi}^{(k)}\big((\bm{\Psi}^{(k)})^TA\bm{\Psi}^{(k)}\big)^{-1}(\bm{\Psi}^{(k)})^TAx\\
=&\ x^T\bm{\Phi}^{(k)}\big((\bm{\Phi}^{(k)})^TA^{-1}\bm{\Phi}^{(k)}\big)^{-1}(\bm{\Phi}^{(k)})^Tx, \quad \forall k\geq0.
\end{align*}
Therefore we have
\[\|x-P_{\bm{\Phi}^{(K)}}x\|_2^2\leq \Big(\sum_{k=0}^{K}\varepsilon(\mathcal{P}^{(k)},q^{(k)})^2\Big)\|x\|_A^2. \]
\end{proof}

\bibliographystyle{siamplain}
\bibliography{reference}
\end{document}

%% file: shared_info.tex
\usepackage{amsopn}
\usepackage{algorithmic}
\Crefname{ALC@unique}{Line}{Lines}
\usepackage{array}
\usepackage{amsfonts}
\usepackage{amsmath}
\usepackage{amssymb}
\usepackage{blkarray}
\usepackage{bm}
\usepackage{bbm}
\usepackage{caption}
\usepackage{color, colortbl}
\usepackage{enumitem}
\definecolor{Gray}{gray}{0.9}
\definecolor{Gray2}{gray}{0.95}
\definecolor{Gray3}{gray}{0.7}
\usepackage{epstopdf}
\usepackage{graphicx}
\usepackage{lipsum}
\usepackage{longtable}
\ifpdf
  \DeclareGraphicsExtensions{.eps,.pdf,.png,.jpg}
\else
  \DeclareGraphicsExtensions{.eps}
\fi
\usepackage{mathtools}
\usepackage{multirow}
\usepackage{ulem}
\usepackage{setspace}
\usepackage{subcaption}
\usepackage{tabu}
\usepackage{xr-hyper}

\newtheorem{thm}{Theorem}[section]

\newtheorem{corollary}[thm]{Corollary}
\newtheorem{proposition}[thm]{Proposition}
\newtheorem{lemma}[thm]{Lemma}
\newtheorem{definition}[thm]{Definition}
\newtheorem{remark}[thm]{Remark}
\newtheorem{example}[thm]{Example}
\newtheorem{construction}[theorem]{Construction}

\newcommand{\Span}{\mathrm{span}}
\newcommand{\argmin}{\arg\min}

\newcommand{\TheTitle}{An adaptive fast solver for a general class of positive definite matrices via energy decomposition} 
\newcommand{\TheAuthors}{Thomas Y. Hou, D. Huang, K.C. Lam, P. Zhang}



\usepackage{fancyhdr}

\title{{\TheTitle}}

\author{
  Thomas Y. Hou
  \and
  De Huang
  \and
  Ka Chun Lam
  \and
  PengChuan Zhang
}

\usepackage{amsopn}